\documentclass[11pt,reqno,oneside]{amsart}
\usepackage{graphicx, subfigure,}
\usepackage{amssymb}
\usepackage{epstopdf}
\usepackage{amssymb}
\usepackage[a4paper, total={6.6in, 9.9in}]{geometry}
\usepackage{multirow}
\usepackage{bm}
\usepackage{amsmath,amsfonts,amsthm,mathrsfs,amssymb,cite}
\usepackage{mathtools}
\usepackage{framed}
\usepackage{exscale}
\usepackage{relsize}
\usepackage{enumerate}
\usepackage{epstopdf}
\usepackage{latexsym}
\usepackage[usenames]{color}

\usepackage{graphics}
\usepackage{framed}

\usepackage{enumerate}
\usepackage{hyperref}
\usepackage{xcolor}
\hypersetup{
  colorlinks,
  linkcolor={blue!80!black},
  urlcolor={blue!80!black},
  citecolor={blue!80!black}
}
\usepackage[capitalize,noabbrev]{cleveref}

\numberwithin{equation}{section}

\newtheorem{definition}{Definition}[section]
\newtheorem{theorem}{Theorem}[section]
\newtheorem{lemma}[theorem]{Lemma}

\theoremstyle{definition}

\theoremstyle{remark}
\newtheorem{remark}[theorem]{Remark}

\numberwithin{equation}{section}

\newcounter{saveeqn}

\title[Spectral patterns of elastic transmission eigenfunctions]{Spectral patterns of elastic transmission eigenfunctions: boundary localisation, surface resonance and stress concentration}

\author{Yan Jiang}
\address{Department of Mathematics, Jilin University, Changchun, Jilin, China.}
\email{jiangyan20@mails.jlu.edu.cn}

\author{Hongyu Liu}
\address{Department of Mathematics, City University of Hong Kong, Hong Kong SAR, China.}
\email{hongyu.liuip@gmail.com, hongyliu@cityu.edu.hk}

\author{Jiachuan Zhang}
\address{School of Physical and Mathematical Sciences, Nanjing Tech University, Nanjing, Jiangsu, China.}
\email{zhangjc@njtech.edu.cn}

\author{Kai Zhang}
\address{Department of Mathematics, Jilin University, Changchun, Jilin, China.}
\email{zhangkaimath@jlu.edu.cn}

\date{} 
\begin{document}

\begin{abstract}

We present a comprehensive study of new discoveries on the spectral patterns of elastic transmission eigenfunctions, including boundary localisation, surface resonance, and stress concentration. In the case where the domain is radial and the underlying parameters are constant, we give rigorous justifications and derive a thorough understanding of those intriguing geometric and physical patterns. We also present numerical examples to verify that the same results hold in general geometric and parameter setups.

\medskip



\end{abstract}

\maketitle


\section{Introduction}

\subsection{Mathematical setup and summary of major findings}

Initially focusing on the mathematics, but not the physics, we present the mathematical formulation of the elastic transmission eigenvalue problem and discuss the major findings in this paper.

Let $\Omega$ be a bounded Lipschitz domain in $\mathbb{R}^N$, $N=2, 3$, and $\lambda, \mu$ be real constants fulfilling $\mu>0$ and $2\lambda+N\mu>0$. Define $\mathbf{C}=(C_{ijkl})_{i,j,k,l=1}^N$ to be a four-order tensor with $C_{ijkl}=\lambda \delta_{ij}\delta_{kl}+\mu(\delta_{ik}\delta_{jl}+\delta_{il}\delta_{jk})$, where $\delta$ is the Kronecker delta. Let $\mathbf{u}(\mathbf{x})$, $\mathbf{x}\in\Omega$, be a $\mathbb{C}^N$-valued function and define
\begin{equation}\label{eq:tensors}
\boldsymbol{\varepsilon}(\mathbf{u}):=\frac 1 2[\nabla\mathbf{u}+(\nabla\mathbf{u})^T],\quad \boldsymbol{\sigma}(\mathbf{u}):=\mathbf{C}:\boldsymbol{\varepsilon}(\mathbf{u}),
\end{equation}
where for $\boldsymbol{\varepsilon}=(\varepsilon_{ij})_{i,j=1}^N$,
\begin{equation}\label{eq:op1}
\mathbf{C}:\boldsymbol{\varepsilon}=(\mathbf{C}:\boldsymbol{\varepsilon})_{i,j=1}^N=\left(\sum_{k,l=1}^N C_{ijkl}\varepsilon_{kl}\right)_{i,j=1}^N.
\end{equation}
$\boldsymbol{\varepsilon}$ and $\boldsymbol{\sigma}$ are, respectively, referred to as the strain and Cauchy stress tensors. Let $\rho, \widetilde{\rho}\in L^\infty(\Omega)$ be positive functions and $\omega\in\mathbb{R}_+$. Consider the following PDE system for $\mathbf{u}, \mathbf{v}\in H^1(\Omega)^N$:
\begin{equation}\label{eq:Elas1}
\nabla\cdot\boldsymbol{\sigma}(\mathbf{u})+\omega^2\rho\mathbf{u}=0,\ \nabla\cdot\boldsymbol{\sigma}(\mathbf{v})+\omega^2\widetilde\rho\mathbf{v}=0\ \ \mbox{in}\ \Omega;\ \ \mathbf{u}=\mathbf{v},\ \ \nu\cdot\boldsymbol{\sigma}(\mathbf{u})=\nu\cdot\boldsymbol{\sigma}(\mathbf{v})\ \ \mbox{on}\ \partial\Omega,
\end{equation}
where $\nu\in\mathbb{S}^{N-1}$ signifies the exterior unit normal to $\partial\Omega$. It is clear that $\mathbf{u}=\mathbf{v}\equiv \mathbf{0}$ are a pair of trivial solutions to \eqref{eq:Elas1}. If there exists nontrivial solutions, $\omega\in\mathbb{R}_+$ is referred to as a (real) transmission eigenvalue and the nontrivial $\mathbf{u}$ and $\mathbf{v}$ are the corresponding transmission eigenfunctions. Henceforth, we set $\mathbb{E}_\Omega$ to signify the set of the transmission eigenvalues to \eqref{eq:Elas1}.

In this paper, we are mainly concerned with the spectral patterns of the elastic transmission eigenfunctions. To that end, we first introduce a notion of boundary localisation as follows  (cf. \cite{CDHLW21}).

\begin{definition}
For a sufficiently small $\epsilon\in\mathbb{R}_+$, we set
\begin{equation*}
\mathcal{N}_{\epsilon}\left(\partial\Omega\right):=\left\{\mathbf{x}\in\Omega;\, \mathrm{dist}\left(\mathbf{x},\partial\Omega\right)<\epsilon\right\}.
\end{equation*}
A function $\boldsymbol{\varphi}\in L^2\left(\Omega\right)^N$ is said to be boundary-localized if there exist sufficiently small $\epsilon, \tilde{\epsilon}\in \mathbb{R}_+$ such that
\begin{equation}\label{eq:def2.1}
\frac{\|\boldsymbol{\varphi}\|_{L^2\left(\Omega\right \backslash\mathcal{N}_{\epsilon}\left(\partial\Omega\right))^N}}
{\|\boldsymbol{\varphi}\|_{L^2\left(\Omega\right)^N}}<\tilde{\epsilon}.
\end{equation}
\end{definition}
By (\ref{eq:def2.1}), if a function is boundary-localized, its $L^2$-energy concentrates around the boundary of the domain $\Omega$. In what follows, we shall show that there exists a sequence of elastic transmission eigenfunctions $\left\{\boldsymbol{\varphi}_{m}\right\}_{m\in\mathbb{N}}$ such that for any given $\epsilon\in\mathbb{R}_+$ it holds that
\begin{equation}\label{eq:limit_1}
\lim\limits_{m\rightarrow+\infty}\frac{\|\boldsymbol{\varphi}_m\|_{L^2\left(\Omega\right \backslash\mathcal{N}_{\epsilon}\left(\partial\Omega\right))^N}}
{\|\boldsymbol{\varphi}_m\|_{L^2\left(\Omega\right)^N}}=0.
\end{equation}
In what follows, we simply refer to the sequence $\left\{\boldsymbol{\varphi}_n\right\}_{n\in\mathbb{N}}$ satisfying (\ref{eq:limit_1}) as boundary-localized.


In what follows, we let $m\in \mathbb{N}$ be a positive integer and $J_m(\mathbf{x})$ be the first kind Bessel function of order $m$. Moreover, we let $j_{m,s}$ and $j'_{m,s}$ be the $s$-th positive zero of $J_m(\mathbf{x})$ and $J'_m(\mathbf{x})$, respectively.

The first main discovery in this paper is contained in the following theorem.
\begin{theorem}\label{monolocalized}
Consider the elastic transmission eigenvalue problem \eqref{eq:Elas1}. Assume that $\Omega$ is the unit ball in $\mathbb{R}^N$, $N=2,3$, and $\rho, \widetilde\rho$ are positive constants with $\rho\neq\widetilde\rho$. Then there exists a sequence $\left\{\omega_m\right\}_{m\in\mathbb{N}}\subset \mathbb{E}_{\Omega}$ such that for sufficiently large $m$, it holds that
\begin{equation}\label{eq:mono_omega_m}
\omega_m\in\left(\sqrt{\frac{\mu}{\rho}}j_{m+R_N,s_{1}(m)}
,\sqrt{\frac{\mu}{\rho}}j_{m+R_N,s_{2}(m)}\right),
\end{equation}
where
\begin{equation}\label{eq:RN}
R_N=0\ \ \mbox{when}\ \ N=2;\ \ 1/2\ \ \mbox{when}\ \ N=3,
\end{equation}
and $j_{m+R_N,s_1(m)}$ and $j_{m+R_N,s_2(m)}$ denote the $s_1(m)$-th and $s_2(m)$-th positive roots of the Bessel function $J_{m+R_N}(\mathbf{x})$ for a fixed order $m+R_N$, respectively.
The parameters $s_1(m)$ and $s_2(m)$ have the following formulas:
\begin{equation}\label{eq:s1_s2}
s_{1}(m):=\left[m^{\gamma_{1}}\right], \quad s_{2}(m)=\left[m^{\gamma_{2}}\right], \quad 0<\gamma_{1}<\gamma_{2}<1,
\end{equation}
where $[t]$ denotes the largest integer smaller than $t$, and the transmission eigenvalues in \eqref{eq:mono_omega_m} fulfil that
\begin{equation}\label{eq:asyn1}
\omega_m=\sqrt{\frac{\mu}{\rho}}m\left (1+\mathcal{O}(m^{-\frac{2}{3}+\frac{2}{3}\gamma})\right), \quad 0<\gamma_{1}\leq\gamma\leq\gamma_{2}<1.
\end{equation}
 Moreover, if we let $(\mathbf{u}_{m}$, $\mathbf{v}_{m})$ be the pair of transmission eigenfunctions associated with $\omega_m$, then it holds that $\{\mathbf{u}_{m}\}_{m\in\mathbb{N}}$ is boundary-localized, whereas $\{\mathbf{v}_{m}\}_{m\in\mathbb{N}}$ is not boundary-localized.
\end{theorem}

In what follows, the sequence of special eigen-modes $\{\mathbf{u}_m, \mathbf{v}_m\}_{m\in\mathbb{N}}$ found in Theorem~\ref{monolocalized} is referred to as mono-localized.

\begin{theorem}\label{bilocalized}
Under the same setup of Theorem~\ref{monolocalized}, for any given $s_0\in\mathbb{N}$, there exists a sequence $\left\{\omega_m\right\}_{m\in\mathbb{N}}\subset \mathbb{E}_{\Omega}$ such that for sufficiently large $m$, it holds that
\begin{equation}\label{eq:bi_omega_m}
\omega_m\in\left(\sqrt{\frac{\mu}{\rho}}\frac{j_{m+R_N,s_{0}}}{\mathbf{n}}
,\sqrt{\frac{\mu}{\rho}}\frac{j_{m+R_N,s_{0}+1}}{\mathbf{n}}\right),\quad\omega_m=\sqrt{\frac{\mu}{\rho}}\frac{m}{\mathbf{n}}\left (1+\mathcal{O}(m^{-\frac{2}{3}})\right),\quad \mathbf{n}:=\sqrt{\frac{\tilde{\rho}}{\rho}},
\end{equation}
where $j_{m+R_N,s_0}$ and $j_{m+R_N,s_0+1}$ denote the $s_0$-th and $(s_0+1)$-th positive roots of the Bessel function $J_{m+R_N}(\mathbf{x})$ for a fixed order $m+R_N$, respectively. Moreover, if we let $(\mathbf{u}_{m}$, $\mathbf{v}_{m})$ be the pair of transmission eigenfunctions associated with $\omega_m$ in \eqref{eq:bi_omega_m}, then it holds that both $\mathbf{u}_{m}$ and $\mathbf{v}_{m}$ are boundary-localized, and $\{\mathbf{u}_m, \mathbf{v}_m\}_{m\in\mathbb{N}}$ is referred to as bi-localized.
\end{theorem}


{It is known that $\mathbf{u}$ and $\mathbf{v}$ possess the following vectorial decomposition in three dimensions:
\begin{equation}\label{eq:decompn1}
\mathbf{w}=\mathbf{w}^p+\mathbf{w}^s,\ \ \mathbf{w}^p:=-\frac{1}{k_p^2}\nabla(\nabla\cdot\mathbf{w}), \ \ \mathbf{w}^s:=\frac{1}{k_s^2}\nabla\times(\nabla\times\mathbf{w}),\quad \mathbf{w}=\mathbf{u}\ \mbox{or}\ \mathbf{v},
\end{equation}
where $k_p:=\omega/\sqrt{\lambda+2\mu}$ and $k_s:=\omega/\sqrt{\mu}$. Here, $\mathbf{w}_p$ and $\mathbf{w}_s$ are known as the compressional and shear parts of the elastic field $\mathbf{w}$ respectively. }
In two dimensions, the decomposition is given as
\begin{equation}\label{eq:decompn2}
\mathbf{w}=\mathbf{w}^p+\mathbf{w}^s,\ \ \mathbf{w}^p:=-\frac{1}{k_p^2}\nabla(\nabla\cdot\mathbf{w}), \ \ \mathbf{w}^s:=\frac{1}{k_s^2}\left(\partial_{xy}^{2}{w}_{2}
-\partial_{y}^{2}{w}_{1},\partial_{xy}^{2}{w}_{1}
-\partial_{x}^{2}{w}_{2}
\right)^{T},
\end{equation}
for $\mathbf{x}=(x, y)\in\mathbb{R}^2$ and $\mathbf{w}=(w_1, w_2)$ with  $\mathbf{w}=\mathbf{u}\ \mbox{or}\ \mathbf{v}$.
For the boundary-localized eigen-modes determined in Theorems~\ref{monolocalized} and \ref{bilocalized}, we have the following boundary-localization properties regarding their compressional and shear parts.

\begin{theorem}\label{shear_and_compressional}
Let $\{\mathbf{u}_m, \mathbf{v}_m\}_{m\in\mathbb{N}}$ be determined in Theorem~\ref{monolocalized} or \ref{bilocalized}, and $\{\mathbf{u}_m^p, \mathbf{v}_m^p\}_{m\in\mathbb{N}}$ and $\{\mathbf{u}_m^s, \mathbf{v}_m^s\}_{m\in\mathbb{N}}$ be respectively the compressional and shear parts. Without loss of generality, we assume that $\mathbf{n}$ introduced in \eqref{eq:bi_omega_m} satisfies that $\mathbf{n}>1$. Table~\ref{tab:1n1} lists the properties of those eigen-modes as well as their respective compressional and shear parts whether they are boundary-localized or not.

\begin{table}[!htbp]
\caption{Boundary-localizing properties of transmission eigenfunctions as well as their compressional and shear parts}\label{tab:1n1}
\begin{tabular}{|c|c|c|c|c|}
\hline \multirow{2}{*}{} & \multicolumn{2}{|c|}{\text{\emph{mono-localized modes} }} & \multicolumn{1}{c|}{{\emph{bi-localized modes} }} \\
\hline &
$\mathbf{n}\sqrt{\frac{\mu}{\lambda+2\mu}}\leq1 $   & $\mathbf{n}\sqrt{\frac{\mu}{\lambda+2\mu}}>1    $   &    $\mathbf{n}>1 $      \\
\hline  $\mathbf{u}_{m}$     & \text{\emph{Yes}} & \text{\emph{Yes}} & \text{\emph{Yes}}\\
\hline  $\mathbf{u}^{s}_{m}$ & \text{\emph{Yes}} & \text{\emph{Yes}} & \text{\emph{Yes}}\\
\hline  $\mathbf{u}^{p}_{m}$ & \text{\emph{Yes}} & \text{\emph{Yes}} & \text{\emph{Yes}}\\
\hline  $\mathbf{v}_{m}$     & \text{\emph{No}}  & \text{\emph{No}}  & \text{\emph{Yes}}\\
\hline  $\mathbf{v}^{s}_{m}$ & \text{\emph{No}}  & \text{\emph{No}}  & \text{\emph{Yes}}\\
\hline  $\mathbf{v}^{p}_{m}$ & \text{\emph{Yes}} & \text{\emph{No}}  & \text{\emph{Yes}}\\
\hline
\end{tabular}
\end{table}
\end{theorem}

\begin{remark}\label{rem:1.4}
It is interesting to observe from Table~\ref{tab:1n1} that the boundary-localizing property of the compressional and shear parts of an elastic field is the same as that of the elastic field itself, except for the mono-localized $(\mathbf{u}_m, \mathbf{v}_m)$ where the non-boundary-localization of $\mathbf{v}_m$ is mainly caused by the non-boundary-localization of $\mathbf{v}_m^s$ but not $\mathbf{v}_m^p$ when $\mathbf{n}\sqrt{\mu/(\lambda+2\mu)}<1$.
\end{remark}

Finally, we present quantitative properties of the boundary-localized eigen-modes when treated as surface waves and show that they are surface resonant waves, accompanying strong stress field concentration. To that end, we introduce the following subdomain $\Sigma$ of $\Omega$ as
\begin{equation}\label{eq:theta12}
\Sigma(\tau,\theta_1,\theta_2)=\left\{(r,\theta)|\tau<r<1,
~\theta_1<\theta<\theta_2\right\},
\end{equation}
where $\mathbf{x}=(r, \theta)$ stands for the polar coordinates in $\mathbb{R}^2$, as well as the following quantities:
\begin{equation}\label{eq:terms}
\begin{array}{ll}
\displaystyle{\left|\nabla \mathbf{u}_m \right|^2_{\Sigma, \infty}:=\max_{\mathbf{x}\in\Sigma}\left|\nabla\mathbf{u}_m(\mathbf{x})\right|}, & \displaystyle{E^{2}_{\mathbf{u}_m, \Sigma}=\int_{\Sigma}\left(\boldsymbol{\sigma}\left(\mathbf{u}_m\right): \nabla \overline{\mathbf{u}_{m}}\right)dr d\theta,}\\
\displaystyle{\left|\nabla \mathbf{v}_m \right|^2_{\Sigma, \infty}=\max_{\mathbf{x}\in\Sigma}\left|\nabla\mathbf{v}_m(\mathbf{x})\right|,} & \displaystyle{E^{2}_{\mathbf{v}_m, \Sigma}=\int_{\Sigma}\left(\boldsymbol{\sigma}\left(\mathbf{v}_m\right): \nabla \overline{\mathbf{v}_{m}}\right)dr d\theta,}
\end{array}
\end{equation}
where the operation $``:"$ is understood as for two matrices $\mathbf{A}=(a_{ij})_{i,j=1}^N$ and $\mathbf{B}=(b_{ij})_{i,j=1}^N$, $\mathbf{A}:\mathbf{B}=\sum_{i,j=1}^N a_{ij}b_{ij}$. Here, we note that with a bit notational abuse, the operation $``:"$ is defined in both \eqref{eq:op1} and \eqref{eq:terms}, which should be clear from the context.

\begin{theorem}\label{surface_resonant}
Consider the bi-localized eigen-modes $\{\mathbf{u}_m, \mathbf{v}_m\}_{m\in\mathbb{N}}$ and suppose that $\mathbf{v}_m$ is $L^2(\Omega)$-normalised, i.e. $\|\mathbf{v}_m\|_{L^2(\Omega)}=1$. Let $\Sigma$ be defined in \eqref{eq:theta12} with $\tau\in (0, 1)$ and $\theta_1, \theta_2\in (0, 2\pi)$ being given and fixed. Suppose that $\lambda, \mu$ and $\rho, \tilde\rho$ are fixed. Then it holds for $m\in\mathbb{N}$ sufficiently large that
\begin{equation}\label{eq:surface_resonant2}
\left|\nabla \mathbf{u}_m \right|^2_{\infty, \Sigma}\geq C_{1}m^{3},\quad \left|\nabla \mathbf{v}_m \right|^2_{\infty, \Sigma}\geq C_{2}m^{\frac{10}{3}},\quad E_{\mathbf{u}_{m}, \Sigma}^{2}\geq C_{3}\mu m^3,\quad E_{\mathbf{v}_{m}, \Sigma}^{2}\geq C_{4}\mu m^{\frac{10}{3}},
\end{equation}
where $C_1$ and $C_2$ are positive constants depending only on $\lambda,\mu, \rho, \tilde\rho$; and $C_3$ and $C_4$ are positive constants depending only on $\lambda, \rho, \tilde\rho$ and $\Sigma$. On the other hand, for fixed $\lambda, \rho, \tilde\rho$ and $m\in\mathbb{R}_+$, it holds for $\mu\in\mathbb{R}_+$ sufficiently large that
\begin{equation}\label{eq:surface_resonant1}
E_{\mathbf{u}_{m}, \Sigma}^{2}\geq C_{5}\mu ,\quad E_{\mathbf{v}_{m}, \Sigma}^{2}\geq C_{6}\mu ,
\end{equation}
where $C_5$ and $C_6$ are positive constants depending only on $\lambda, \rho, \tilde\rho$ and $\Sigma$.
\end{theorem}

\begin{remark}\label{rem:n1}
Theorem~\ref{surface_resonant} shows that if treated as surface waves propagating along $\partial\Omega$, the bi-localized transmission eigenfunctions $\mathbf{u}_m$ and $\mathbf{v}_m$ form certain resonant modes, manifesting highly-oscillatory patterns along with the energy blowup. In fact, according to \eqref{eq:bi_omega_m}, we see that the ``natural" frequency for $\mathbf{u}_m$ and $\mathbf{v}_m$ is $\omega_m\sim m$, whereas they oscillate much more severely along $\partial\Omega$, with the surface-oscillating frequency bigger than $m^{3/2}$ for $\mathbf{u}_m$ and $m^{5/3}$ for $\mathbf{v}_m$. The resonance phenomenon is further corroborated by the blowup of the stress energies ${E}_{\mathbf{u}_m, \Sigma}$ and ${E}_{\mathbf{v}_m, \Sigma}$. Moreover, the estimates in \eqref{eq:surface_resonant1} indicate that the surface resonance can be more evident if the Lam\'e parameter $\mu$ is large. In fact, since the generic constants $C_5$ and $C_6$ are independent of $m$, and hence $\omega_m$, the estimates in \eqref{eq:surface_resonant1} indicate that even for those low-mode-number transmission eigenfunctions, if they are boundary-localized, they exhibit highly oscillatory surface-resonant behaviours provided $\mu$ is large.
Finally, in Section~\ref{sect:numerics}, we further verify such surface resonance properties by numerics.
\end{remark}

\begin{remark}\label{rem:n2}
In Theorem~\ref{surface_resonant}, we only consider the bi-localised eigenfunctions determined in Theorem~\ref{bilocalized}. Nevertheless, by following a similar argument and in principle, one can show similar surface-resonant properties of the mono-localized eigenfunctions determined in Theorem~\ref{monolocalized}. That is, as long as the transmission eigen-mode is boundary-localised, no matter the $\mathbf{u}$-part or the $\mathbf{v}$-part, it exhibits highly oscillatory surface-resonant patterns. Moreover, the surface-resonance results also holds for the three-dimensional cases. However, a complete description of those results is  lengthy and moreover the corresponding verifications involve tedious and sometime repeating calculations. Hence, we choose to stick to the two-dimensional bi-localised transmission eigen-modes to study the surface-resonance properties.
\end{remark}

%


\subsection{Physical relevance and background discussion}

To motivate the current study, we consider the time-harmonic elastic scattering from an inhomogeneous inclusion embedded in a homogeneous background space. Let $(\mathbb{R}^N; \lambda, \mu, \rho)$ and $(\Omega; \lambda,\mu, \tilde\rho)$ be respectively specify the medium configurations of the background space and the inhomogeneous medium inclusion. Here, $\lambda, \mu$ and $\rho$ are constants, fulling $\mu>0, 2\lambda+N\mu>0$ and $\rho>0$, which characterise the bulk moduli and the density of the elastic material. We also assume that $\tilde\rho\in L^\infty(\Omega)$ and $\tilde\rho>0$. Let $\mathbf{u}^i$ be an entire solution to $\nabla\cdot\boldsymbol{\sigma}(\mathbf{u}^i)+\omega^2\rho\mathbf{u}^i={0}$ in $\mathbb{R}^N$, which signifies an incident field. The impingement of $\mathbf{u}^i$ on $(\Omega; \lambda,\mu, \tilde\rho)$ generates the elastic scattering with the physical wave field $\mathbf{v}\in H_{loc}^1(\mathbb{R}^N)^N$ fulfilling the following transmission problem:
 \[
 \nabla\cdot\boldsymbol{\sigma}(\mathbf{v})+\omega^2\tilde\rho\mathbf{v}=0\ \mbox{in}\ \Omega;\ \nabla\cdot\boldsymbol{\sigma}(\mathbf{v})+\omega^2\rho\mathbf{v}=0\ \mbox{in}\ \mathbb{R}^N\backslash\overline{\Omega};\ \mathbf{v}|^+=\mathbf{v}|^-,\ \nu\cdot\boldsymbol{\sigma}(\mathbf{v})|^+=\nu\cdot\boldsymbol{\sigma}(\mathbf{v})|^-\ \mbox{on}\ \partial\Omega,
 \]
 where $\pm$ signify the traces on $\partial\Omega$ taken from $\mathbb{R}^N\backslash\overline{\Omega}$ and $\Omega$, respectively. An inverse problem of practical importance is to recover $(\Omega; \tilde\rho)$ by knowledge of the scattering pattern outside the scatterer, namely $(\mathbf{v}-\mathbf{u}^i)|_{\mathbb{R}^N\backslash\overline{\Omega}}$. We refer to \cite{MBDL} for more related discussions on the forward and inverse elastic problems. However, we are curious about the wave patterns when invisibility/transparency occurs, namely $(\mathbf{v}-\mathbf{u}^i)|_{\mathbb{R}^N\backslash\overline{\Omega}}\equiv 0$, or equivalently $\mathbf{v}=\mathbf{u}^i$ in $\mathbb{R}^N\backslash\overline{\Omega}$. In such a case, one can directly verify that $\mathbf{u}=\mathbf{u}^i|_{\Omega}$ and $\mathbf{v}|_{\Omega}$ fulfils \eqref{eq:Elas1}. That is, when invisibility/transparency occurs, the scattering patterns are trapped inside the scatterer to form the transmission eigenfunctions.

The spectral study of transmission eigenvalue problems arising in the wave scattering theory has a long and colourful history; see \cite{CCH,Liu} and the references cited therein. However, the spectral patterns of transmission eigenfunctions were only unveiled recently. In \cite{CDHLW21,CDLS,JLZZ22}, it is shown that acoustic transmission eigenfunctions exhibit the boundary-localisation phenomenon. The boundary-localising properties were further extended to the electromagnetic and acoustic-elastic transmission eigenfunctions in \cite{DLWW22,DLLT}. It is noted that in all of the aforementioned literature, the boundary-localisation was rigorously justified for the radial geometry, whereas for the case with general geometries, it is mainly verified numerically with the only exception in \cite{CDLS} where the boundary-localisation was theoretically justified in two different senses for the acoustic transmission eigenfunctions. In addition, it is shown in \cite{B,BL1,BL,BLX,DCL,DLS} that transmission eigenfunctions associated with different wave systems exhibit locally vanishing patterns around corners or high-curvature places around $\partial\Omega$.

The current study follows a similar spirit to that in \cite{CDHLW21,CDLS,DLLT} on the boundary-localisation of transmission eigenfunctions. However, we would like to highlight several novel mathematical and physical developments due to the new setting. First, in addition to revealing the boundary-localisation of the elastic transmission eigenfunctions in Theorems~\ref{monolocalized} and \ref{bilocalized}, we further explore the boundary-localising properties of the corresponding shear and compressional parts in Theorem~\ref{shear_and_compressional} and Remark~\ref{rem:1.4}, which provide a more in-depth and physically relevant understanding of the boundary-localisation of elastic transmission eigenfunctions. Second, it is the first time in the literature to discover that the boundary-localised transmission eigenfunctions exhibit surface-resonant behaviours, accompanying highly-oscillatory pattern as well as  strong stress energy concentration; that is, Theorem~\ref{surface_resonant} and Remarks~\ref{rem:n1} and \ref{rem:n2}.

Finally, we briefly discuss two practical implications of our results. In fact, the spectral patterns of transmission eigenfunctions have already produced several interesting applications of practical importance. In \cite{CDLS}, a new interpretation of the invisibility cloaking was given in terms of the boundary-localisation of the acoustic transmission eigenfunctions. In \cite{DLWW22}, a scheme of generating artificial mirage was proposed based on using the boundary localisation of the electromagnetic transmission eigenfunctions. In \cite{CDHLW21}, a novel acoustic wave imaging scheme was proposed by using the boundary-localising properties of the acoustic transmission eigenfunctions, and it was numerically observed that super-resolution effects can be achieved. By following a similar spirit, one can make use of the boundary-localising properties in Theorems~\ref{monolocalized} and \ref{bilocalized} for the inverse elastic problem discussed above of imaging $\Omega$ by knowledge of
$(\mathbf{v}-\mathbf{u}^i)|_{\mathbb{R}^N\backslash\overline{\Omega}}$. The general procedure can be roughly described as recovering those trapped transmission eigen-modes by using the exterior wave measurement, and then using the boundary-localising behaviours to identify $\partial\Omega$. However, we can provide a rigorous justification on the super-resolution imaging effect that it can produce. In fact, it is known that if one uses $(\mathbf{v}-\mathbf{u}^i)|_{\mathbb{R}^N\backslash\overline{\Omega}}$ to image $\Omega$, the resolution limit is determined by the wavelength of the signal collected, which is in turn determined by the ``natural" frequency $\omega$. However, by Theorem~\ref{surface_resonant}, we know those trapped transmission eigen-modes oscillate with a much smaller wavelength due to the surface resonance (see also the numerical demonstrations in Section~\ref{sect:numerics}), and hence using them, one should be able to see much finer details of $\partial\Omega$. We shall present a more comprehensive study along this direction in a forthcoming paper. The other interesting implication is related to the stress concentration in Theorem~\ref{surface_resonant}. In fact, it is known that a strong stress concentration can cause the failure of an elastic structure. Hence, our result in Theorem~\ref{surface_resonant} indicate that one can generate desired stress concentrations to crack elastic structures. We shall explore more on this aspect in our future work.

The rest of the paper is organised as follows. Sections~\ref{sect:2} and \ref{sect:3} are devoted to the proofs of the main theorems in two and three dimensions, respectively. In Section~\ref{sect:numerics}, we present the numerical results as well as some relevant discussions.


\section{Proofs of main theorems in two dimensions }\label{sect:2}

In this section, we prove the main theorems in two dimensions.


\subsection{Auxiliary results}

We introduce some properties of the Bessel function.
First, the roots of Bessel function have the following sharp bounds and relationships.
\begin{lemma}(\hspace{-0.2mm}\cite{AS72,WQ99}) For $m\in \mathbb{N}$, it holds that
\begin{equation}\label{eq:boundj_ms}
m-\frac{a_{s}}{2^{1 / 3}}m^{1 / 3}<j_{m, s}<m-\frac{a_{s}}{2^{1 / 3}}m^{1 / 3}+\frac{3}{20} a_{s}^{2} \frac{2^{1 / 3}}{m^{1 / 3}},
\end{equation}
where $a_{s}$ is the $s$-th negative zero of the Airy function and has the representation
\begin{equation}\label{eq:bounda_s}
a_{s}=-\left(\frac{3 \pi}{8}(4 s-1)\right)^{2 / 3}\left(1+\sigma_{s}\right), \quad 0 \leq \sigma_{s} \leq 0.130\left(\frac{3 \pi}{8}(4 s-1.051)\right)^{-2}.
\end{equation}
Moreover, the roots satisfy following estimate
\begin{equation}\label{eq:roots}
m<j'_{m,1}<j_{m,1}<j'_{m,2}<j_{m,2}<\cdots.
\end{equation}
\end{lemma}

Second, we give several formulas for Bessel functions.
\begin{lemma}(\hspace{-0.2mm}\cite{AS72,KRA06,KOR02,PAR84}) For $m$ sufficiently large and $0<x<1$, it holds that
\begin{equation}\label{eq_Jm_mx}
J_m(mx)=\left(\frac{4\zeta}{1-x^2}\right)^{\frac{1}{4}}\frac{Ai\left(
m^{\frac{2}{3}}\zeta\right)}{m^{\frac{1}{3}}}\left(
1+\mathcal{O}\left(\frac{1}{m^{\frac{4}{3}}}\right)\right),
\end{equation}
where $Ai(x)$ is the Airy function and $\zeta=\left(\displaystyle{\frac{3}{2}}
\left(\ln\left(\frac{1+\sqrt{1-x^2}}{x}\right)-
\sqrt{1-x^2}\right)\right)^{\frac{2}{3}}$.

For $m$ sufficiently large and $x>m$, following asymptotic formulas hold
\begin{eqnarray}\label{eq:asympt_Jm_x_large_m}
J_{m}(x)&=&\sqrt{\frac{2}{\pi \sqrt{x^{2}-m^{2}}}} \cos \left(\sqrt{x^{2}-m^{2}}-\frac{m \pi}{2}+m \arcsin (m / x)-\frac{\pi}{4}\right)\left(1+\mathcal{O}\left(\frac{1}{m}
\right)\right),\\
\label{eq:asympt_Jm'_x_large_m}
J_{m}'(x)&=&-\sqrt{\frac{2\sqrt{x^{2}-m^{2}}}{\pi x^2}} \cos \left(\sqrt{x^{2}-m^{2}}+m \arccos (m / x)-\frac{3\pi}{4}\right)\left(1+\mathcal{O}\left(\frac{1}{m}
\right)\right).
\end{eqnarray}

Furthermore, one has
\begin{equation}\label{eq:paris1}
0<\frac{x m}{2 m+2}<\frac{J_{m+1}(m x)}{J_{m}(m x)}<\frac{x m}{m+2}<1, \quad \quad 0<x \leq 1, m=1,2,\cdots
\end{equation}
and
\begin{equation}\label{eq:lowerboundofJm'_Jm}
\frac{\sqrt{(2 m+1)^{2}-4 x^{2}}-1}{2 x} \leq \frac{J_{m}'(x)}{J_{m}(x)} <\frac{m}{x}  , \quad 0 \leq x \leq m+\frac{1}{2}, m=0,1,\cdots.
\end{equation}
\end{lemma}

Third, using the integration by parts and the recurrence relation of the Bessel functions (formula (9.1.27) in \cite{AS72}),
we can deduce two useful integral formulas.
\begin{lemma}
For $m\in \mathbb{N}$, following formulas hold
\begin{equation}\label{eq:integrabyparts}
\begin{array}{lll}
\hspace{-1mm}\displaystyle{\int_{0}^{\tau}\left(k^2J_{m}^{'2}(kr)r
+m^2\frac{J_{m}^2(kr)}{r}\right)\mathrm{d}r}=\displaystyle{\int_{0}^{\tau}k^2J_{m-1}^{2}(kr)
r\mathrm{d}r}-mJ_{m}^2(k\tau),~\forall~0<\tau<1,~k>0,
\end{array}
\end{equation}
and
\begin{equation}\label{eq:integrabyparts2}
\int_{0}^{\tau}rJ^2_m\left(r\right)\mathrm{d}r=\frac{1}{2}
\tau^2J_{m}^{'2}\left(\tau\right)
+\frac{1}{2}\left(m^2-\tau^2\right)J_{m}^{2}\left(\tau\right),~\forall~\tau>0.
\end{equation}
\end{lemma}

Finally, we give an expansion for spherical harmonic functions.
\begin{lemma}(\hspace{-0.2mm}\cite{AS72})
Let $\{Y_m^n\}_{m=0,\cdots, \infty}^{n=-m,\cdots,m}$ be the spherical harmonic functions. Then
\begin{equation}\label{eq:shfexpansion}
Y_m^n\left(\theta,\phi\right):=\sqrt{\frac{2m+1}{4\pi}\frac{\left(m-|n|\right)!}{\left(m+|n|\right)!}}P_m^{|n|}\left(\cos\theta\right)e^{\mathrm{i}n\phi},
\end{equation}
where the associated Legendre functions satisfy
\begin{equation}\label{eq:Lpro}
\int_{-1}^{1} \frac{P_{m}^{n_1} P_{m}^{n_2}}{1-x^{2}} d x
= \begin{cases}0, & \text { if } n_1 \neq n_2, \\
\frac{(m+n_1) !}{n_1(m-n_1) !}, & \text { if } n_1=n_2 \neq 0, \\
\infty, & \text { if } n_1=n_2=0.
\end{cases}
\end{equation}
\end{lemma}

In the rest of this subsection, we construct a function whose roots are transmission eigenvalues in $\mathbb{R}^2$.
Consider the Helmholtz decomposition of $\mathbf{u}$ and $\mathbf{v}$, i.e. $\mathbf{u}=\nabla\phi+\nabla\times\xi$ and $\mathbf{v}=\nabla\psi+\nabla\times\zeta$,
then
\begin{equation}\label{eq:auxiliary1}
\left\{
\begin{array}{llll}
\Delta\phi+\frac{\rho\omega^2}{2\mu+\lambda}\phi=0, &\Delta\psi+\frac{\widetilde{\rho}\omega^2}{2\widetilde{\mu}
+\widetilde{\lambda}}\psi=0  & \text{in} \  \Omega, \\
\Delta\xi+\frac{\rho\omega^2}{\mu}\xi=0, &\Delta\zeta+\frac{\widetilde{\rho}\omega^2}{\widetilde{\mu}}\zeta=0  & \text{in} \  \Omega.
\end{array}
\right.
\end{equation}
Since $\rho$ and $\tilde{\rho}$ are  positive constants, the Fourier
series of the above quantities yield
\begin{equation}\label{eq:auxiliary notation}
\left\{
\begin{array}{llll}
&\phi(x)=\sum\limits_{m=0}^\infty \alpha_m J_{m}(k_1|x|)e^{\mathrm{i}m\theta},
&\psi(x)=\sum\limits_{m=0}^\infty \beta_m J_{m}(\widetilde{k}_1|x|) e^{\mathrm{i}m\theta},\\
&\xi(x)=\sum\limits_{m=0}^\infty \gamma_m J_{m}(k_2|x|) e^{\mathrm{i}m\theta},
&\zeta(x)=\sum\limits_{m=0}^\infty \delta_m J_{m}(\widetilde{k}_2|x|) e^{\mathrm{i}m\theta},\\
\end{array}
\right.
\end{equation}
where $k_1^2=\frac{\rho\omega^2}{2\mu+\lambda}$, $\widetilde{k}_1^2=\frac{\widetilde{\rho}\omega^2}{2\widetilde{\mu}
+\widetilde{\lambda}}$, $k_2^2=\frac{\rho\omega^2}{\mu}$, and $\widetilde{k}_2^2=\frac{\widetilde{\rho}\omega^2}{\widetilde{\mu}}$.

For any $m\geq 1$, we set
\begin{equation*}
\left\{
\begin{array}{llll}
&\phi_m(x)=\alpha_m J_{m}(k_1|x|) e^{\mathrm{i}m\theta},
&\psi_m(x)=\beta_m J_{m}(\widetilde{k}_1|x|) e^{\mathrm{i}m\theta},
\\
&\xi_m(x)=\gamma_m J_{m}(k_2|x|) e^{\mathrm{i}m\theta},
&\zeta_m(x)=\delta_m J_{m}(\widetilde{k}_2|x|) e^{\mathrm{i}m\theta},\\
&\mathbf{u}_m(x)=\nabla\phi_m(x)+\nabla\times\xi_m(x),
&\mathbf{v}_m(x)=\nabla\psi_m(x)+\nabla\times\zeta_m(x).
\end{array}
\right.
\end{equation*}
It follows from the first boundary condition in \eqref{eq:Elas1} and \eqref{eq:auxiliary notation} that
\begin{equation}\label{eq:boundary1}
\left\{
\begin{array}{lll}
\alpha_m k_1 J^{\prime}_{m}(k_1)-\gamma_m \mathrm{i}mJ_{m}(k_2)&=&\beta_m \widetilde{k}_1 J^{\prime}_{m}(\widetilde{k}_1)-\delta_m \mathrm{i}mJ_{m}(\widetilde{k}_2),\\
\alpha_m \mathrm{i}mJ_{m}(k_1)+\gamma_m k_2 J^{\prime}_{m}(k_2)&=&\beta_m \mathrm{i}mJ_{m}(\widetilde{k}_1)+\delta_m \widetilde{k}_2 J^{\prime}_{m}(\widetilde{k}_2),
\end{array}
\right.
\end{equation}
while the second boundary condition in \eqref{eq:Elas1} implies
\begin{equation}\label{eq:boundary2}
\left\{
\begin{array}{ll}
&\alpha_m((\rho\omega^2-2\mu m^2)J_m(k_1)+2\mu k_1J_m^{\prime}(k_1))+\gamma_m2\mu\mathrm{i}m (k_2J_m^{\prime}(k_2)-J_m(k_2))\\
=&\beta_m((\widetilde{\rho}\omega^2-2\widetilde{\mu} m^2)J_m(\widetilde{k}_1)+2\widetilde{\mu} \widetilde{k}_1J_m^{\prime}(\widetilde{k}_1))
+\delta_m2\widetilde{\mu}\mathrm{i}m (\widetilde{k}_2J_m^{\prime}(\widetilde{k}_2)-J_m(\widetilde{k}_2)),\\
&\alpha_m(-2\mu\mathrm{i}m (k_1J_m^{\prime}(k_1)-J_m(k_1)))+\gamma_m((\rho\omega^2-2\mu m^2)J_m(k_2)+2\mu k_2J_m^{\prime}(k_2))\\
=&\beta_m(-2\widetilde{\mu}\mathrm{i}m (\widetilde{k}_1J_m^{\prime}(\widetilde{k}_1)
-J_m(\widetilde{k}_1)))+\delta_m((\widetilde{\rho}\omega^2
-2\widetilde{\mu} m^2)J_m(\widetilde{k}_2)+2\widetilde{\mu} \widetilde{k}_2J_m^{\prime}(\widetilde{k}_2)).\\
\end{array}
\right.
\end{equation}
Combining \eqref{eq:boundary1} and \eqref{eq:boundary2}, one can see that $\omega$ is a transmission eigenvalue if
\begin{equation}\label{eq:fmomega}
f_m(\omega):=
\left|
\begin{array}{llll}
k_1 J^{\prime}_{m}(k_1) & - \mathrm{i}mJ_{m}(k_2) & -\widetilde{k}_1 J^{\prime}_{m}(\widetilde{k}_1)& \mathrm{i}mJ_{m}(\widetilde{k}_2)\\
\mathrm{i}mJ_{m}(k_1) & k_2 J^{\prime}_{m}(k_2) & -\mathrm{i}mJ_{m}(\widetilde{k}_1) & -\widetilde{k}_2 J^{\prime}_{m}(\widetilde{k}_2)\\
a & b & c & d\\
e & f & g & h
\end{array}
\right|=0,
\end{equation}
where the parameters are specified as
\begin{equation*}
\begin{array}{lll}
&a=(\rho\omega^2-2\mu m^2)J_m(k_1)+2\mu k_1J_m^{\prime}(k_1),
&b=2\mu\mathrm{i}m (k_2J_m^{\prime}(k_2)-J_m(k_2)),\\
&c=-((\widetilde{\rho}\omega^2-2\widetilde{\mu} m^2)J_m(\widetilde{k}_1)+2\widetilde{\mu} \widetilde{k}_1J_m^{\prime}(\widetilde{k}_1)),
&d=-2\widetilde{\mu}\mathrm{i}m (\widetilde{k}_2J_m^{\prime}(\widetilde{k}_2)-J_m(\widetilde{k}_2)),\\
&e=-2\mu\mathrm{i}m (k_1J_m^{\prime}(k_1)-J_m(k_1)),
&f=(\rho\omega^2-2\mu m^2)J_m(k_2)+2\mu k_2J_m^{\prime}(k_2),\\
&g=2\widetilde{\mu}\mathrm{i}m (\widetilde{k}_1J_m^{\prime}(\widetilde{k}_1)-J_m(\widetilde{k}_1)),
&h=-((\widetilde{\rho}\omega^2-2\widetilde{\mu} m^2)J_m(\widetilde{k}_2)+2\widetilde{\mu} \widetilde{k}_2J_m^{\prime}(\widetilde{k}_2)).
\end{array}
\end{equation*}


\subsection{Mono-localization when $N=2$}


We are in a position to prove Theorem~\ref{monolocalized} when $N=2$. The proof is divided into two parts. The first part is to estimate the interval where the eigenvalues are located, and the second part is to prove that the boundary-localising pattern of the eigenfunctions.

\begin{proof}
Part 1. We first study the eigenvalue distribution. We consider the following three cases one by one
$$
1<\mathbf{n}<\sqrt{\frac{2\mu+\lambda}{\mu}}; \quad
1<\mathbf{n}=\sqrt{\frac{2\mu+\lambda}{\mu}}; \quad
1<\sqrt{\frac{2\mu+\lambda}{\mu}}<\mathbf{n}.
$$

$\mbox{Case 1:}~1<\mathbf{n}<\sqrt{\frac{2\mu+\lambda}{\mu}}$.~~
From \eqref{eq:fmomega}, we have
\begin{equation}\label{eq:fm_s1fm_s2_case1}
\begin{array}{lll}
&&f_{m}\left(\sqrt{\frac{\mu}{\rho}} j_{m, s_{1}}\right)f_{m}\left(\sqrt{\frac{\mu}{\rho}} j_{m, s_{2}}\right)\\

&=&\frac{\mu}{2\mu+\lambda}j_{m,s_1}^5j_{m,s_2}^5 \mathbf{n}^{6} J_{m}^{\prime}\left(j_{m, s_{1}}\right) J_{m}\left(j_{m, s_{1}}\mathbf{n}\right)J_{m}^{\prime}\left(j_{m, s_{2}}\right) J_{m}\left(j_{m, s_{2}}\mathbf{n}\right)\\

&\times& J_m\left(\sqrt{\frac{\mu}{2\mu+\lambda}}j_{m,s_1}\right)
J_m\left(\sqrt{\frac{\mu}{2\mu+\lambda}}\mathbf{n}j_{m,s_1}\right)\\

&\times& \left(\frac{J^{\prime}_m\left(
\sqrt{\frac{\mu}{2\mu+\lambda}}j_{m,s_1}\right)}
{\sqrt{\frac{\mu}{2\mu+\lambda}}j_{m,s_1}
J_m\left(\sqrt{\frac{\mu}{2\mu+\lambda}}j_{m,s_1}\right)}
-\frac{J_m^{\prime}\left(
\sqrt{\frac{\mu}{2\mu+\lambda}}\mathbf{n}j_{m,s_1}\right)}
{\sqrt{\frac{\mu}{2\mu+\lambda}}
\mathbf{n}j_{m,s_1}J_m\left(\sqrt{\frac{\mu}{2\mu+\lambda}}
\mathbf{n}j_{m,s_1}\right)}\right)\\

&\times& J_m\left(\sqrt{\frac{\mu}{2\mu+\lambda}}j_{m,s_2}\right)
J_m\left(\sqrt{\frac{\mu}{2\mu+\lambda}}\mathbf{n}j_{m,s_2}\right)\\

&\times& \left(\frac{J^{\prime}_m\left(
\sqrt{\frac{\mu}{2\mu+\lambda}}j_{m,s_2}\right)}
{\sqrt{\frac{\mu}{2\mu+\lambda}}j_{m,s_2}
J_m\left(\sqrt{\frac{\mu}{2\mu+\lambda}}j_{m,s_2}\right)}
-\frac{J_m^{\prime}\left(
\sqrt{\frac{\mu}{2\mu+\lambda}}\mathbf{n}j_{m,s_2}\right)}
{\sqrt{\frac{\mu}{2\mu+\lambda}}
\mathbf{n}j_{m,s_2}J_m\left(\sqrt{\frac{\mu}{2\mu+\lambda}}
\mathbf{n}j_{m,s_2}\right)}\right)\\
&:=&I_1I_2I_3I_4I_5.
\end{array}
\end{equation}
Combining \eqref{eq:s1_s2}, \eqref{eq:boundj_ms}, and \eqref{eq:bounda_s}, we have
\begin{equation}\label{eq:asympt_jms_case1}
j_{m,s_i}=m+\mathcal{O}\left(m^{\frac{1+2\gamma_i}{3}}\right), \quad \text{for}~i=1,2.
\end{equation}
Therefore, if $m$ is large enough, we can derive the following formulas
\begin{equation*}
\begin{array}{ll}
\sqrt{\frac{\mu}{2\mu+\lambda}}j_{m,s_i}<m, \quad  \sqrt{\frac{\mu}{2\mu+\lambda}}\mathbf{n}j_{m,s_i}<m, \quad \text{for}~i=1,2.
\end{array}
\end{equation*}
Noting that the Bessel function is monotone between the origin and the first extreme point as well as using the above estimates, one can deduce
\begin{equation}\label{EQ:I2I4}
I_2>0,\quad I_4>0.
\end{equation}
With the auxiliary function $f(x)=\frac{1}{x}\frac{J'_{m}(x)}{J_{m}(x)}$ and by following a similar argument in the proof of Theorem 1 in \cite{KRA06}, one can show that
\begin{equation*}
\left(\frac{J'_{m}(x)}{J_{m}(x)}\right)' \leq 0\quad\text{for}\quad x\in\left(0,j_{m,1}\right).
\end{equation*}
This, together with $\frac{1}{x}$ being monotonically decreasing and nonnegative, implies that $f$ is monotonically decreasing and positive in the interval $(0,j_{m,1})$.
Therefore, for $x_i=\sqrt{\frac{\mu}{2\mu+\lambda}}j_{m,s_i}$, $y_i=\sqrt{\frac{\mu}{2\mu+\lambda}}
j_{m,s_i}\mathbf{n}\in (0,j_{m,1})(i=1,2)$, it can be verified that for sufficiently large $m$ we have
\begin{equation}\label{eq:fx}
J_{m}'(x_i)y_iJ_{m}(y_i)- J_{m}'(y_i)x_iJ_{m}(x_i)>0.
\end{equation}
Thus we have
\begin{equation}\label{EQ:I3I5}
I_3>0,\quad I_5>0.
\end{equation}

Without loss of generality, we assume that $J_m^{\prime}(j_{m,s_1})J_m(j_{m,s_1}\mathbf{n})>0$.
Based on \eqref{eq:asympt_Jm_x_large_m} and \eqref{eq:asympt_Jm'_x_large_m}, there exists at least one choice of $s_{2}=\left[m^{\gamma_{2}}\right]$ such that
\begin{equation*}\label{case1}
\begin{array}{lll}
J_m'(j_{m,s_2})J_m(j_{m,s_2}\mathbf{n})&=&
\displaystyle{-\frac{2}{\pi}\left(\frac{j^2_{m,s_2}-m^2}
{j^2_{m,s_2}\mathbf{n}^2-m^2}\right)^{\frac{1}{4}}}
\cos\left(m\mathcal{O}\left(m^{\frac{2\left(\gamma_2-1
\right)}{3}}\right)\right)\\
&\times &\cos \left(m\left(\sqrt{\mathbf{n}^2-1}-\frac{\pi}{2}+
\arcsin\left(\frac{1}{\mathbf{n}}\right)
+\mathcal{O}\left(m^{\frac{2\left(\gamma_2-1\right)}{3}}\right)
\right)\right)\left(1+\mathcal{O}\left(\frac{1}{m}\right)\right)\\
&<&0.
\end{array}
\end{equation*}
The last estimate is based on the fact that the above two cosine functions never have the same frequency, which
together with (\ref{eq:fm_s1fm_s2_case1}), (\ref{EQ:I2I4}), and (\ref{EQ:I3I5}) implies that
\begin{equation}\label{eq:mono less 0}
f_{m}\left(\sqrt{\frac{\mu}{\rho}} j_{m, s_{1}}\right)f_{m}\left(\sqrt{\frac{\mu}{\rho}} j_{m, s_{2}}\right)<0.
\end{equation}

$\mbox{Case 2:}~1<\mathbf{n}=\sqrt{\frac{2\mu+\lambda}{\mu}}$. By similar arguments as in Case 1, we have
\begin{equation*}\label{eq:fm_s1fm_s2_case2}
\begin{array}{lll}
&&f_{m}\left(\sqrt{\frac{\mu}{\rho}} j_{m, s_{1}}\right)f_{m}\left(\sqrt{\frac{\mu}{\rho}} j_{m, s_{2}}\right)\\
&=&j_{m,s_1}^5j_{m,s_2}^5 \mathbf{n}^{4} J'_{m}\left(j_{m, s_{1}}\right) J_{m}\left(j_{m, s_{1}}\mathbf{n}\right)J'_{m}\left(j_{m, s_{2}}\right) J_{m}\left(j_{m, s_{2}}\mathbf{n}\right)\\

&\times& J_m\left(\sqrt{\frac{\mu}{2\mu+\lambda}}j_{m,s_1}\right)
J_m'\left(\sqrt{\frac{\mu}{2\mu+\lambda}}j_{m,s_1}\right)
\times \left(\mathbf{n}\frac{J_m\left(\mathbf{n}
\sqrt{\frac{\mu}{2\mu+\lambda}}j_{m,s_1}\right)}
{J_m\left(\sqrt{\frac{\mu}{2\mu+\lambda}}j_{m,s_1}\right)}
-\frac{J_m'\left(\mathbf{n}
\sqrt{\frac{\mu}{2\mu+\lambda}}j_{m,s_1}\right)}
{J_m'\left(\sqrt{\frac{\mu}{2\mu+\lambda}}
j_{m,s_1}\right)}\right)\\

&\times& J_m\left(\sqrt{\frac{\mu}{2\mu+\lambda}}j_{m,s_2}\right)
J_m'\left(\sqrt{\frac{\mu}{2\mu+\lambda}}j_{m,s_2}\right)\times \left(\mathbf{n}\frac{J_m\left(\mathbf{n}
\sqrt{\frac{\mu}{2\mu+\lambda}}j_{m,s_2}\right)}
{J_m\left(\sqrt{\frac{\mu}{2\mu+\lambda}}j_{m,s_2}\right)}
-\frac{J_m'\left(\mathbf{n}
\sqrt{\frac{\mu}{2\mu+\lambda}}j_{m,s_2}\right)}
{J_m'\left(\sqrt{\frac{\mu}{2\mu+\lambda}}
j_{m,s_2}\right)}\right)\\

&=&j_{m,s_1}^5j_{m,s_2}^5 \mathbf{n}^{4}
J_m\left(\sqrt{\frac{\mu}{2\mu+\lambda}}j_{m,s_1}\right)
J_m\left(\sqrt{\frac{\mu}{2\mu+\lambda}}j_{m,s_2}\right)\\
&\times &J_{m}^{\prime2}\left(j_{m, s_{1}}\right)
J_{m}^{\prime2}\left(j_{m, s_{2}}\right)
J_{m}\left(j_{m, s_{1}}\mathbf{n}\right)
J_{m}\left(j_{m, s_{2}}\mathbf{n}\right).
\end{array}
\end{equation*}
Using the estimate in \eqref{eq:asympt_jms_case1} and the monotonicity of the Bessel function on the interval between {the origin and the first extreme point}, one has that
\begin{equation*}
J_{m}\left(\sqrt{\frac{\mu}{2 \mu+\lambda}} j_{m, s_{1}}\right) J_{m}\left(\sqrt{\frac{\mu}{2 \mu+\lambda}} j_{m, s_{2}}\right)>0.
\end{equation*}
Therefore, the sign of $f_{m}\left(\sqrt{\frac{\mu}{\rho}} j_{m, s_{1}}\right) f_{m}\left(\sqrt{\frac{\mu}{\rho}} j_{m, s_{2}}\right)$ depends only on $J_{m}\left(j_{m, s_{1}} \mathbf{n}\right) J_{m}\left(j_{m, s_{2}} \mathbf{n}\right)$.

By straightforward calculations, one can verify that
\begin{equation*}\label{case2}
\begin{array}{lll}
J_{m}\left(j_{m, s_{1}} \mathbf{n}\right) J_{m}\left(j_{m, s_{2}} \mathbf{n}\right)&=&\displaystyle{\frac{2}{\pi}
\left(\left(j_{m,s_1}^2\mathbf{n}^2-m^2\right)
\left(j_{m,s_2}^2\mathbf{n}^2-m^2\right)\right)^{-\frac{1}{4}}}\\
&\times&\cos \left(m\left(\sqrt{\mathbf{n}^{2}-1}-\frac{\pi}{2}+\arcsin \left(\frac{1}{\mathbf{n}}\right)+\mathcal{O}
\left(m^{\frac{2\left(\gamma_{1}-1\right)}{3}}\right)\right)\right)\\
&\times&\cos \left(m\left(\sqrt{\mathbf{n}^{2}-1}-\frac{\pi}{2}+\arcsin \left(\frac{1}{\mathbf{n}}\right)+\mathcal{O}
\left(m^{\frac{2\left(\gamma_{2}-1\right)}{3}}\right)\right)\right)
\left(1+\mathcal{O}\left(\frac{1}{m}\right)\right)\\
&=&\displaystyle{\frac{2}{\pi}
\left(\left(j_{m,s_1}^2\mathbf{n}^2-m^2\right)
\left(j_{m,s_2}^2\mathbf{n}^2-m^2\right)\right)^{-\frac{1}{4}}}\\
&\times&\displaystyle{\frac{1}{2}}\left(
\cos\left(2m\left(\sqrt{\mathbf{n}^{2}-1}-\frac{\pi}{2}+\arcsin \left(\frac{1}{\mathbf{n}}\right)+\mathcal{O}
\left(m^{\frac{2\left(\gamma_{2}-1\right)}{3}}\right)\right)\right)
\right.\\
&&+\left.\cos\left(m\left(\mathcal{O}
\left(m^{\frac{2\left(\gamma_{2}-1\right)}{3}}\right)\right)\right)
\right)\\
&<&0,
\end{array}
\end{equation*}
which implies \eqref{eq:mono less 0}.

$\mbox{Case 3:}~1<\sqrt{\frac{2\mu+\lambda}{\mu}}<\mathbf{n}$.~
By a similar argument as in \eqref{eq:fm_s1fm_s2_case1}, we have
\begin{equation}\label{eq:fm_s1fm_s2_case3}
\begin{array}{lll}
&&f_{m}\left(\sqrt{\frac{\mu}{\rho}} j_{m, s_{1}}\right)f_{m}\left(\sqrt{\frac{\mu}{\rho}} j_{m, s_{2}}\right)\\
&=&j_{m,s_1}^5j_{m,s_2}^5 \mathbf{n}^{4} J_{m}'\left(j_{m, s_{1}}\right) J_{m}\left(j_{m, s_{1}}\mathbf{n}\right)J_{m}'\left(j_{m, s_{2}}\right) J_{m}\left(j_{m, s_{2}}\mathbf{n}\right)\\

&\times& J_m\left(\sqrt{\frac{\mu}{2\mu+\lambda}}j_{m,s_1}\right)
J_m'\left(\sqrt{\frac{\mu}{2\mu+\lambda}}j_{m,s_1}\right)\times \left(\mathbf{n}\frac{J_m\left(\mathbf{n}
\sqrt{\frac{\mu}{2\mu+\lambda}}j_{m,s_1}\right)}
{J_m\left(\sqrt{\frac{\mu}{2\mu+\lambda}}j_{m,s_1}\right)}
-\frac{J_m'\left(\mathbf{n}
\sqrt{\frac{\mu}{2\mu+\lambda}}j_{m,s_1}\right)}
{J_m'\left(\sqrt{\frac{\mu}{2\mu+\lambda}}
j_{m,s_1}\right)}\right)\\

&\times& J_m\left(\sqrt{\frac{\mu}{2\mu+\lambda}}j_{m,s_2}\right)
J_m'\left(\sqrt{\frac{\mu}{2\mu+\lambda}}j_{m,s_2}\right)\times \left(\mathbf{n}\frac{J_m\left(\mathbf{n}
\sqrt{\frac{\mu}{2\mu+\lambda}}j_{m,s_2}\right)}
{J_m\left(\sqrt{\frac{\mu}{2\mu+\lambda}}j_{m,s_2}\right)}
-\frac{J_m'\left(\mathbf{n}
\sqrt{\frac{\mu}{2\mu+\lambda}}j_{m,s_2}\right)}
{J_m'\left(\sqrt{\frac{\mu}{2\mu+\lambda}}
j_{m,s_2}\right)}\right)\\
&:=&I_6I_7I_8I_9I_{10}.
\end{array}
\end{equation}
By using the recurrence relation of the Bessel function (formula (9.1.27) in \cite{AS72}) and \eqref{eq_Jm_mx}, one has
\begin{equation*}\label{I7I8}
\begin{array}{lll}
I_7I_8&=&\left\{\mathbf{n}\left(J_{m-1}\left(\sqrt{\frac{\mu}
{2\mu+\lambda}}j_{m,s_1}\right)-\sqrt{\frac{\mu}
{2\mu+\lambda}}\left(1-\frac{1}{\mathbf{n}^2}\right)J_m
\left(\sqrt{\frac{\mu}{2\mu+\lambda}}j_{m,s_1}\right)\right)
J_m\left(\mathbf{n}\sqrt{\frac{\mu}{2\mu+\lambda}}
j_{m,s_1}\right)\right.\\

&-&\left.J_m\left(\sqrt{\frac{\mu}{2\mu+\lambda}}
j_{m,s_1}\right)J_{m-1}\left(\mathbf{n}
\sqrt{\frac{\mu}{2\mu+\lambda}}
j_{m,s_1}\right)\right\}\left(1+\mathcal{O}\left(m^{\frac{2\left(
\gamma_1-1\right)}{3}}\right)\right)\\

&=&J_m\left(\sqrt{\frac{\mu}{2\mu+\lambda}}
j_{m,s_1}\right)\left\{\left(\frac{\mathbf{n}}
{\sqrt{\frac{\mu}{2\mu+\lambda}}}-
\sqrt{\frac{\mu}{2\mu+\lambda}}\left(
\mathbf{n}-\frac{1}{\mathbf{n}}\right)\right)
J_m\left(\mathbf{n}\sqrt{\frac{\mu}{2\mu+\lambda}}
j_{m,s_1}\right)\right.\\

&-&\left.J_{m-1}\left(\mathbf{n}
\sqrt{\frac{\mu}{2\mu+\lambda}}
j_{m,s_1}\right)\right\}\left(1+\mathcal{O}\left(m^{\frac{2\left(
\gamma_1-1\right)}{3}}\right)\right)\\

&=&J_m\left(\sqrt{\frac{\mu}{2\mu+\lambda}}
j_{m,s_1}\right)\left(\frac{\mathbf{n}}
{\sqrt{\frac{\mu}{2\mu+\lambda}}}-
\sqrt{\frac{\mu}{2\mu+\lambda}}\left(
\mathbf{n}-\frac{1}{\mathbf{n}}\right)-1\right)\\

&\times&\sqrt{\frac{2}{\pi \sqrt{x^{2}-m^{2}}}} \cos \left(m\left(\frac{\mu}{2\mu+\lambda}
\sqrt{\mathbf{n}^{2}-1}-\frac{\pi}{2}+\arcsin \left(\frac{1}{\mathbf{n}\sqrt{\frac{\mu}
{2\mu+\lambda}}}\right)+\mathcal{O}
\left(m^{\frac{2\left(\gamma_{1}-1\right)}{3}}
\right)\right)\right)\\

&\times&\left(1+\mathcal{O}\left(m^{\frac{2\left(
\gamma_1-1\right)}{3}}\right)\right).
\end{array}
\end{equation*}
Similarly, we have
\begin{equation*}
\begin{array}{lll}
I_9I_{10}&=&J_m\left(\sqrt{\frac{\mu}{2\mu+\lambda}}
j_{m,s_2}\right)\left(\frac{\mathbf{n}}
{\sqrt{\frac{\mu}{2\mu+\lambda}}}-
\sqrt{\frac{\mu}{2\mu+\lambda}}\left(
\mathbf{n}-\frac{1}{\mathbf{n}}\right)-1\right)\\

&\times&\sqrt{\frac{2}{\pi \sqrt{x^{2}-m^{2}}}} \cos \left(m\left(\frac{\mu}{2\mu+\lambda}
\sqrt{\mathbf{n}^{2}-1}-\frac{\pi}{2}+\arcsin \left(\frac{1}{\mathbf{n}\sqrt{\frac{\mu}
{2\mu+\lambda}}}\right)+\mathcal{O}
\left(m^{\frac{2\left(\gamma_{2}-1\right)}{3}}
\right)\right)\right)\\

&\times&\left(1+\mathcal{O}\left(m^{\frac{2\left(
\gamma_2-1\right)}{3}}\right)\right).
\end{array}
\end{equation*}
Using the above two estimates and (\ref{eq:fm_s1fm_s2_case3}), we know that the sign of $f_{m}\left(\sqrt{\frac{\mu}{\rho}} j_{m, s_{1}}\right) f_{m}\left(\sqrt{\frac{\mu}{\rho}} j_{m, s_{2}}\right)$ depends only on the
oscillation part of (\ref{eq:fm_s1fm_s2_case3}), that is
\begin{equation*}\label{sign_case3}
\begin{array}{lll}
&&\left\{\cos \left(m\left(\sqrt{\mathbf{n}^{2}-1}-\frac{\pi}{2}+\arcsin \left(\frac{1}{\mathbf{n}}\right)+\mathcal{O}\left(m^{\frac{2
\left(\gamma_{1}-1\right)}{3}}\right)\right)\right)\right.\\

&\times&\left.\cos \left(m\left(\sqrt{\mathbf{n}^{2}-1}-\frac{\pi}{2}+\arcsin \left(\frac{1}{\mathbf{n}}\right)+\mathcal{O}\left(m^{\frac{2
\left(\gamma_{2}-1\right)}{3}}\right)\right)\right)\right\}\\

&\times&\cos \left(m\left(\mathcal{O}
\left(m^{\frac{2\left(\gamma_{1}-1\right)}{3}}
\right)\right)\right)\cos \left(m\left(\mathcal{O}
\left(m^{\frac{2\left(\gamma_{2}-1\right)}{3}}
\right)\right)\right)\\

&\times&\left\{\cos \left(m\left(\frac{\mu}{2\mu+\lambda}
\sqrt{\mathbf{n}^{2}-1}-\frac{\pi}{2}+\arcsin \left(\frac{1}{\mathbf{n}\sqrt{\frac{\mu}
{2\mu+\lambda}}}\right)+\mathcal{O}
\left(m^{\frac{2\left(\gamma_{1}-1\right)}{3}}
\right)\right)\right)\right.\\

&\times&\left.\cos \left(m\left(\frac{\mu}{2\mu+\lambda}
\sqrt{\mathbf{n}^{2}-1}-\frac{\pi}{2}+\arcsin \left(\frac{1}{\mathbf{n}\sqrt{\frac{\mu}
{2\mu+\lambda}}}\right)+\mathcal{O}
\left(m^{\frac{2\left(\gamma_{2}-1\right)}{3}}
\right)\right)\right)\right\}\\
&:=&I_{11}I_{12}I_{13},
\end{array}
\end{equation*}
where
\begin{equation*}
\begin{array}{lll}
I_{11}&:=&\frac{1}{4}\left\{\cos \left(2m\left(\sqrt{\mathbf{n}^{2}-1}-\frac{\pi}{2}+\arcsin \left(\frac{1}{\mathbf{n}}\right)+\mathcal{O}\left(m^{\frac{2
\left(\gamma_{2}-1\right)}{3}}\right)\right)\right)+\cos \left(m\left(\mathcal{O}
\left(m^{\frac{2\left(\gamma_{2}-1\right)}{3}}
\right)\right)\right)\right\},\\

I_{12}&:=&\cos \left(m\left(\mathcal{O}
\left(m^{\frac{2\left(\gamma_{1}-1\right)}{3}}
\right)\right)\right)\cos \left(m\left(\mathcal{O}
\left(m^{\frac{2\left(\gamma_{2}-1\right)}{3}}
\right)\right)\right),\\

I_{13}&:=&\left\{\cos \left(2m\left(\frac{\mu}{2\mu+\lambda}
\sqrt{\mathbf{n}^{2}-1}-\frac{\pi}{2}+\arcsin \left(\frac{1}{\mathbf{n}\sqrt{\frac{\mu}
{2\mu+\lambda}}}\right)+\mathcal{O}
\left(m^{\frac{2\left(\gamma_{2}-1\right)}{3}}
\right)\right)\right)\right.\\
&+&\left.\cos \left(m\left(\mathcal{O}
\left(m^{\frac{2\left(\gamma_{2}-1\right)}{3}}
\right)\right)\right)\right\}.
\end{array}
\end{equation*}
Since $I_{11}$ and $I_{13}$ never have the same frequency, there exists a $\gamma_2$ such that $I_{11}I_{13}<0$.
Moreover, given $\gamma_2$, one can choose a $\gamma_1$ such that $I_{12}>0$, which implies \eqref{eq:mono less 0}.

In what follows, for given functions $s_1(m),s_2(m)$, there exists a sufficiently large $M$ such that if $m>M$, we have a transmission eigenvalue denoted by $\omega_m(s_1(m),s_2(m))$ satisfying
\begin{equation}\label{eq:mono_eig}
\omega_m(s_1(m),s_2(m)) \in\left(\sqrt{\frac{\mu}{\rho}}j_{m,s_1(m)}, \sqrt{\frac{\mu}{\rho}}j_{m,s_2(m)}\right), \quad m=M+1, M+2, \cdots,
\end{equation}
which implies \eqref{eq:mono_omega_m}.

Part 2. We prove the boundary-localising properties of the
corresponding transmission eigenfunctions $(\mathbf{u}_{m}$,
$\mathbf{v}_{m})$. Let $\Omega_\tau=\left\{\mathbf{x} \in
\mathbb{R}^{2}:|\mathbf{x}|<\tau\right\}$ and $\omega_m$ be the
eigenvalues defined in \eqref{eq:mono_eig}. The associated
eigenfunctions are given by
\begin{equation*}
\begin{array}{lll}
\mathbf{u}_{m}&=&\alpha_m\nabla J_m(\sqrt{\frac{\rho}{2\mu+\lambda}}\omega_m|x|)e^{\mathrm{i}m\theta}
+\gamma_m\nabla\times J_m(\sqrt{\frac{\rho}{\mu}}\omega_m|x|)e^{\mathrm{i}m\theta},\\
\mathbf{v}_{m}&=&\beta_m\nabla J_m(\sqrt{\frac{\tilde{\rho}}{2\mu+\lambda}}\omega_m|x|)e^{\mathrm{i}m\theta}
+\delta_m\nabla\times J_m(\sqrt{\frac{\tilde{\rho}}{\mu}}\omega_m|x|)e^{\mathrm{i}m\theta}.
\end{array}
\end{equation*}

Firstly, for a fixed $\tau$, we next prove that there exists a sufficiently large $m$, such that $\sqrt{\frac{\rho}{2\mu+\lambda}}\omega_m\tau<m$ and $\sqrt{\frac{\rho}{\mu}}\omega_m\tau<m$.
Combining\eqref{eq:boundj_ms}, \eqref{eq:bounda_s}, and \eqref{eq:mono_eig}, one can derive
\begin{equation}\label{eq:estimate_k1}
\begin{array}{ll}
\sqrt{\frac{\rho}{2\mu+\lambda}}\omega_m\tau
<\sqrt{\frac{\rho}{2\mu+\lambda}}
\sqrt{\frac{\mu}{\rho}}(m+o(m))\tau<\sqrt{\frac{\mu}{2\mu+\lambda}}\tau m+o(m)<m,
\end{array}
\end{equation}
and
\begin{equation}\label{eq:estimate_k2}
\begin{array}{ll}
\sqrt{\frac{\rho}{\mu}}\omega_m\tau<\sqrt{\frac{\rho}{\mu}}
\sqrt{\frac{\mu}{\rho}}(m+o(m))\tau<\tau m+o(m)<m.
\end{array}
\end{equation}

Secondly, we shall prove that the transmission eigenfunction $\mathbf{u}_{m}$ are boundary-localized.
By the definition of $\mathbf{u}_{m}$, we have
\begin{equation}\label{eq:um1}
\begin{array}{lll}
\displaystyle{\frac{\|\mathbf{u}_m\|^2_{L^2(\Omega_{\tau})}}
{\|\mathbf{u}_m\|^2_{L^{2}(\Omega)}}}
&=&\frac{\displaystyle{\int_{\Omega_{\tau}}
\mathbf{u}_m^{T}\overline{\mathbf{u}}_m\mathrm{d}S}}
{\displaystyle{\int_{\Omega}
\mathbf{u}_m^{T}\overline{\mathbf{u}}_m\mathrm{d}S}}\\
&=&\frac{\displaystyle{\int_{0}^{2 \pi} \int_{0}^{\tau}\left(\begin{array}{l}
\alpha_m\frac{\partial\phi_m}{\partial x}-\gamma_m\frac{\partial\xi_m}{\partial y} \\
\alpha_m\frac{\partial\phi_m}{\partial y}+\gamma_m\frac{\partial\xi_m}{\partial x}
\end{array}\right)^{T}\left(\begin{array}{l}
\overline{\alpha_m\frac{\partial\phi_m}{\partial x}}-\overline{\gamma_m\frac{\partial\xi_m}{\partial y}} \\
\overline{\alpha_m\frac{\partial\phi_m}{\partial y}}+\overline{\gamma_m\frac{\partial\xi_m}{\partial x}}
\end{array}\right) r \mathrm{d} r \mathrm{d} \theta}}{\displaystyle{\int_{0}^{2 \pi} \int_{0}^{1}\left(\begin{array}{l}
\alpha_m\frac{\partial\phi_m}{\partial x}-\gamma_m\frac{\partial\xi_m}{\partial y} \\
\alpha_m\frac{\partial\phi_m}{\partial y}+\gamma_m\frac{\partial\xi_m}{\partial x}
\end{array}\right)^{T}\left(\begin{array}{l}
\overline{\alpha_m\frac{\partial\phi_m}{\partial x}}-\overline{\gamma_m\frac{\partial\xi_m}{\partial y}} \\
\overline{\alpha_m\frac{\partial\phi_m}{\partial y}}+\overline{\gamma_m\frac{\partial\xi_m}{\partial x}}
\end{array}\right) r \mathrm{d} r \mathrm{d} \theta}}\\
&=&\displaystyle \frac{\int_{0}^{\tau}\left[II_1
-2\Im{(\alpha_m\overline{\gamma_m})}mh'_{1,m}\right]\mathrm{d}r}
{\int_{0}^{1}\left[II_1
-2\Im{(\alpha_m\overline{\gamma_m})}mh'_{1,m}\right]\mathrm{d}r},
\end{array}
\end{equation}
where
\begin{equation*}
\begin{array}{lll}
II_1&=&\displaystyle{
C_{1,m}\left(\left|\frac{\partial\phi_m}{\partial x}\right|^2+\left|\frac{\partial\phi_m}{\partial y}\right|^2\right)r+C_{2,m}\left(\left|\frac{\partial\xi_m}{\partial x}\right|^2+\left|\frac{\partial\xi_m}{\partial y}\right|^2\right)r}\\
&=&\displaystyle{
\sum_{i=1}^{2}C_{i,m}\left(k_{i,m}^2J_{m}^{'2}(k_{i,m}r)r
+m^2\frac{J_{m}^2(k_{i,m}r)}{r}\right)},\\
C_{1,m}&=&|\alpha_m|^2, \quad C_{2,m}=|\gamma_m|^2, \quad h_{1,m}(r)=J_m(k_{1,m}r)J_m(k_{2,m}r),\\
k_{1,m}&=&\sqrt{\frac{\rho}{2\mu+\lambda}}\omega_m, \quad k_{2,m}=\sqrt{\frac{\rho}{\mu}} \omega_{m}.
\end{array}
\end{equation*}
For the term $II_1$, combining \eqref{eq:estimate_k1}, \eqref{eq:estimate_k2}, \eqref{eq:um1} and \eqref{eq:integrabyparts}, together with $k_{1,m}<k_{2,m}<m$ and $k_{i,m}^2J_{m-1}^{2}(k_{i,m}r)r$ is monotonically increasing, we have
\begin{equation*}
\begin{array}{lll}
&&\displaystyle{\frac{\|\mathbf{u}_m\|^2_{L^2(\Omega_{\tau})}}
{\|\mathbf{u}_m\|^2_{L^{2}(\Omega)}}}\\
&\leq& \displaystyle{\frac{m^2[|\alpha_m|^2J^2_{m-1}(k_{1,m}\tau)
+|\gamma_m|^2J_{m-1}(k_{2,m}\tau)]+2m(|\Im(\alpha_m\overline{\gamma_m})|)
J_m(k_{1,m}\tau)J_m(k_{2,m}\tau)}{\displaystyle{m^2\int_{0}^{1}\left[\sum\limits_{i=1}^2C_{i,m}
\frac{J_{m}^2(k_{i,m}r)}{r}\right]\mathrm{d}r}-2m\Im(\alpha_m\overline{\gamma_m}) J_m(k_{1,m})J_m(k_{2,m})}}.
\end{array}
\end{equation*}
Let $g_2(r)=\frac{J_{m}^2(k_{2,m}r)}{r}$, it follows from \eqref{eq:lowerboundofJm'_Jm} that $g_2(r)$ is monotonically increasing. Then if $m$ is large enough, one derives
\begin{equation*}
\begin{array}{lll}
&&\displaystyle{\frac{\|\mathbf{u}_m\|^2_{L^2(\Omega_{\tau})}}
{\|\mathbf{u}_m\|^2_{L^{2}(\Omega)}}}\\
 &\leq&
  \frac{\displaystyle{
  m^2\left(|\alpha_m|^2J^2_{m-1}\left(k_{1,m}\tau\right)
  +|\gamma_m|^2J_{m-1}\left(k_{2,m}\tau\right)\right)
  +2m(|\Im\left(\alpha_m\overline{\gamma_m}\right)|)
  J_m\left(k_{1,m}\tau\right)J_m\left(k_{2,m}\tau\right)}}
  {\displaystyle{m^2\left(1-\frac{2}{\sqrt{m}}\right)C_{2,m}
  \frac{1-\tau}{2}J_m^2\left(k_{2,m}\frac{1+\tau}{2}\right)
  +m^2\left(1-\frac{2}{\sqrt{m}}\right)C_{1,m}
  \frac{1-\tau}{2}J_m^2\left(k_{1,m}\frac{1+\tau}{2}\right)}}\\
\end{array}
\end{equation*}

\begin{equation*}\label{eq:um2}
\begin{array}{ll}
 &\leq
\begin{array}{ll}
 C_{7}(\tau)\left[\frac{\displaystyle{
 |\alpha_m|^2J_{m-1}^2\left(k_{1,m}\tau\right)
 +|\gamma_m|^2J_{m-1}^2\left(k_{2,m}\tau\right)}}
 {\displaystyle{|\alpha_m|^2J_m^2\left(k_{1,m}\frac{1+\tau}{2}\right)
 +|\gamma_m|^2J_m^2\left(k_{2,m}\frac{1+\tau}{2}\right)}}+\right.\\
 \left.\frac{\displaystyle{
 2\left(|\Im\left(\alpha_m\overline{\gamma_m}\right)|\right)
 J_m\left(k_{1,m}\tau\right)J_m\left(k_{2,m}\tau\right)}}
 {\displaystyle{
 m\left(|\alpha_m|^2J_m^2\left(k_{1,m}\frac{1+\tau}{2}\right)
 +|\gamma_m|^2J_m^2\left(k_{2,m}\frac{1+\tau}{2}\right)\right)}}\right]
 \end{array}\\
 &\leq
 C_{7}(\tau)\left[\left(\displaystyle{
 \frac{J_{m-1}\left(k_{1,m}\tau\right)}
 {J_m\left(k_{1,m}\frac{1+\tau}{2}\right)}}\right)^2
 +\left(\displaystyle{
 \frac{J_{m-1}\left(k_{2,m}\tau\right)}
 {J_m\left(k_{2,m}\frac{1+\tau}{2}\right)}}\right)^2
 +\displaystyle{\frac{1}{m}}\right].
\end{array}
\end{equation*}
Similar to the arguments in the proof of Thereom\,2.6 in \cite{DJLZ22}, for any $\tau \in (0,1)$, there exist $\delta_1(\tau,\mathbf{n})$ and $\delta_2(\tau,\mathbf{n})>0$ such that
\begin{equation*}
\displaystyle{\frac{\|\mathbf{u}_m\|^2_{L^2(\Omega_{\tau})}}
{\|\mathbf{u}_m\|^2_{L^{2}(\Omega)}}}\leq C_{7}(\tau)\left(1-\delta_1(\tau,\mathbf{n})\right)^{2m}
+C_{7}(\tau)\left(1-\delta_2(\tau,\mathbf{n})\right)^{2m}
+C_{7}(\tau)\frac{1}{m},
\end{equation*}
which gives
\begin{equation*}
\lim\limits_{m\rightarrow\infty}\displaystyle{\frac{\|\mathbf{u}_m\|^2_{L^2(\Omega_{\tau})}}{\|\mathbf{u}_m\|^2_{L^{2}(\Omega)}}}=0.
\end{equation*}
Hence, the transmission eigenfunction of $\mathbf{u}_m$ is boundary-localized.

Finally, it remains to prove that the eigenfunctions $\mathbf{v}_m$ is not boundary-localized.
A similar argument to \eqref{eq:um1} yields that
\begin{equation*}\label{eq:vm1}
\begin{array}{lll}
&&\displaystyle{\frac{\|\mathbf{v}_m\|^2_{L^2(\Omega_{\tau})}}
{\|\mathbf{v}_m\|^2_{L^{2}(\Omega)}}}
=\frac{\displaystyle{\int_{\Omega_{\tau}}
\mathbf{v}_m^{T}\overline{\mathbf{v}}_m\mathrm{d}S}}
{\displaystyle{\int_{\Omega}
\mathbf{v}_m^{T}\overline{\mathbf{v}}_m\mathrm{d}S}}\\
&=&\frac{\displaystyle{\int_{0}^{2 \pi} \int_{0}^{\tau}\left(\begin{array}{l}
\beta_m\frac{\partial\psi_m}{\partial x}-\delta_m\frac{\partial\zeta_m}{\partial y} \\
\beta_m\frac{\partial\psi_m}{\partial y}+\delta_m\frac{\partial\zeta_m}{\partial x}
\end{array}\right)^{T}\left(\begin{array}{l}
\overline{\beta_m\frac{\partial\psi_m}{\partial x}}-\overline{\delta_m\frac{\partial\zeta_m}{\partial y}} \\
\overline{\beta_m\frac{\partial\psi_m}{\partial y}}+\overline{\delta_m\frac{\partial\zeta_m}{\partial x}}
\end{array}\right) r \mathrm{d} r \mathrm{d} \theta}}{\displaystyle{\int_{0}^{2 \pi} \int_{0}^{1}\left(\begin{array}{l}
\beta_m\frac{\partial\psi_m}{\partial x}-\delta_m\frac{\partial\zeta_m}{\partial y} \\
\beta_m\frac{\partial\psi_m}{\partial y}+\delta_m\frac{\partial\zeta_m}{\partial x}
\end{array}\right)^{T}\left(\begin{array}{l}
\overline{\beta_m\frac{\partial\psi_m}{\partial x}}-\overline{\delta_m\frac{\partial\zeta_m}{\partial y}} \\
\overline{\beta_m\frac{\partial\psi_m}{\partial y}}+\overline{\delta_m\frac{\partial\zeta_m}{\partial x}}
\end{array}\right) r \mathrm{d} r \mathrm{d} \theta}}\\
&=&\frac{\displaystyle{\int_{0}^{\tau}
\left[\widetilde{II}_1-2\Im{(\beta_m\overline{\delta_m})}m
\widetilde{h}'_{1,m}\right] \mathrm{d}r}}
{\displaystyle{\int_{0}^{1}\left[\widetilde{II}_1-2\Im{(\beta_m\overline{\delta_m})}m
\widetilde{h}'_{1,m}\right]\mathrm{d}r}},\\
\end{array}
\end{equation*}
where
\begin{equation*}
\begin{array}{lll}
\widetilde{II}_1&=&\displaystyle{
\widetilde{C}_{1,m}\left(\left|\frac{\partial\psi_m}{\partial x}\right|^2+\left|\frac{\partial\psi_m}{\partial y}\right|^2\right)+\widetilde{C}_{2,m}\left(\left|\frac{\partial\zeta_m}
{\partial x}\right|^2+\left|\frac{\partial\zeta_m}{\partial y}\right|^2\right)}\\
&=&\displaystyle{
\sum_{i=1}^{2}\widetilde{C}_{i,m} \left(\widetilde{k}_{i,m}^2 J_{m}^{'2}(\widetilde{k}_{i,m}r)r+m^2\frac{J_{m}^2(\widetilde{k}_{i,m}r)}{r}\right)},\\
\widetilde{C}_{1,m}&=&|\beta_m|^2, \quad \widetilde{C}_{2,m}=|\delta_m|^2, \quad \widetilde{h}_{1,m}(r)=J_m(\widetilde{k}_{1,m}r)
J_m(\widetilde{k}_{2,m}r),\\
\end{array}
\end{equation*}
\begin{equation*}
\begin{array}{lll}
\displaystyle \widetilde{k}_{1,m}&=&\sqrt{\frac{\widetilde{\rho}}
{2\mu+\lambda}}\omega_m, \quad \widetilde{k}_{2,m}=\sqrt{\frac{\widetilde{\rho}}{\mu}} \omega_{m}.
\end{array}
\end{equation*}
It follows from \eqref{eq:asympt_Jm_x_large_m} and the Cauchy inequality that
\begin{equation}\label{eq:vm2}
\begin{array}{lll}
\displaystyle{\frac{\|\mathbf{v}_m\|^2_{L^2(\Omega_{\tau})}}
{\|\mathbf{v}_m\|^2_{L^{2}(\Omega)}}}
&\geq& \frac{\displaystyle{\int_{0}^{\tau}\left[\widetilde{II}_1
-2\Im{\left(\beta_m\overline{\delta_m}\right)}m
\widetilde{h}'_{1,m}(r)\right]\mathrm{d}r}}
{\displaystyle{\int_{0}^{1}\left[\widetilde{II}_1
-2\Im{\left(\beta_m\overline{\delta_m}\right)}m
\widetilde{h}'_{1,m}(r)\right]\mathrm{d}r}}\\

&\geq&  \frac{\displaystyle\int_{0}^{\tau}\widetilde{II}_1\mathrm{d}r-m\left(\sum_{i=1}^{2}\widetilde{C}_{i,m}
J^2_m\left(\widetilde{k}_{i,m}\tau\right)\right)}
{\displaystyle\int_{0}^{1}\widetilde{II}_1\mathrm{d}r+m\left(\sum_{i=1}^{2}\widetilde{C}_{i,m}
J^2_m\left(\widetilde{k}_{i,m}\right)\right)}\\

&\geq&\frac{\displaystyle{\sum\limits_{i=1}^{2}\widetilde{C}_{i,m}
\left(\int_{0}^{\tau} k_{i, m}^{2} J_{m-1}^{2}\left(k_{i, m} r\right) r \mathrm{~d} r-2m J_{m}^{2}\left(k_{i, m} \tau\right)\right)}}
{\displaystyle{\sum\limits_{i=1}^{2}\widetilde{C}_{i,m}
\left(\int_{0}^{1} k_{i, m}^{2} J_{m-1}^{2}\left(k_{i, m} r\right) r \mathrm{~d} r\right)}},
\end{array}
\end{equation}
where the last inequality follows from \eqref{eq:integrabyparts}.
Suppose that $\tau>\frac{1}{\mathbf{n}}$, then we have
\begin{equation*}
\tau k_{2,m}>j'_{m,1}.
\end{equation*}
It follows from \eqref{eq:integrabyparts2} that
\begin{equation}\label{eq:vm_rterm}
\frac{\displaystyle{mJ_m^2\left(k_{2,m} \tau\right)}}
{\displaystyle{\int_{0}^{R}rJ^2_m\left(r\right)\mathrm{d}r}}
=\mathcal{O}\left(\frac{1}{m}\right).
\end{equation}
Furthermore, from \eqref{eq:asympt_Jm_x_large_m}, we have
\begin{equation}\label{eq:int_m_to_nm_Jm2r}
\begin{aligned}
& \int_{j'_{m,1}}^{k_{2,m}\tau} J_{m}^{2}\left(r\right) r \mathrm{d} r \\
= & \frac{1}{\pi\sqrt{\left(\mathbf{n}
\tau\right)^2-1}}\int_{j'_{m,1}}^{k_{2,m}\tau} \cos ^{2}\left(
\sqrt{r^2-m^{2}}-\frac{m \pi}{2}+m \arcsin \left(\frac{m}{r}\right)-\frac{\pi}{4}\right) \mathrm{d} r
\left(1+\mathcal{O}\left(\frac{1}{m}\right)\right)\\
= &\frac{1}{2\pi\sqrt{\left(\mathbf{n}
\tau\right)^2-1}}m\int_{\frac{j'_{m,1}}{m}}^{\frac{k_{2,m}\tau}{m}}
\left(1+\sin \left(
2m\left(\sqrt{r^2-1}-\frac{\pi}{2}+\arcsin \left(\frac{1}{r}\right)\right)\right)\right)\mathrm{d} r\left(1+\mathcal{O}\left(\frac{1}{m}\right)\right)\\
= &\frac{\left(\mathbf{n}
\tau-1\right) m}{2\pi\sqrt{\left(\mathbf{n}
\tau\right)^2-1}}\left(1+\mathcal{O}\left(\frac{1}{m}\right)\right).
\end{aligned}
\end{equation}
Using \eqref{eq:integrabyparts2} and formula (9.5.16) in \cite{AS72}, we have
\begin{equation}\label{eq:int_0_to_m_Jm2r}
\begin{aligned}
\int_{0}^{j'_{m,1}} J_{m}^{2}\left(r\right) r \mathrm{d} r&=\frac{1}{2}\left(j^{\prime2}_{m,1}-m^2\right)J^2_m\left(j'_{m,1}\right)
=C_6m^{\frac{2}{3}}\left(1+\mathcal{O}\left(\frac{1}{m^{\frac{2}{3}}}\right)\right).
\end{aligned}
\end{equation}
Substituting \eqref{eq:vm_rterm}, \eqref{eq:int_m_to_nm_Jm2r} and \eqref{eq:int_0_to_m_Jm2r} into \eqref{eq:vm2}, one can derive
\begin{equation*}
\begin{array}{lll}
\displaystyle{\frac{\|\mathbf{v}_m\|^2_{L^2(\Omega_{\tau})}}
{\|\mathbf{v}_m\|^2_{L^{2}(\Omega)}}}
&\gtrsim & \displaystyle{\frac{\min\{\frac{(\mathbf{n} \tau-1)}{2 \pi \sqrt{(\mathbf{n} \tau)^{2}-1}},\frac{\sqrt{\mu+\lambda}}{2 \pi (\sqrt{2\mu+\lambda}+\sqrt{\mu})}\}m}
{2C_{8}m^{\frac{2}{3}}\left(1+\mathcal{O}\left(\frac{1}{m^{\frac{2}{3}}}
\right)\right)+\frac{(\mathbf{n} -1) }{2 \pi \sqrt{\mathbf{n}^{2}-1}}m}\left(1+\mathcal{O}\left(\frac{1}{m}\right)\right)}\\
&\gtrsim & C_{9} \displaystyle{\min\{\frac{(\mathbf{n} \tau-1)}{ \sqrt{(\mathbf{n} \tau)^{2}-1}},\frac{\sqrt{\mu+\lambda}}{ (\sqrt{2\mu+\lambda}+\sqrt{\mu})}\}
\frac{\sqrt{\mathbf{n}^{2}-1}}{(\mathbf{n}-1)}
\left(1+\mathcal{O}\left(\frac{1}{m^{\frac{1}{3}}}\right)\right)}.
\end{array}
\end{equation*}
The proof is complete.
\end{proof}


\subsection{Bi-localisation when $N=2$}

We prove the bi-localisation result when $N=2$, namely Theorem~\ref{bilocalized}.

\begin{proof}
Part 1.~It follows from \eqref{eq:fmomega} that
\begin{equation}\label{eq:bi_fm_s0fm_s0+1}
\begin{split}
&f_{m}\left(\sqrt{\frac{\mu}{\rho}} \frac{j_{m, s_{0}}}{\mathbf{n}}\right)f_{m}\left(\sqrt{\frac{\mu}{\rho}} \frac{j_{m, s_{0}+1}}{\mathbf{n}}\right)
=k_{2}^{10} \mathbf{n}^{4} J_{m}'\left(j_{m, s_{0}}\right) J_{m}\left(\frac{j_{m, s_{0}}}{\mathbf{n}}\right)J_{m}'\left(j_{m, s_{0}+1}\right) J_{m}\left(\frac{j_{m, s_{0}+1}}{\mathbf{n}}\right)\\
&\times\left(J_{m}'\left(x_{1}\right) y_{1} J_{m}\left(y_{1}\right)-J_{m}'\left(y_{1}\right) x_{1} J_{m}\left(x_{1}\right)\right)
\times\left(J_{m}'\left(x_{2}\right) y_{2} J_{m}\left(y_{2}\right)-J_{m}'\left(y_{2}\right) x_{2} J_{m}\left(x_{2}\right)\right),
\end{split}
\end{equation}
where $x_1=\sqrt{\frac{\mu}{2\mu+\lambda}}\frac{j_{m,s_0}}{\mathbf{n}}$, $y_1=\sqrt{\frac{\mu}{2\mu+\lambda}}j_{m,s_0}$,
$x_2=\sqrt{\frac{\mu}{2\mu+\lambda}}\frac{j_{m,s_0+1}}{\mathbf{n}}$, and $y_2=\sqrt{\frac{\mu}{2\mu+\lambda}}j_{m,s_0+1}$. \\
We claim that
\begin{equation}\label{eq:s0s01}
f_{m}\left(\sqrt{\frac{\mu}{\rho}} \frac{j_{m, s_{0}}}{\mathbf{n}}\right)f_{m}\left(\sqrt{\frac{\mu}{\rho}} \frac{j_{m, s_{0}+1}}{\mathbf{n}}\right)<0.
\end{equation}
On the one hand, based on \eqref{eq:boundj_ms} and \eqref{eq:bounda_s}, we have
\begin{equation*}\label{eq:asympt_jms}
j_{m,s_0}=m+o(m), \quad j_{m,s_0+1}=m+o(m),
\end{equation*}
and $\frac{j_{m,s_0+1}}{\mathbf{n}}<m< j'_{m,1}$ provided $m$ is sufficiently large.
Furthermore, by the monotonicity of the Bessel function on the interval between the origin and the first extreme point, one deduces that
\begin{equation}\label{eq:monoJm}
0<J_{m}\left(\frac{j_{m,s_0}}{\mathbf{n}}\right)<J_{m}\left(\frac{j_{m,s_0+1}}{\mathbf{n}}\right).
\end{equation}
On the other hand, it follows from \eqref{eq:roots} that
\begin{equation}\label{eq:Jm's0Jm's0+1}
J_{m}'\left(j_{m, s_0}\right)J_{m}'\left(j_{m, s_0+1}\right)<0.
\end{equation}
Substituting \eqref{eq:fx}, \eqref{eq:monoJm} and \eqref{eq:Jm's0Jm's0+1} into \eqref{eq:bi_fm_s0fm_s0+1}, we can obtain (\ref{eq:s0s01}).

For given constant $s_0$, there exists a sufficiently large $M$ such that if $m>M$, we have a transmission eigenvalue denoted by $\omega_m(m,s_0)$ satisfying
\begin{equation}\label{eq:bi_eig}
\omega_m(m,s_0) \in\left(\sqrt{\frac{\mu}{\rho}}\frac{j_{m,s_0}}{\mathbf{n}}, \sqrt{\frac{\mu}{\rho}}\frac{j_{m,s_0+1}}{\mathbf{n}}\right), \quad m=M+1, M+2, \cdots,
\end{equation}
which implies (\ref{eq:bi_omega_m}).

Part 2.~ Let $(\mathbf{u}_m,\mathbf{v}_m)$ be the pair of transmission eigenfunctions associated with ${\omega_m}$ in \eqref{eq:bi_eig}.
By the definition of $\mathbf{u}_{m}$, we have
\begin{equation}\label{eq:bi_um1}
\begin{array}{lll}
\displaystyle{\frac{\|\mathbf{u}_m\|^2_{L^2(\Omega_{\tau})}}
{\|\mathbf{u}_m\|^2_{L^{2}(\Omega)}}}
&=&\frac{\displaystyle{\int_{\Omega_{\tau}}
\mathbf{u}_m^{T}\overline{\mathbf{u}}_m\mathrm{d}S}}
{\displaystyle{\int_{\Omega}
\mathbf{u}_m^{T}\overline{\mathbf{u}}_m\mathrm{d}S}}\\
&=&\frac{\displaystyle{\int_{0}^{2 \pi} \int_{0}^{\tau}\left(\begin{array}{l}
\alpha_m\frac{\partial\phi_m}{\partial x}-\gamma_m\frac{\partial\xi_m}{\partial y} \\
\alpha_m\frac{\partial\phi_m}{\partial y}+\gamma_m\frac{\partial\xi_m}{\partial x}
\end{array}\right)^{T}\left(\begin{array}{l}
\overline{\alpha_m\frac{\partial\phi_m}{\partial x}}-\overline{\gamma_m\frac{\partial\xi_m}{\partial y}} \\
\overline{\alpha_m\frac{\partial\phi_m}{\partial y}}+\overline{\gamma_m\frac{\partial\xi_m}{\partial x}}
\end{array}\right) r \mathrm{d} r \mathrm{d} \theta}}{\displaystyle{\int_{0}^{2 \pi} \int_{0}^{1}\left(\begin{array}{l}
\alpha_m\frac{\partial\phi_m}{\partial x}-\gamma_m\frac{\partial\xi_m}{\partial y} \\
\alpha_m\frac{\partial\phi_m}{\partial y}+\gamma_m\frac{\partial\xi_m}{\partial x}
\end{array}\right)^{T}\left(\begin{array}{l}
\overline{\alpha_m\frac{\partial\phi_m}{\partial x}}-\overline{\gamma_m\frac{\partial\xi_m}{\partial y}} \\
\overline{\alpha_m\frac{\partial\phi_m}{\partial y}}+\overline{\gamma_m\frac{\partial\xi_m}{\partial x}}
\end{array}\right) r \mathrm{d} r \mathrm{d} \theta}}\\
&=&\frac{\displaystyle{\int_{0}^{\tau}
\left[II_2-2\Im{(\alpha_m\overline{\gamma_m})}m
h_{2,m}'\right] \mathrm{d}r}}
{\displaystyle{\int_{0}^{1}\left[II_2-2\Im{(\alpha_m\overline{\gamma_m})}m
h_{2,m}'\right] \mathrm{d}r}},
\end{array}
\end{equation}
where
\begin{equation*}
\begin{array}{lll}
{II}_2&=&\displaystyle{C_{1,m}\left(\left|\frac{\partial\phi_m}{\partial x}\right|^2+\left|\frac{\partial\phi_m}{\partial y}\right|^2\right)+C_{2,m}\left(\left|\frac{\partial\xi_m}{\partial x}\right|^2+\left|\frac{\partial\xi_m}{\partial y}\right|^2\right)}\\
&=&\displaystyle{\sum_{i=3}^{4}C_{i,m}\left(k_{i,m}^2J_{m}^{'2}(k_{i,m}r)r
+m^2\frac{J_{m}^2(k_{i,m}r)}{r}\right)},\\
k_{3,m}&=&\sqrt{\frac{\rho}{2\mu+\lambda}}{\omega_m}, \quad k_{4,m}=\sqrt{\frac{\rho}{\mu}} {\omega_m}, \quad
h_{2,m}(r)=J_m(k_{3,m}r)J_m(k_{4,m}r).
\end{array}
\end{equation*}
For a fixed $\tau$, we have
\begin{equation}\label{eq:estimate_k3}
\begin{aligned}
\sqrt{\frac{\rho}{2\mu+\lambda}}{\omega_m}\tau&<\sqrt{\frac{\rho}{2 \mu+\lambda}} \sqrt{\frac{\mu}{\rho}} \frac{1}{\mathbf{n}}(m+o(m)) \tau
<\sqrt{\frac{\mu}{2 \mu+\lambda}} \frac{\tau}{\mathbf{n}} m+o(m) <m,
\end{aligned}
\end{equation}
and
\begin{equation}\label{eq:estimate_k4}
\begin{aligned}
\sqrt{\frac{\rho}{\mu}}{\omega_m}\tau&<\sqrt{\frac{\rho}{ \mu}} \sqrt{\frac{\mu}{\rho}} \frac{1}{\mathbf{n}}(m+o(m)) \tau
<\frac{\tau}{\mathbf{n}} m+o(m) <m.
\end{aligned}
\end{equation}
Substituting \eqref{eq:integrabyparts}, \eqref{eq:estimate_k3} and \eqref{eq:estimate_k4} into \eqref{eq:bi_um1},
together with $k_{3,m}<k_{4,m}<m$ and $k_{i,m}^2J_{m-1}^{2}(k_{i,m}r)r$ being monotonically increasing, we have
\begin{eqnarray*}
&&\displaystyle{\frac{\|\mathbf{u}_m\|^2_{L^2(\Omega_{\tau})}}
{\|\mathbf{u}_m\|^2_{L^{2}(\Omega)}}}\\
&\leq& \displaystyle{\frac{m^2(|\alpha_m|^2J^2_{m-1}(k_{3,m}\tau)
+|\gamma_m|^2J_{m-1}(k_{4,m}\tau))
+2m(|\Im(\alpha_m\overline{\gamma_m})|)
J_m(k_{3,m}\tau)J_m(k_{4,m}\tau)}
{\displaystyle{m^2\int_{0}^{1}\left[\sum\limits_{i=3}^{4}C_{i,m}
\frac{J_{m}^2(k_{i,m}r)}{r}\right]\mathrm{d}r}
-2m\Im(\alpha_m\overline{\gamma_m}) J_m(k_{3,m})J_m(k_{4,m})}}.
\end{eqnarray*}
Using \eqref{eq:lowerboundofJm'_Jm}, if $m$ is large enough, one can obtain
\begin{equation*}\label{eq:um2}
\begin{array}{lll}
&&\displaystyle{\frac{\|\mathbf{u}_m\|^2_{L^2(\Omega_{\tau})}}
{\|\mathbf{u}_m\|^2_{L^{2}(\Omega)}}}\\
 &\leq&
  \frac{\displaystyle{
  m^2\left(|\alpha_m|^2J^2_{m-1}\left(k_{3,m}\tau\right)
  +|\gamma_m|^2J_{m-1}\left(k_{4,m}\tau\right)\right)
  +2m(|\Im\left(\alpha_m\overline{\gamma_m}\right)|)
  J_m\left(k_{3,m}\tau\right)J_m\left(k_{4,m}\tau\right)}}
  {\displaystyle{m^2\left(1-\frac{2}{\sqrt{m}}\right)C_{2,m}
  \frac{1-\tau}{2}J_m^2\left(k_{4,m}\frac{1+\tau}{2}\right)
  +m^2\left(1-\frac{2}{\sqrt{m}}\right)C_{1,m}
  \frac{1-\tau}{2}J_m^2\left(k_{3,m}\frac{1+\tau}{2}\right)}}\\
 &\leq&
\begin{array}{ll}C_{10}(\tau)
 \left(\frac{\displaystyle{
 |\alpha_m|^2J_{m-1}^2\left(k_{3,m}\tau\right)
 +|\gamma_m|^2J_{m-1}^2\left(k_{4,m}\tau\right)}}
 {\displaystyle{|\alpha_m|^2J_m^2\left(k_{3,m}\frac{1+\tau}{2}\right)
 +|\gamma_m|^2J_m^2\left(k_{4,m}\frac{1+\tau}{2}\right)}}+\right.\\
 \left.\frac{\displaystyle{
 2\left(|\Im\left(\alpha_m\overline{\gamma_m}\right)|\right)
 J_m\left(k_{3,m}\tau\right)J_m\left(k_{4,m}\tau\right)}}
 {\displaystyle{
 m\left(|\alpha_m|^2J_m^2\left(k_{3,m}\frac{1+\tau}{2}\right)
 +|\gamma_m|^2J_m^2\left(k_{4,m}\frac{1+\tau}{2}\right)\right)}}\right)
 \end{array}\\
 &\leq&
 C_{10}(\tau)\left(\left(\displaystyle{
 \frac{J_{m-1}\left(k_{3,m}\tau\right)}
 {J_m\left(k_{3,m}\frac{1+\tau}{2}\right)}}\right)^2
 +\left(\displaystyle{
 \frac{J_{m-1}\left(k_{4,m}\tau\right)}
 {J_m\left(k_{4,m}\frac{1+\tau}{2}\right)}}\right)^2
 +\displaystyle{\frac{1}{m}}\right).
\end{array}
\end{equation*}
Similar to the arguments in the proof of Thereom\,2.6 in \cite{DJLZ22}, for any $\tau \in (0,1)$,  there exists $\delta_3(\tau,\mathbf{n}),\delta_4(\tau,\mathbf{n})>0$ such that
\begin{equation*}
\displaystyle{\frac{\|\mathbf{u}_m\|^2_{L^2(\Omega_{\tau})}}
{\|\mathbf{u}_m\|^2_{L^{2}(\Omega)}}}\leq C_{10}(\tau)\left[\left(1-\delta_3(\tau,\mathbf{n})\right)^{2m}
+\left(1-\delta_4(\tau,\mathbf{n})\right)^{2m}+\frac{1}{m}\right],
\end{equation*}
which implies
\begin{equation*}
\lim\limits_{m\rightarrow\infty}\displaystyle{\frac{\|\mathbf{u}_m\|^2_{L^2(\Omega_{\tau})}}{\|\mathbf{u}_m\|^2_{L^{2}(\Omega)}}}=0.
\end{equation*}
Hence, the transmission eigenfunction of $\mathbf{u}_m$ is boundary-localized.

Part 3.~ The proof of boundary-localization of $\mathbf{v}_m$ is similar to Part 2, except the estimations of \eqref{eq:estimate_k3} and \eqref{eq:estimate_k4}.
Replacing $k_{3,m}$, $k_{4,m}$, and $h_{2, m}(r)$ by $\widetilde{k}_{3,m}=\sqrt{\frac{\widetilde{\rho}}{2 \mu+\lambda}} \widetilde{\omega}_{m}$,
$\widetilde{k}_{4, m}=\sqrt{\frac{\widetilde{\rho}}{\mu}} \widetilde{\omega}_{m}$,
and $\widetilde{h}_{2, m}(r)=J_{m}\left(\widetilde{k}_{3,m} r\right) J_{m}\left(\widetilde{k}_{4, m}r\right)$, respectively, we have
\begin{equation*}\label{estimate_tildek3}
\begin{aligned}
\sqrt{\frac{\widetilde{\rho}}{2\mu+\lambda}}{\omega_m}\tau <\sqrt{\frac{\widetilde{\rho}}{2 \mu+\lambda}} \sqrt{\frac{\mu}{\rho}} \frac{1}{\mathbf{n}}(m+o(m)) \tau
<\sqrt{\frac{\mu}{2 \mu+\lambda}} \tau m+o(m)<m,
\end{aligned}
\end{equation*}
and
\begin{equation*}\label{estimate_tildek4}
\begin{aligned}
\sqrt{\frac{\widetilde{\rho}}{\mu}}{\omega_m}\tau <\sqrt{\frac{\widetilde{\rho}}{ \mu}} \sqrt{\frac{\mu}{\rho}} \frac{1}{\mathbf{n}}(m+o(m)) \tau
<\tau m+o(m) <m.
\end{aligned}
\end{equation*}
By tedious but similar calculations as in Part 2, we can derive
\begin{equation*}
\begin{aligned}
\frac{\left\|\mathbf{v}_{m}\right\|_{L^{2}\left(\Omega_{\tau}\right)}^{2}}
{\left\|\mathbf{v}_{m}\right\|_{L^{2}(\Omega)}^{2}} \leq C_{11}(\tau)\left[\left(\frac{J_{m-1}\left(\widetilde{k}_{3,m}  \tau\right)}{J_{m}\left(\widetilde{k}_{3,m} \frac{1+\tau}{2}\right)}\right)^{2}+\left(\frac{J_{m-1}\left(\widetilde{k}_{4, m} \tau\right)}{J_{m}\left(\widetilde{k}_{4, m} \frac{1+\tau}{2}\right)}\right)^{2}+\frac{1}{m}\right],
\end{aligned}
\end{equation*}
which implies
\begin{equation*}
\lim _{m \rightarrow \infty} \frac{\left\|\mathbf{v}_{m}\right\|_{L^{2}\left(\Omega_{\tau}\right)}^{2}}
{\left\|\mathbf{v}_{m}\right\|_{L^{2}(\Omega)}^{2}}=0.
\end{equation*}
The proof is complete.
\end{proof}


\subsection{Boundary-localisations for compressional and shear waves in $\mathbb{R}^2$}
In the mono-localization case, we next prove the compressional and shear parts of $\mathbf{u}_{m}$ are boundar-localised.
\begin{proof}
We only consider the case for the compressional part, and the shear part can be proved in a similar manner.
By using the definitions, (\ref{eq:integrabyparts}), and the proof of Theorem 2.6 in \cite{DJLZ22}, we have
\begin{equation}\label{eq:ump}
\begin{array}{lll}
&&\frac{\left\|\mathbf{u}^{p}_{m}\right\|_{L^{2}\left(\Omega_{\tau}\right)}^{2}}
{\left\|\mathbf{u}^{p}_{m}\right\|_{L^{2}(\Omega)}^{2}}=\frac{\displaystyle{
\int_{\Omega_{\tau}} |\nabla \times \xi_{m}|^2 \mathrm{~d} S}}{\displaystyle{\int_{\Omega} | \nabla \times \xi_{m} |^2 \mathrm{~d} S}}
=\frac{\displaystyle{\int_{0}^{2 \pi} \int_{0}^{\tau}\left(\left|\frac{\partial \xi_{m}}{\partial x}\right|^{2}+\left|\frac{\partial \xi_{m}}{\partial y}\right|^{2}\right)r \mathrm{~d} r \mathrm{~d} \theta}}
{\displaystyle{\int_{0}^{2 \pi} \int_{0}^{1}\left(\left|\frac{\partial \xi_{m}}{\partial x}\right|^{2}+\left|\frac{\partial \xi_{m}}{\partial y}\right|^{2}\right)r \mathrm{~d} r \mathrm{~d} \theta}}\\
&=&\frac{\displaystyle{\int_{0}^{\tau}\left(k_{2, m}^{2} J_{m}^{\prime 2}\left(k_{2, m} r\right) r+m^{2} \frac{J_{m}^{2}\left(k_{2, m} r\right)}{r}\right) \mathrm{~d} r }}
{\displaystyle{\int_{0}^{1}\left(k_{2, m}^{2} J_{m}^{\prime 2}\left(k_{2, m} r\right) r+m^{2} \frac{J_{m}^{2}\left(k_{2, m} r\right)}{r}\right) \mathrm{~d} r }}=\frac{\displaystyle{\int_{0}^{\tau} k_{2, m}^{2} J_{m-1}^{2}\left(k_{2, m} r\right) r \mathrm{~d} r-m J_{m}^{2}\left(k_{2, m} \tau\right)}}
{\displaystyle{\int_{0}^{1} k_{2, m}^{2} J_{m-1}^{2}\left(k_{2, m} r\right) r \mathrm{~d} r-m J_{m}^{2}\left(k_{2, m}\right)}}\\
&\leq&\frac{\displaystyle{k_{2, m}^{2} J_{m-1}^{2}\left(k_{2, m} \tau\right) \tau^2 }}
{\displaystyle{m^{2}\int_{0}^{1} J_{m}^{2}\left(k_{2, m} r\right)r \mathrm{~d} r }}
\leq\frac{\displaystyle{k_{2, m}^{2} J_{m-1}^{2}\left(k_{2, m} \tau\right) \tau^2 }}
{\displaystyle{m^{2}\frac{\frac{1}{2} J_{m}^{3}\left(k_{2, m}\right)}{J_{m}\left(k_{2, m}\right)+2 k_{2, m} J_{m}'\left(k_{2, m}\right)} }}\\
&\leq&\frac{\displaystyle{2k_{2, m}^{2} J_{m-1}^{2}\left(k_{2, m} \tau\right) \tau^2 }}
{\displaystyle{m^{2} J_{m}^{2}\left(k_{2, m}\right) }}\left(1+2 k_{2, m} \frac{J_{m}'\left(k_{2, m}\right)}{J_{m}\left(k_{2, m}\right)}\right)
\leq\left(\frac{2k_{2,m}^2\tau^2}{m^2}\right)\left(\frac{J_{m-1}\left(k_{2, m} \tau\right)}{J_{m}\left(k_{2, m}\right)}\right)^{2}\left(1+2 k_{2, m} \frac{J_{m}'\left(k_{2, m}\right)}{J_{m}\left(k_{2, m}\right)}\right).
\end{array}
\end{equation}
Substituting \eqref{eq:paris1} and \eqref{eq:lowerboundofJm'_Jm} into \eqref{eq:ump}, through a straightforward calculation, one can obtain
\begin{equation}
\frac{\left\|\mathbf{u}^{p}_{m}\right\|_{L^{2}\left(\Omega_{\tau}\right)}^{2}}
{\left\|\mathbf{u}^{p}_{m}\right\|_{L^{2}(\Omega)}^{2}}\leq
C_{12}(\tau,\mathbf{n})\left(\frac{J_{m}\left(k_{2, m} \tau\right)}{J_{m}\left(k_{2, m}\right)}\right)^{2}.
\end{equation}
Hence, we have
\begin{equation}
\lim _{m \rightarrow \infty} \frac{\left\|\mathbf{u}^{p}_{m}\right\|_{L^{2}\left(\Omega_{\tau}\right)}^{2}}
{\left\|\mathbf{u}^{p}_{m}\right\|_{L^{2}(\Omega)}^{2}}=0.
\end{equation}

Similarly, we can prove the mono-localisation of $\mathbf{u}^{s}_{m}$ as well as $\mathbf{v}^{p}_{m}$ when $\mathbf{n}\sqrt{\frac{\mu}{\lambda+2\mu}}\leq 1$.
On the other hand, the proof of non-localising properties of $\mathbf{v}^{s}_{m}$ is similar to the arguments in the proof of Theorem~\ref{monolocalized}, as well as $\mathbf{v}^{p}_{m}$ when $\mathbf{n}\sqrt{\frac{\mu}{\lambda+2\mu}}>1$.
Finally, the proofs of the bi-localising cases are similar to the mono-localising cases.
\end{proof}


\subsection{Surface resonance and stress concentration}

We present the proof of Theorem~\ref{surface_resonant} in $\mathbb{R}^2$.

\begin{proof}
Suppose {$\|\mathbf{u}_m\|_{L^{2}(\Omega)}=1$}. By direct calculations, one can derive
\begin{equation}\label{eq:Eum}
\hspace{-5mm}\begin{array}{lll}
&& E_{\mathbf{u}_{m},\Sigma}^{2}\\
&=&\displaystyle{\int_{\Sigma}}\left|\alpha_m\right|^2
\left[\left((\lambda+2 \mu) k_{1,m}^{4}+ \frac{4\mu m^{2}}{r^{4}}+\frac{4\mu m^4}{r^4}- \frac{4\mu m^2k_{1,m}^{2}}{r^2}\right)J_m^2\left(k_{1,m}r\right)\right.\\
&+&\left(\frac{4 \mu}{r}\left(k_{1,m}^{2}-\frac{4 m^{2}}{r^{2}}\right)\right)
J_m\left(k_{1,m}r\right)k_{1,m}J_m'\left(k_{1,m}r\right)+\left.\frac{4\mu\left(m^2+1\right)}{r}
k_{1,m}^2J_m^{\prime2}\left(k_{1,m}r\right)\right]\\
&+&|\gamma_m|^2\left[\mu\left(k_{2,m}^{4}+\frac{4 m^2(m^2+1)}{r^4}- \frac{4 m^2k_{2,m}^{2}}{r^2}\right)J_m^2\left(k_{2,m}r\right)\right.\\
&+&\frac{4\mu}{r}\left(k_{2,m}^{2}-\frac{4 m^{2}}{r^{2}}\right)
J_m\left(k_{2,m}r\right)k_{2,m}J_m'\left(k_{2,m}r\right)+\left.\mu\frac{4 (m^{2}+1)}{r^{2}}
k_{2,m}^2J_m^{\prime2}\left(k_{2,m}r\right)\right]\\
&+&2\mu\Im\left(\alpha_m\overline{\gamma_m}\right)\left[
\frac{2m}{r}\left(\left(k_{1,m}J'_m\left(k_{1,m}r\right)
-\frac{J_m\left(k_{1,m}r\right)}{r}\right)
\left(\frac{2k_{2,m}}{r}J_m'\left(k_{2,m}r\right)
+\left(k_{2,m}^2-\frac{2m^2}{r^2}\right)J_m\left(k_{2,m}r\right)\right)\right.\right.\\
&+&\left.\left(k_{2,m}J'_m\left(k_{2,m}r\right)
-\frac{J_m\left(k_{2,m}r\right)}{r}\right)
\left(\frac{2k_{1,m}}{r}J'_m\left(k_{1,m}r\right)
+\left(k_{1,m}^2-\frac{2m^2}{r^2}\right)J_m\left(k_{1,m}r\right)\right)\right)\\
&-&\left.m\left(\frac{(J_m(k_{1,m}r)J_m(k_{2,m}r))'}{r}\right)\right]\mathrm{d}S,\\
\end{array}
\end{equation}
where
\begin{equation}\label{eq:coefficient}
\begin{split}
&\alpha_m=\frac{1}{C_{13}}, \quad \gamma_m=\frac{-\mathrm{i}}{C_{13}}\cdot\frac{\mathbf{n}}{\mathbf{n}^{2}-1} \cdot \frac{k_{1, m}}{m} \cdot \frac{J_{m}\left(k_{1, m}\right)}{J_{m}\left(k_{2,m}\right)} \cdot\left(\mathbf{n} \frac{J_{m}'\left(k_{1, m}\right)}{J_{m}\left(k_{1, m}\right)}-\frac{J_{m}'\left(k_{1, m} \mathbf{n}\right)}{J_{m}\left(k_{1, m} \mathbf{n}\right)}\right),\\
& C^2_{13}\\
&=\frac{J^2_m(k_{1,m})}{\mathbf{n}^4J^2_m(k_{1,m}\mathbf{n})}
\left[\frac{1}{2}\left(k_{1,m} \mathbf{n}\right)^{2}\left(J_{m}^{2}\left(k_{1,m} \mathbf{n}\right)+J_{m-1}^{2}\left(k_{1,m} \mathbf{n}\right)\right)\right.\\
&-\left.m k_{1,m} \mathbf{n} J_{m}\left(k_{1,m} \mathbf{n}\right) J_{m-1}\left(k_{1,m} \mathbf{n}\right)-m J_{m}^{2}\left(k_{1,m} \mathbf{n}\right)\right]\\
&+\left\{\frac{\mathbf{n}^2}{\left(\mathbf{n}^{2}-1\right)^2} \times \frac{k^2_{1, m}}{m^2} \times \frac{J^2_{m}\left(k_{1, m}\right)}{\mathbf{n}^4J^2_{m}\left(k_{2, m}\mathbf{n}\right)} \times \left(\mathbf{n} \frac{J_{m}'\left(k_{1, m}\right)}{J_{m}\left(k_{1, m}\right)}-\frac{J_{m}'\left(k_{1, m} \mathbf{n}\right)}{J_{m}\left(k_{1, m} \mathbf{n}\right)}\right)^2\right.\\
&\times\left[\frac{1}{2}\left(k_{2,m} \mathbf{n}\right)^{2}\left(J_{m}^{2}\left(k_{2,m} \mathbf{n}\right)+J_{m-1}^{2}\left(k_{2,m} \mathbf{n}\right)\right)\right.-\left.\left.m \mathbf{n}k_{2,m}  J_{m}\left(k_{2,m} \mathbf{n}\right) J_{m-1}\left(k_{2,m} \mathbf{n}\right)-m J_{m}^{2}\left(k_{2,m} \mathbf{n}\right)\right]\right\}.
\end{split}
\end{equation}
It is clear that $\Im(\alpha_m\overline{\gamma_m})=0$. Substituting \eqref{eq:estimate_k3}, \eqref{eq:estimate_k4} and \eqref{eq:lowerboundofJm'_Jm}
into \eqref{eq:Eum}, we deduce
\begin{equation}\label{eq:Eum}
\begin{array}{lll}
E_{\mathbf{u}_{m},\Sigma}^{2}&\geq&C_{14}\left(\mathbf{n},\mu,\lambda\right)
\left(m^4+\mathcal{O}\left(m^{\frac{11}{3}}\right)\right)
|\alpha_m|^2\left(\theta_2-\theta_1\right)
\displaystyle{\int_{\tau}^{1}}J_m^2\left(k_{1,m}r\right)r\mathrm{d}r\\
&+&C_{15}\left(\mathbf{n},\mu,\lambda\right)
\left(m^4+\mathcal{O}\left(m^{\frac{11}{3}}\right)\right)
|\gamma_m|^2\left(\theta_2-\theta_1\right)
\displaystyle{\int_{\tau}^{1}}J_m^2\left(k_{2,m}r\right)r\mathrm{d}r,
\end{array}
\end{equation}
where $C_{14}=\frac{2\mu^2}{2\mu+\lambda}\left((\frac{1}{\mathbf{n}}-1)^2+\frac{3\mu+2\lambda}{\mu}\right)>0$,
$C_{15}=2\mu\left((\frac{1}{\mathbf{n}}-1)^2+1\right)>0$, and
$\theta_{i}(i=1,2)$ are defined in (\ref{eq:theta12}). {It follows from \eqref{eq:coefficient}, \eqref{eq:paris1}, \eqref{eq:lowerboundofJm'_Jm}, and
\eqref{eq:integrabyparts2} that for $m$ sufficiently large, there exists a constant $C_{14}>0$ such
that
\begin{equation}\label{eq:integral estimate}
\displaystyle{\frac{\int_{\tau_2}^{\tau_1}J_m^2\left(mr\right)r\mathrm{d}r}{J_m^2\left(m\tau_1\right)}}\geq C_{16}\left\{\begin{array}{ll}
m^{-1}, ~\mbox{if}~\tau_2<\tau_1<1,\\
m^{-\frac{2}{3}}, ~\mbox{if}~ \tau_2<\tau_1=1+\mathcal{O}(m^{-\frac{2}{3}}).
\end{array}
\right.
\end{equation}
Substituting \eqref{eq:integral estimate} into \eqref{eq:Eum}, we arrive at
\begin{equation*}
\begin{array}{lll}
E_{\mathbf{u}_{m},\Sigma}^{2}\geq C_{1}(\lambda,m,\mathbf{n},\tau,\Sigma) \mu.
\end{array}
\end{equation*}
Similarly, for any $\lambda$ and $m$ sufficiently large, one can derive
\begin{equation*}
\begin{split}
&E_{\mathbf{v}_{m},\Sigma}^{2}\geq C_{2}(\lambda,m,\mathbf{n},\tau,\Sigma) \mu ,\quad E_{\mathbf{u}_{m},\Sigma}^{2}\geq C_{3}(\lambda,\mu,\mathbf{n},\tau,\Sigma) m^{3},\quad
E_{\mathbf{v}_{m},\Sigma}^{2}\geq C_{4}(\lambda,\mu,\mathbf{n},\tau,\Sigma) m^{\frac{10}{3}},\\
&\hspace*{.6cm}\left|\nabla
\mathbf{u}_m \right|^2_{\infty}\geq C_{5}(\lambda,\mu,\mathbf{n},\tau)m^3,\quad \left|\nabla \mathbf{v}_m \right|^2_{\infty}\geq C_{6}(\lambda,\mu,\mathbf{n},\tau)m^{\frac{10}{3}}.
\end{split}
\end{equation*}}
\end{proof}


\section{Proofs of main theorems in three dimensions}\label{sect:3}


We deal with the main theorems in three dimensions. We shall only sketch the necessary modifications compared to the two-dimensional treatments in the previous section.


\subsection{Preliminaries}

Using the Fourier series, the solutions of the system \eqref{eq:Elas1} have the following forms:
\begin{equation}
\begin{split}
\mathbf{u}_{m}^{n}&=a_m^n\nabla j_m(k_1r)Y_m^n+b_m^n\nabla\times \left(xj_m(k_2r)Y_m^n\right)+c_m^n\nabla\times\nabla\times \left(xj_m(k_2r)Y_m^n\right),\\
\mathbf{v}_{m}^{n}&=d_m^n\nabla j_m(\widetilde{k}_1r)Y_m^n+e_m^n\nabla\times \left(xj_m(\widetilde{k}_2r)Y_m^n\right)+f_m^n\nabla\times\nabla\times \left(xj_m(\widetilde{k}_2r)Y_m^n\right).
\end{split}
\end{equation}
By using the boundary condition in \eqref{eq:Elas1}, we see that $\omega$ is a transmission eigenvalue if
\begin{equation}
F_{m}^{n}(\omega;\theta,\phi):=\det(A)=0,
\end{equation}
where the entries of $A\in \mathbb{R}^{6\times6}$ are given by
\begin{equation*}
\begin{array}{llll}
A_{1,1}&=k_1j'_m\left(k_1\right)Y_m^n,
& A_{2,1}&=j_m\left(k_1\right)\partial_{\theta}Y_m^n,\\
A_{1,2}&=m\left(m+1\right)j_m\left(k_2\right)Y_m^n,
& A_{2,2}&=\left(j_m\left(k_2\right)
+k_2j'_m\left(k_2\right)\right)\partial_{\theta}Y_m^n,\\
A_{1,3}&=0,
& A_{2,3}&=\frac{1}{\sin \theta}j_m\left(k_2\right)\partial_{\phi}Y_m^n,\\
A_{1,4}&=\widetilde{k}_1j'_m(\widetilde{k}_1)Y_m^n,
& A_{2,4}&=j_m(\widetilde{k}_1)\partial_{\theta}Y_m^n,\\
A_{1,5}&=m\left(m+1\right)j_m(\widetilde{k}_2)Y_m^n,
& A_{2,5}&=\left(j_m(\widetilde{k}_2)
+\widetilde{k}_2j'_m(\widetilde{k}_2)\right)\partial_{\theta}Y_m^n,\\
A_{1,6}&=0,
& A_{2,6}&=\frac{1}{\sin \theta}j_m(\widetilde{k}_2)\partial_{\phi}Y_m^n,\\
A_{3,1}&=\frac{1}{\sin \theta}j_m\left(k_1\right)\partial_{\phi}Y_m^n,
& A_{4,1}&=\left(2\mu k_1^2j_{m}^{\prime\prime}\left(
k_1\right)-\lambda k_1^2j_m\left(k_1\right)\right)Y_m^n,\\
A_{3,2}&=\frac{1}{\sin \theta}\left(j_m\left(k_2\right)
+k_2j'_m\left(k_2\right)\right)\partial_{\phi}Y_m^n,
& A_{4,2}&=-2\mu m\left(m+1\right)\left( j_{m}\left(
k_2\right)-k_2j'_m\left(k_2\right)\right)Y_m^n,\\
A_{3,3}&=-j_m\left(k_2\right)\partial_{\theta}Y_m^n,
& A_{4,3}&=0,\\
A_{3,4}&=\frac{1}{\sin \theta}j_m(\widetilde{k}_1)\partial_{\phi}Y_m^n,
& A_{4,4}&=\left(2\mu \widetilde{k}_1^2j_{m}^{''}(
\widetilde{k}_1)-\lambda \widetilde{k}_1^2j_m(\widetilde{k}_1)\right)Y_m^n, \\
A_{3,5}&=\frac{1}{\sin \theta}\left(j_m(\widetilde{k}_2)
+\widetilde{k}_2j'_m(\widetilde{k}_2)\right)\partial_{\phi}Y_m^n,
& A_{4,5}&=-2\mu m\left(m+1\right)\left(j_{m}\left(
\widetilde{k}_2\right)-\widetilde{k}_2j'_m(\widetilde{k}_2)\right)Y_m^n, \\
A_{3,6}&=-j_m(\widetilde{k}_2)\partial_{\theta}Y_m^n,
& A_{4,6}&=0,
\end{array}
\end{equation*}
\vspace{-2mm}
\begin{equation*}
\begin{array}{lllll}
A_{5,1}&=-2\mu \left( j_{m}\left(
k_1\right)-k_1j'_m\left(k_1\right)\right)\partial_{\theta}Y_m^n,\\
A_{5,2}&=2\mu\left(-\left(j_m\left(k_2\right)+k_2j'_m\left(k_2\right)\right)\partial_{\theta}Y_m^n+
j_m\left(k_2\right)\left(\frac{1}{\sin^2 \theta}\partial_{\theta}Y_m^n+
m(m+1)\partial_{\theta}Y_m^n+\frac{2\cot \theta }{\sin^2 \theta}\partial^2_{\phi}Y_m^n\right)
\right)\\
&+\mu j_m\left(k_2\right)\left(\frac{1}{\sin^2 \theta}\partial_{\theta}Y_m^n+
k_2^2\partial_{\theta}Y_m^n+\frac{2\cot \theta }{\sin^2 \theta}\partial^2_{\phi}Y_m^n\right),\\
A_{5,3}&=-\frac{\mu}{\sin \theta}\left(2j_{m}\left(
k_2\right)\partial_{\phi}Y_m^n+\left( j_{m}\left(k_2\right)+k_2j'_m\left(k_2\right)\right)\partial_{\theta}Y_m^n\right),\\
A_{5,4}&=-2\mu \left( j_{m}(\widetilde{k}_1)-\widetilde{k}_1j'_m(\widetilde{k}_1)\right)\partial_{\theta}Y_m^n,\\
A_{5,5}&=2\mu\left[-\left(j_m(\widetilde{k}_2)+\widetilde{k}_2j'_m(\widetilde{k}_2)\right)\partial_{\theta}Y_m^n+
j_m(\widetilde{k}_2)\left(\frac{1}{\sin^2 \theta}\partial_{\theta}Y_m^n+
m(m+1)\partial_{\theta}Y_m^n+\frac{2\cot \theta }{\sin^2 \theta}\partial^2_{\phi}Y_m^n\right)
\right]\\
&+\mu j_m(\widetilde{k}_2)\left(\frac{1}{\sin^2 \theta}\partial_{\theta}Y_m^n+
\widetilde{k}_2^2\partial_{\theta}Y_m^n+\frac{2\cot \theta }{\sin^2 \theta}\partial^2_{\phi}Y_m^n\right),\\
A_{5,6}&=-\frac{\mu}{\sin \theta}\left(2j_{m}(
\widetilde{k}_2)\partial_{\phi}Y_m^n+\left( j_{m}(
\widetilde{k}_2)+\widetilde{k}_2j'_m(\widetilde{k}_2)\right)\partial_{\theta}Y_m^n\right),\\
A_{6,1}&=-\frac{2\mu}{\sin \theta} \left( j_{m}\left(
k_1\right)-k_1j'_m\left(k_1\right)\right)\partial_{\phi}Y_m^n,\\
A_{6,2}&=\frac{\mu}{\sin \theta}\left(2m\left(m+1\right)j_{m}\left(
k_2\right)-\left(j_{m}\left(
k_2\right)+k_2j'_m\left(k_2\right)\right)+k_2^2j_{m}\left(
k_2\right)\right)\partial_{\phi}Y_m^n,\\
A_{6,3}&=\mu\left(2j_{m}\left(
k_1\right)+j_{m}\left(
k_2\right)+k_2j'_m\left(k_2\right)\right)\partial_{\theta}Y_m^n,\\
A_{6,4}&=-\frac{2\mu}{\sin \theta} \left( j_{m}(\widetilde{k}_1)-\widetilde{k}_1j'_m(\widetilde{k}_1)\right)\partial_{\phi}Y_m^n,\\
A_{6,5}&=\frac{\mu}{\sin \theta}\left(2m\left(m+1\right)j_{m}(\widetilde{k}_2)-\left(j_{m}(\widetilde{k}_2)+\widetilde{k}_2
j'_m(\widetilde{k}_2)\right)+\widetilde{k}_2^2j_{m}(\widetilde{k}_2)\right)\partial_{\phi}Y_m^n,\\
A_{6,6}&=\mu\left(2j_{m}(\widetilde{k}_2)+j_{m}(\widetilde{k}_2)+\widetilde{k}_2j'_m(\widetilde{k}_2)\right)\partial_{\theta}Y_m^n.\\
\end{array}
\end{equation*}
Here, $j_m(\cdot)$ is the spherical Bessel function and $Y_m^n$ is the spherical harmonics for $m\in\mathbb{N}\cup\{0\}$ and $n=-m,\ldots,m$. By direct calculations, we can derive
\begin{equation*}
\begin{array}{lll}
F_m^n(\omega;\theta,\phi)&=&-\frac{3\rho\omega^2\mu^2}{\sin^4 \theta}
\left(Y_m^n\right)^2\left(\partial_{\phi}Y_m^n\right)^2
\left(\partial_{\theta}Y_m^n\right)^2\left(1-2n^2\cot\theta
\frac{Y_m^n}{\partial_{\theta}Y_m^n}\right)\\
&\times&\left(\sin\theta
\frac{\partial_{\theta}Y_m^n}{\partial_{\phi}Y_m^n}-\frac{1}{\sin\theta}
\frac{\partial_{\phi}Y_m^n}{\partial_{\theta}Y_m^n}\right)
\left(k_1j'_m\left(k_1\right)j_m\left(k_1\mathbf{n}\right)
-k_1\mathbf{n}j_m\left(k_1\right)j'_m\left(k_1\mathbf{n}\right)\right)\\
&\times&\left(k_2j'_m\left(k_2\right)j_m\left(k_2\mathbf{n}\right)
-k_2\mathbf{n}j_m\left(k_2\right)j'_m\left(k_2\mathbf{n}\right)\right)\\
&\times&\left(\frac{\frac{1}{\sin\theta}+\sin\theta
\frac{\partial_{\theta}Y_m^n}{\partial_{\phi}Y_m^n}}
{\frac{3}{\sin^2\theta}\left(1-2n^2\cot\theta
\frac{Y_m^n}{\partial_{\theta}Y_m^n}\right)}
\left|\begin{array}{ll}
j_m\left(k_2\right)+k_2j'_m\left(k_2\right) & j_m\left(k_2\mathbf{n}\right)+k_2\mathbf{n}j'_m\left(k_2\mathbf{n}\right) \\
k_2^2j_m\left(k_2\right) & \left(k_2\mathbf{n}\right)^2j_m\left(k_2\mathbf{n}\right)
\end{array}\right|\right.\\
&&\left.-\sin\theta
\frac{\partial_{\theta}Y_m^n}{\partial_{\phi}Y_m^n}
\left|\begin{array}{ll}
k_2j'_m\left(k_2\right) & k_2\mathbf{n}j'_m\left(k_2\mathbf{n}\right) \\
j_m\left(k_2\right) & j_m\left(k_2\mathbf{n}\right)
\end{array}\right|
\right).
\end{array}
\end{equation*}


\subsection{Mono-localisation when $N=3$}

Now, we are ready to prove Theorem~\ref{monolocalized} in case of $N=3$.

\begin{proof}
Part 1.~  Using Lemma 2.2 in \cite{DLWW22}, we can show that there exists at least one zero point $\omega_m$ of $F^{n}_{m}(\omega)$ in $\left(\sqrt{\frac{\mu}{\rho}} j_{m+1/2, s_{1}}, \sqrt{\frac{\mu}{\rho}} j_{m+1/2, s_{2}}\right)$,
where $s_{1}(m)$ and $s_{2}(m)$ have the following asymptotic formula:
\begin{equation*}
s_{1}(m):=\left[m^{\gamma_{1}}\right], \quad s_{2}(m)=\left[m^{\gamma_{2}}\right], \quad 0<\gamma_{1}<\gamma_{2}<1.
\end{equation*}
Part 2.~  We next prove the boundary-localising properties of the
corresponding transmission eigenfunctions $\left(\mathbf{u}^{n}_{m},
\mathbf{v}^{n}_{m}\right)$. Suppose $\mathbf{n}\tau>1$. Here, we only prove that $\mathbf{v}^{n}_{m}$ is not boundary-localized, and the boundary-localisation of $\mathbf{u}^n_m$ can be proved by following a similar argument to the 2D case. By the definition of
$\mathbf{v}^{n}_{m}$, we have
\begin{equation}\label{III1}
\frac{\displaystyle{\left\|\mathbf{v}^{n}_{m}\right\|_{L^{2}\left(\Omega_{\tau}\right)}^{2}}}
{\displaystyle{\left\|\mathbf{v}^{n}_{m}\right\|_{L^{2}(\Omega)}^{2}}}
=\frac{\displaystyle{\int_{\Omega_{\tau}} (\mathbf{v}^{n}_{m})^T \overline{\mathbf{v}}^{n}_{m} \mathrm{~d} V}}
{\displaystyle{\int_{\Omega} (\mathbf{v}^{n}_{m})^T \overline{\mathbf{v}}^{n}_{m} \mathrm{~d} V}}
=\frac{\displaystyle{\int_{\Omega_{\tau}} III_1\mathrm{~d} V}}{\displaystyle{\int_{\Omega} III_1 \mathrm{~d} V}},
\end{equation}
where
\begin{equation*}
\begin{split}
III_1:=&\left|d_m^n k_1\mathbf{n}j'_{m}\left(k_1\mathbf{n}r\right)+f_m^n \frac{m\left(m+1\right)}{r}j_{m}\left(k_2\mathbf{n}r\right)\right|^2
Y_m^n\overline{Y_m^n}+\frac{1}{\sin^2\theta}\left|e_m^nj_{m}\left(k_2\mathbf{n}r\right)\right|^2
\partial_{\phi}Y_m^n\overline{\partial_{\phi}Y_m^n}\\
+&\frac{1}{r^2\sin^2\theta}\left|d_m^n j_{m}\left(k_1\mathbf{n}r\right)+f_m^n \left(j_{m}\left(k_2\mathbf{n}r\right)
+rk_2\mathbf{n}j'_{m}\left(k_2\mathbf{n}r\right)\right)\right|^2
\partial_{\phi}Y_m^n\overline{\partial_{\phi}Y_m^n}.
\end{split}
\end{equation*}
Substituting \eqref{eq:integrabyparts2}, \eqref{eq:shfexpansion}, \eqref{eq:vm_rterm}, and \eqref{eq:int_0_to_m_Jm2r}  into \eqref{III1}, we can derive
\begin{equation*}
\begin{array}{lll}
\frac{\displaystyle{\left\|\mathbf{v}^{n}_{m}\right\|_{L^{2}\left(\Omega_{\tau}\right)}^{2}}}
{\displaystyle{\left\|\mathbf{v}^{n}_{m}\right\|_{L^{2}(\Omega)}^{2}}}
&=&
\frac{\displaystyle{\int_{0}^{\tau}III_2 \mathrm{~d} r}}
{\displaystyle{\int_{0}^{1}III_2\mathrm{~d} r}}\\
&\geq& C_{8}\begin{cases}\frac{\sqrt{\left(\mathbf{n}\tau-1\right)}}
{\sqrt{\left(\mathbf{n}\tau\right)^2-1}}\frac{\sqrt{\mathbf{n}^2-1}}
{\sqrt{\mathbf{n}-1}}, & \text { if } \sqrt{\frac{\mu}{2\mu+\lambda}}\mathbf{n}<1, \\
\min\{\frac{\sqrt{\left(\mathbf{n}\tau-1\right)}}
{\sqrt{\left(\mathbf{n}\tau\right)^2-1}}\frac{\sqrt{\mathbf{n}^2-1}}
{\sqrt{\mathbf{n}-1}},\frac{\sqrt{\left(\mathbf{n}\tau\sqrt{\frac{\mu}{2\mu+\lambda}}-1\right)}}
{\sqrt{\left(\mathbf{n}\tau\sqrt{\frac{\mu}{2\mu+\lambda}}\right)^2-1}}\frac{\sqrt{\mathbf{n}^2-1}}
{\sqrt{\mathbf{n}-1}}\},
& \text { if } \sqrt{\frac{\mu}{2\mu+\lambda}}\mathbf{n}\geq1, \end{cases}
\end{array}
\end{equation*}
where
\begin{equation*}
\begin{split}
III_2:=&
2\pi r^2\left|d_m^n k_1\mathbf{n}j'_{m}\left(k_1\mathbf{n}r\right)+f_m^n \frac{m\left(m+1\right)}{r}j_{m}\left(k_2\mathbf{n}r\right)\right|^2
+\left(2m+1\right)m\pi r^2\left|e_m^nj_{m}\left(k_2\mathbf{n}r\right)\right|^2\\
+&
\left(2m+1\right)m\pi\left|d_m^n j_{m}\left(k_1\mathbf{n}r\right)+f_m^n \left(j_{m}\left(k_2\mathbf{n}r\right)
+rk_2\mathbf{n}j'_{m}\left(k_2\mathbf{n}r\right)\right)\right|^2.
\end{split}
\end{equation*}
The proof is complete.
\end{proof}


\subsection{Bi-localized modes when $N=3$}

In this subsection, we will give the proof of theorem~\ref{bilocalized} in case of $N=3$.

\begin{proof}
Let
\begin{equation*}
\widetilde{F}^{n}_{m}(\omega)=k_2j'_m\left(k_2\right)j_m\left(k_2\mathbf{n}\right)
-k_2\mathbf{n}j_m\left(k_2\right)j'_m\left(k_2\mathbf{n}\right).
\end{equation*}
Then we have
\begin{equation*}
\widetilde{F}^{n}_{m}(j_{m+1/2,s_0})\widetilde{F}^{n}_{m}(j_{m+1/2,s_0+1})<0.
\end{equation*}
By IVT, there exists at least one zero point $\omega_m$ of $\widetilde{F}^{n}_{m}(\omega)$ in
$$
\left(\sqrt{\frac{\mu}{\rho}} \frac{j_{m+1/2, s_{0}}}{\mathbf{n}}, \sqrt{\frac{\mu}{\rho}} \frac{j_{m+1/2, s_{0}+1}}{\mathbf{n}}\right),
$$
which implies $\widetilde{F}^{n}_{m}(\omega_m;\theta,\phi)=0$.

Next, we show $\mathbf{u}^{n}_{m}$ is boundar-localized.
By the definition of $\mathbf{u}^{n}_{m}$ and similar arguments for (\ref{III1}), we have
\begin{equation*}
\begin{array}{lll}
\frac{\displaystyle{\left\|\mathbf{u}^{n}_{m}\right\|_{L^{2}\left(\Omega_{\tau}\right)}^{2}}}
{\displaystyle{\left\|\mathbf{u}^{n}_{m}\right\|_{L^{2}(\Omega)}^{2}}}
&=&\frac{\displaystyle{\int_{\Omega_{\tau}} (\mathbf{u}_{m}^{n})^{T} \overline{\mathbf{u}}^{n}_{m} \mathrm{~d} V}}
{\displaystyle{\int_{\Omega} (\mathbf{u}_{m}^{n})^{T} \overline{\mathbf{u}}^{n}_{m} \mathrm{~d} V}}
\leq
\frac{\displaystyle{\int_{0}^{\tau}j_m^{\prime2}(k_1r)}}
{\displaystyle{\int_{0}^{1}j_m^{\prime2}(k_1r)}}+
\frac{\displaystyle{\int_{0}^{\tau}j_m^{2}(k_2r)}}
{\displaystyle{\int_{0}^{1}j_m^{2}(k_2r)}}+
\frac{\displaystyle{\int_{0}^{\tau}j_m^{\prime2}(k_2r)}}
{\displaystyle{\int_{0}^{1}j_m^{\prime2}(k_2r)}}\\
&+&
\frac{\displaystyle{\int_{0}^{\tau}\frac{1}{r}j'_m(k_1r)j_m(k_2r)}}
{\displaystyle{\int_{0}^{1}\frac{1}{r}j'_m(k_1r)j_m(k_2r)}}+
\frac{\displaystyle{\int_{0}^{\tau}\frac{1}{r}j_m(k_1r)j'_m(k_2r)}}
{\displaystyle{\int_{0}^{1}\frac{1}{r}j_m(k_1r)j'_m(k_2r)}}+
\frac{\displaystyle{\int_{0}^{\tau}\frac{1}{r}j'_m(k_2r)j_m(k_2r)}}
{\displaystyle{\int_{0}^{1}\frac{1}{r}j'_m(k_2r)j_m(k_2r)}}\\
&+&
\frac{\displaystyle{\int_{0}^{\tau}\frac{1}{r^2}j^2_n(k_1r)}}
{\displaystyle{\int_{0}^{1}\frac{1}{r^2}j^2_n(k_1r)}}+
\frac{\displaystyle{2\int_{0}^{\tau}\frac{1}{r^2}j^2_n(k_2r)}}
{\displaystyle{\int_{0}^{1}\frac{1}{r^2}j^2_n(k_2r)}}+
\frac{\displaystyle{\int_{0}^{\tau}\frac{1}{r^2}j_m(k_1r)j_m(k_2r)}}
{\displaystyle{\int_{0}^{1}\frac{1}{r^2}j_m(k_1r)j_m(k_2r)}},
\end{array}
\end{equation*}
where the last inequality is based on $\frac{a+b}{c+d}\leq \frac{a}{c}+\frac{b}{d}$ with $a,b,c,d>0$.
Using an argument similar to that in the proof of Theorem~\ref{bilocalized} when $N=2$, we have
\begin{equation*}
\lim _{n \rightarrow \infty} \frac{\left\|\mathbf{u}^{n}_{m}\right\|_{L^{2}\left(\Omega_{\tau}\right)}^{2}}
{\left\|\mathbf{u}^{n}_{m}\right\|_{L^{2}(\Omega)}^{2}}=0.
\end{equation*}
In a similar manner, we can prove $\mathbf{v}^{n}_{m}$ is boundary-localized.
\end{proof}


\subsection{Boundary-localisations for compressional and shear wave in $\mathbb{R}^3$}

%
We only consider the compressional wave, and the shear wave can be proved in a similar manner.
By using the definitions, (\ref{eq:integrabyparts}), and the proof of Theorem 2.6 in \cite{DJLZ22}, we have
\begin{align*}
&\frac{\left\|\mathbf{u}^{p}_{m}\right\|_{L^{2}\left(\Omega_{\tau}\right)}^{2}}
{\left\|\mathbf{u}^{p}_{m}\right\|_{L^{2}(\Omega)}^{2}}=\frac{\displaystyle{
\int_{\Omega_{\tau}} |\nabla \times \xi_{m}|^2 \mathrm{~d} S}}{\displaystyle{\int_{\Omega} | \nabla \times \xi_{m} |^2 \mathrm{~d} S}}
=\frac{\displaystyle{\int_{0}^{2 \pi} \int_{0}^{\tau}\left(\left|\frac{\partial \xi_{m}}{\partial x}\right|^{2}+\left|\frac{\partial \xi_{m}}{\partial y}\right|^{2}\right)r \mathrm{~d} r \mathrm{~d} \theta}}
{\displaystyle{\int_{0}^{2 \pi} \int_{0}^{1}\left(\left|\frac{\partial \xi_{m}}{\partial x}\right|^{2}+\left|\frac{\partial \xi_{m}}{\partial y}\right|^{2}\right)r \mathrm{~d} r \mathrm{~d} \theta}}\\
=&\frac{\displaystyle{\int_{0}^{\tau}\left(k_{2, m}^{2} J_{m}^{\prime 2}\left(k_{2, m} r\right) r+m^{2} \frac{J_{m}^{2}\left(k_{2, m} r\right)}{r}\right) \mathrm{~d} r }}
{\displaystyle{\int_{0}^{1}\left(k_{2, m}^{2} J_{m}^{\prime 2}\left(k_{2, m} r\right) r+m^{2} \frac{J_{m}^{2}\left(k_{2, m} r\right)}{r}\right) \mathrm{~d} r }}=\frac{\displaystyle{\int_{0}^{\tau} k_{2, m}^{2} J_{m-1}^{2}\left(k_{2, m} r\right) r \mathrm{~d} r-m J_{m}^{2}\left(k_{2, m} \tau\right)}}
{\displaystyle{\int_{0}^{1} k_{2, m}^{2} J_{m-1}^{2}\left(k_{2, m} r\right) r \mathrm{~d} r-m J_{m}^{2}\left(k_{2, m}\right)}}\\
\leq&\frac{\displaystyle{k_{2, m}^{2} J_{m-1}^{2}\left(k_{2, m} \tau\right) \tau^2 }}
{\displaystyle{m^{2}\int_{0}^{1} J_{m}^{2}\left(k_{2, m} r\right)r \mathrm{~d} r }}
\leq\frac{\displaystyle{k_{2, m}^{2} J_{m-1}^{2}\left(k_{2, m} \tau\right) \tau^2 }}
{\displaystyle{m^{2}\frac{\frac{1}{2} J_{m}^{3}\left(k_{2, m}\right)}{J_{m}\left(k_{2, m}\right)+2 k_{2, m} J_{m}'\left(k_{2, m}\right)} }}\\
\leq&\frac{\displaystyle{2k_{2, m}^{2} J_{m-1}^{2}\left(k_{2, m} \tau\right) \tau^2 }}
{\displaystyle{m^{2} J_{m}^{2}\left(k_{2, m}\right) }}\left(1+2 k_{2, m} \frac{J_{m}'\left(k_{2, m}\right)}{J_{m}\left(k_{2, m}\right)}\right)\\
\leq&\left(\frac{2k_{2,m}^2\tau^2}{m^2}\right)\left(\frac{J_{m-1}\left(k_{2, m} \tau\right)}{J_{m}\left(k_{2, m}\right)}\right)^{2}\left(1+2 k_{2, m} \frac{J_{m}'\left(k_{2, m}\right)}{J_{m}\left(k_{2, m}\right)}\right).
\end{align*}
Substituting \eqref{eq:paris1} and \eqref{eq:lowerboundofJm'_Jm} into \eqref{eq:ump}, along with straightforward calculations, one can obtain
\begin{equation}
\frac{\left\|\mathbf{u}^{p}_{m}\right\|_{L^{2}\left(\Omega_{\tau}\right)}^{2}}
{\left\|\mathbf{u}^{p}_{m}\right\|_{L^{2}(\Omega)}^{2}}\lesssim
C_8(\tau,\mathbf{n})\left(\frac{J_{m}\left(k_{2, m} \tau\right)}{J_{m}\left(k_{2, m}\right)}\right)^{2}.
\end{equation}
Hence, we have
\begin{equation}
\lim _{m \rightarrow \infty} \frac{\left\|\mathbf{u}^{p}_{m}\right\|_{L^{2}\left(\Omega_{\tau}\right)}^{2}}
{\left\|\mathbf{u}^{p}_{m}\right\|_{L^{2}(\Omega)}^{2}}=0.
\end{equation}

Similarly, we can prove the mono-localising properties of $\mathbf{u}^{s}_{m}$ as well as $\mathbf{v}^{p}_{m}$ when $\mathbf{n}\sqrt{\frac{\mu}{\lambda+2\mu}}\leq 1$.
On the other hand, the proof of the non-localisation of $\mathbf{v}^{s}_{m}$ is similar to that of Theorem~\ref{monolocalized}, as well as $\mathbf{v}^{p}_{m}$ when $\mathbf{n}\sqrt{\frac{\mu}{\lambda+2\mu}}>1$.
Finally, the proofs of the bi-localising cases are similar to the mono-localising cases.

\section{Numerics and discussions}\label{sect:numerics}

In this section, we present several representative numerical examples to corroborate our theoretical findings in the previous sections. All these simulations are implemented in {MATLAB 2021b, and run on a workstation with a 2.9 GHz Intel(R) Xeon(R) Platinum 8268 CPU and 2 TB RAM}.

The tested domains include circle, triangle, square, and kite-2d in two dimensions, and ball, kite-3d, and cuboid in three dimensions as follows
\begin{align*}
&\mbox{Circle}&&\{(x,y)|x^2+y^2\leq 1\},\\
&\mbox{Triangle}&&\{(x,y)|\mbox{Triangle with vertices~} (-1,0), (1,0),(0,\sqrt{3})\},\\
&\mbox{Square}&&\{(x,y)|0\leq x\leq 1, 0\leq y\leq 1\},\\
&\mbox{Kite-2d}&&\{(x,y)|\mbox{Kite with boundary~} s_1(t),~0\leq t\leq 2\pi\},\\
&\mbox{Ball}&&\{(x,y,z)|x^2+y^2+z^2\leq 1\},\\
&\mbox{Kite-3d}&&\{(x,y,z)|\mbox{Kite with boundary~} s_2(t,\theta),~0\leq t\leq 2\pi,~0\leq \theta\leq 2\pi\},\\
&\mbox{Cuboid}&&\{(x,y,z)|0\leq x\leq 1, 0\leq y\leq 0.2, 0\leq z\leq 0.2\},
\end{align*}
where
\begin{align*}
s_1(t)&=(\frac{2}{3}\cos{t}+0.6\cos{2t}-0.6,\sin{t}),\\
s_2(t,\theta)&=(\frac{2}{3}\cos{t}+0.6\cos{2t}-0.6,\sin{t}\sin{\theta},\sin{t}\cos{\theta}).
\end{align*}
In the 2D case, standard triangular meshes with mesh length about $0.01$ are generated for all domains and cubic Lagrangian finite element are used to discrete problem (\ref{eq:Elas1}). For the 3D case,  standard tetrahedron meshes with mesh length about $0.025$ are generated for all domains and quadratic Lagrangian finite element are used to discretize problem (\ref{eq:Elas1}).

\begin{figure}[htbp]
  \centering
  \includegraphics[width=0.32\textwidth]{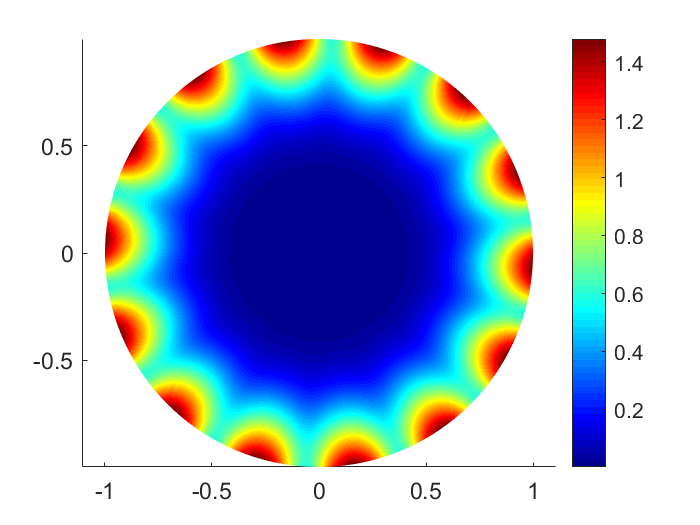}
  \includegraphics[width=0.32\textwidth]{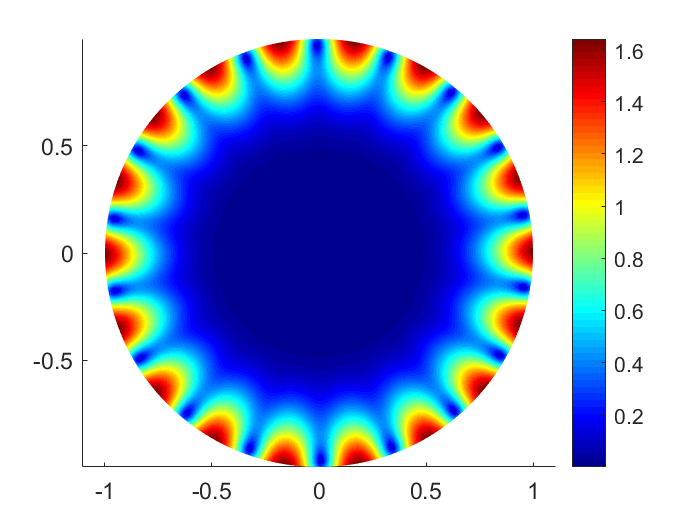}
  \includegraphics[width=0.32\textwidth]{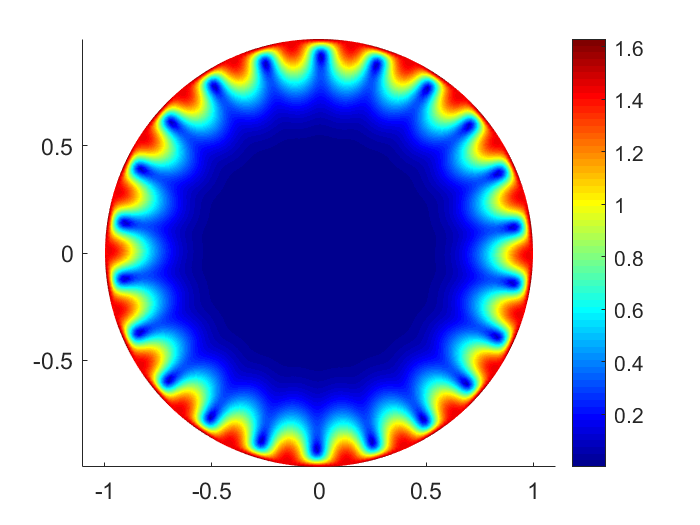}\\
  \includegraphics[width=0.32\textwidth]{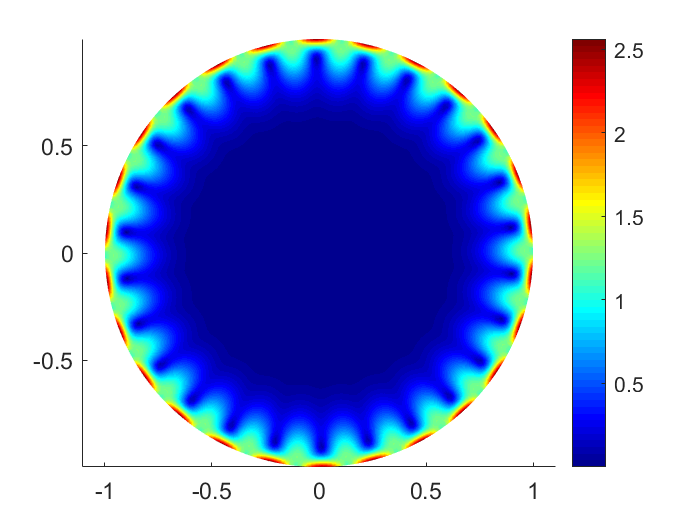}
  \includegraphics[width=0.32\textwidth]{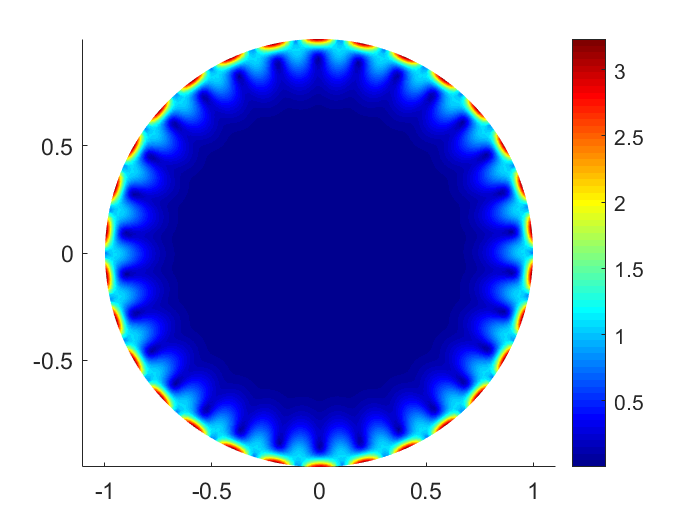}
  \includegraphics[width=0.32\textwidth]{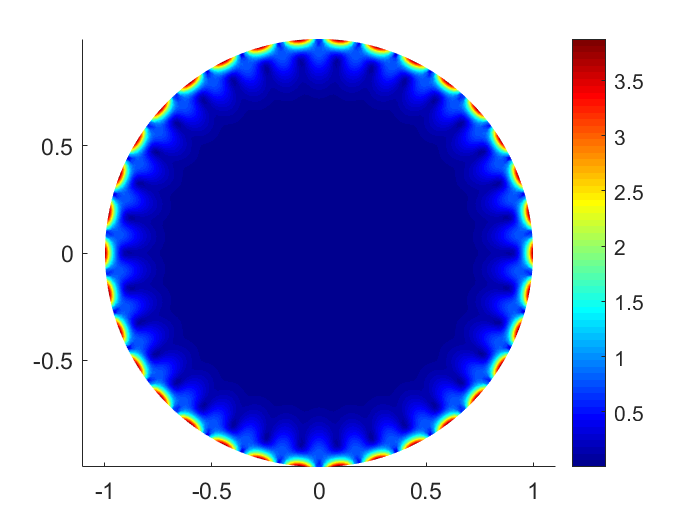}\\
\caption{The eigenfunctions $|\mathbf{u}_m|$ for mono-localized modes with
$\omega_m=5.71$, $6.56$, $7.28$, $7.88$, $8.46$, and $9.05$.}
 \label{fig:radial mono u}
\end{figure}

\begin{figure}[htbp]
  \centering
  \includegraphics[width=0.32\textwidth]{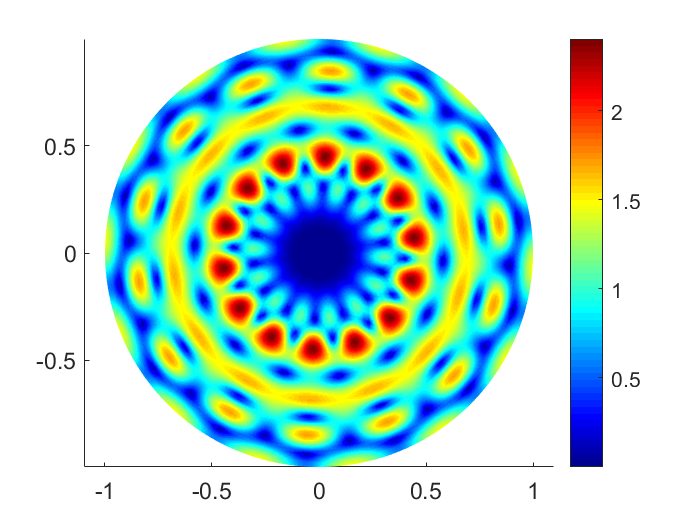}
  \includegraphics[width=0.32\textwidth]{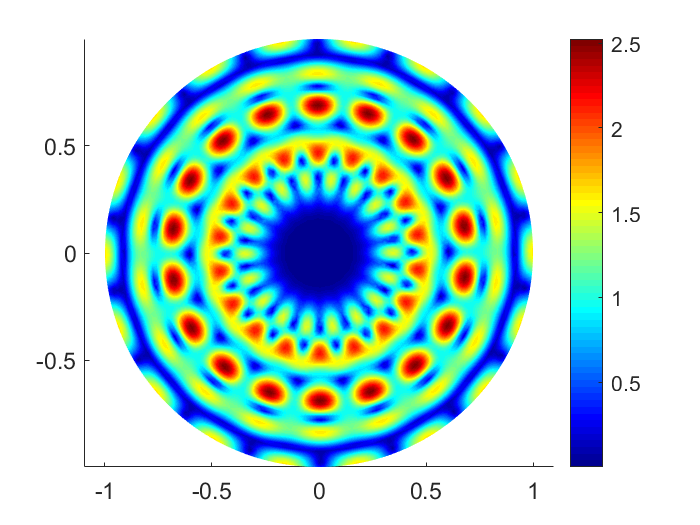}
  \includegraphics[width=0.32\textwidth]{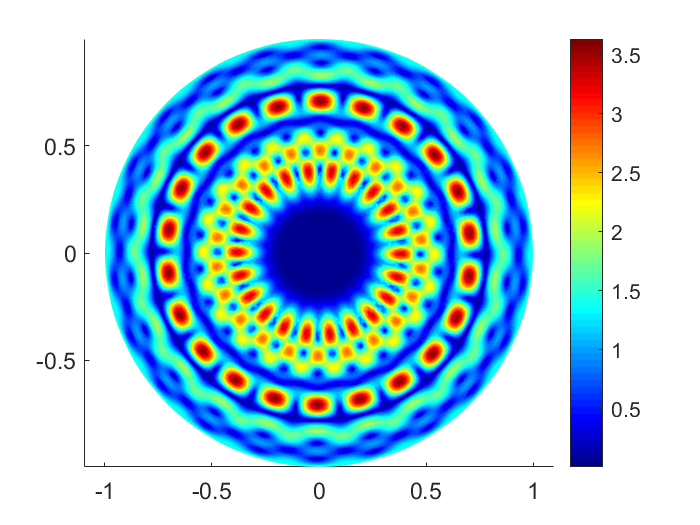}\\
  \includegraphics[width=0.32\textwidth]{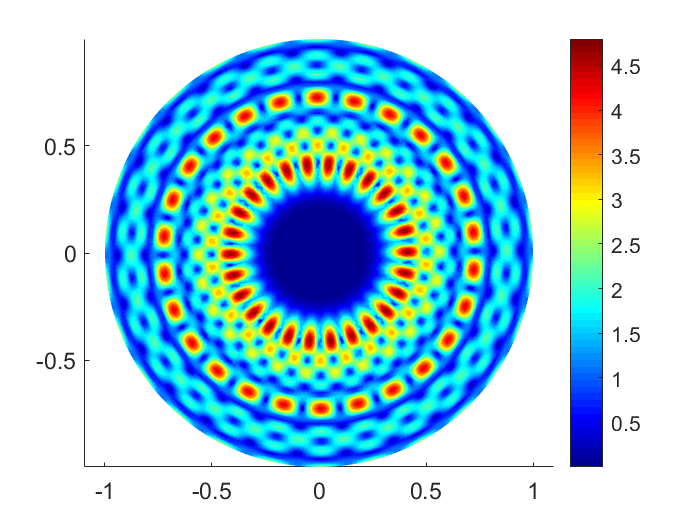}
  \includegraphics[width=0.32\textwidth]{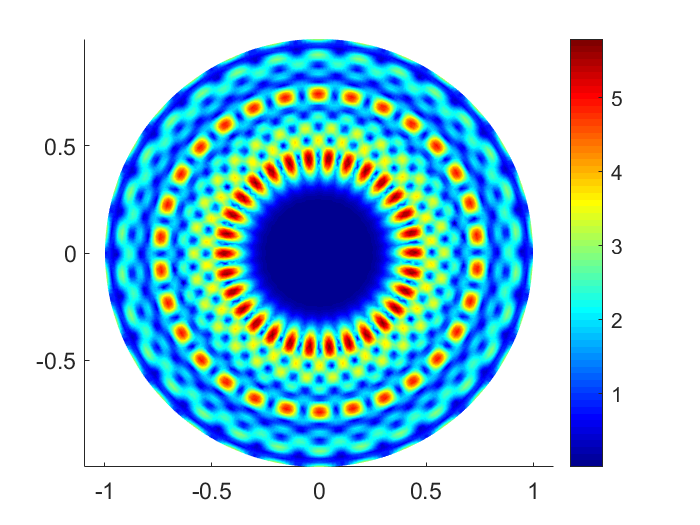}
  \includegraphics[width=0.32\textwidth]{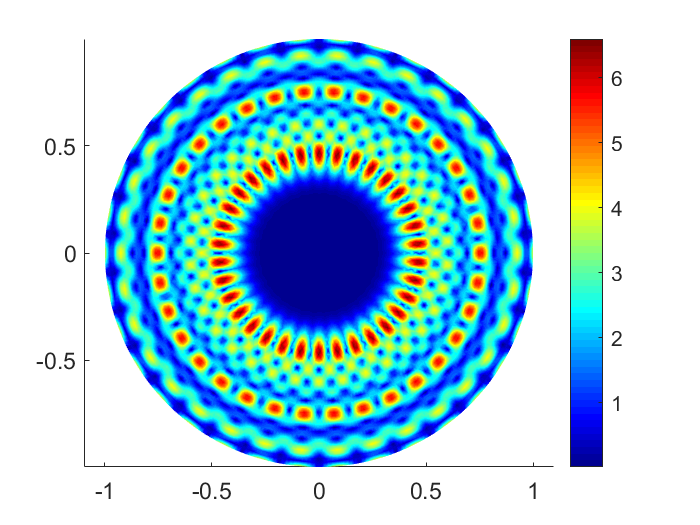}\\
\caption{The eigenfunctions $|\mathbf{v}_m|$ for mono-localized modes with $\omega_m=5.71$, $6.56$, $7.28$, $7.88$, $8.46$, and $9.05$.}
 \label{fig:radial mono v}
\end{figure}

\begin{figure}[htbp]
  \centering
  \includegraphics[width=0.32\textwidth]{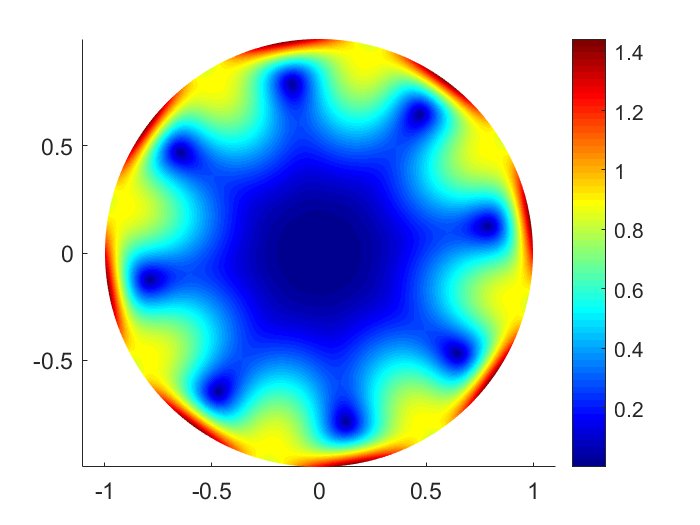}
  \includegraphics[width=0.32\textwidth]{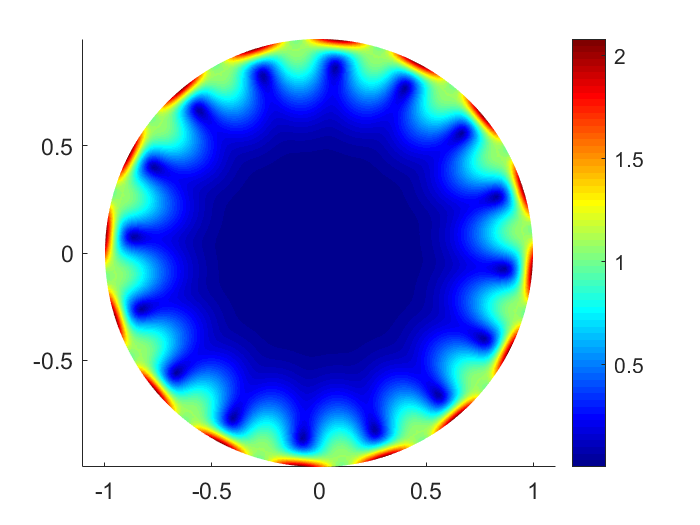}
  \includegraphics[width=0.32\textwidth]{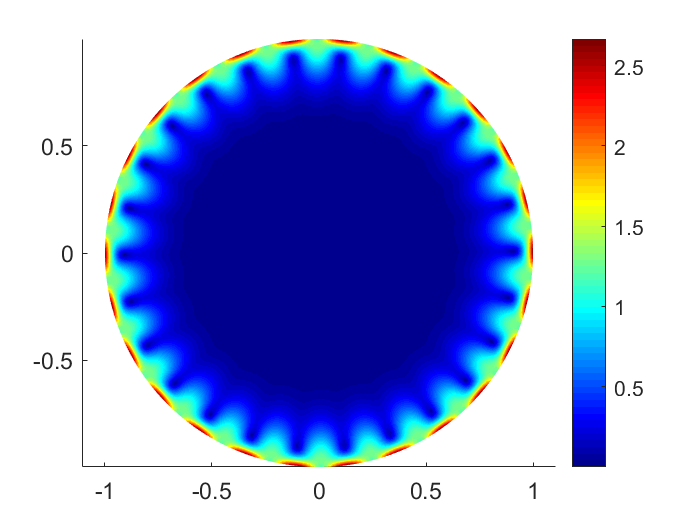}\\
  \includegraphics[width=0.32\textwidth]{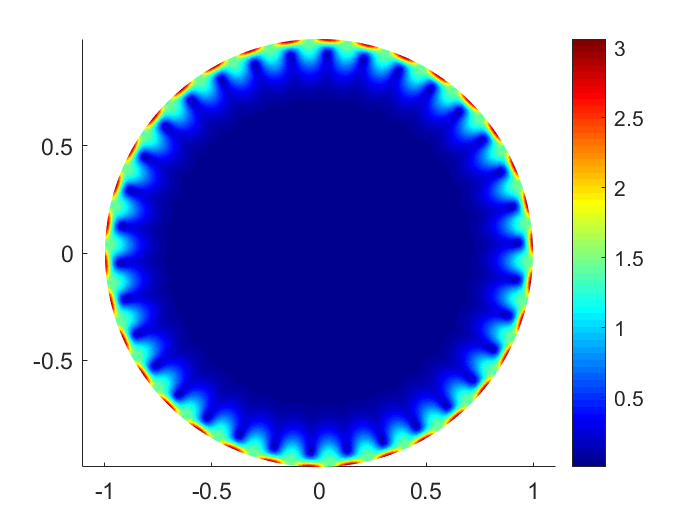}
  \includegraphics[width=0.32\textwidth]{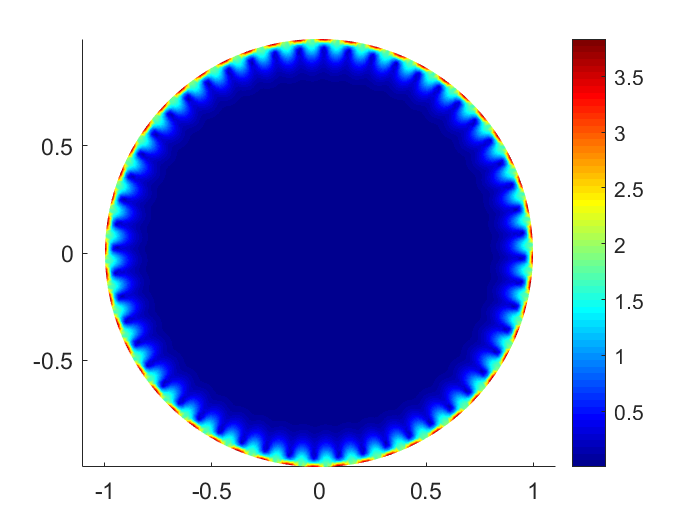}
  \includegraphics[width=0.32\textwidth]{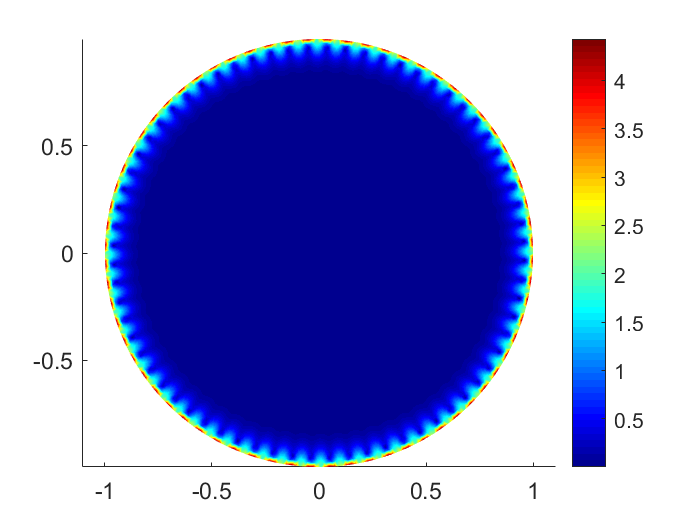}\\
\caption{The eigenfunctions $|\mathbf{u}_m|$ for bi-localized modes with
$\omega_m=2.19$, $3.21$, $4.46$, $5.43$, $7.58$, and $9.48$.}
 \label{fig:radial bi u}
\end{figure}

\begin{figure}[htbp]
  \centering
  \includegraphics[width=0.32\textwidth]{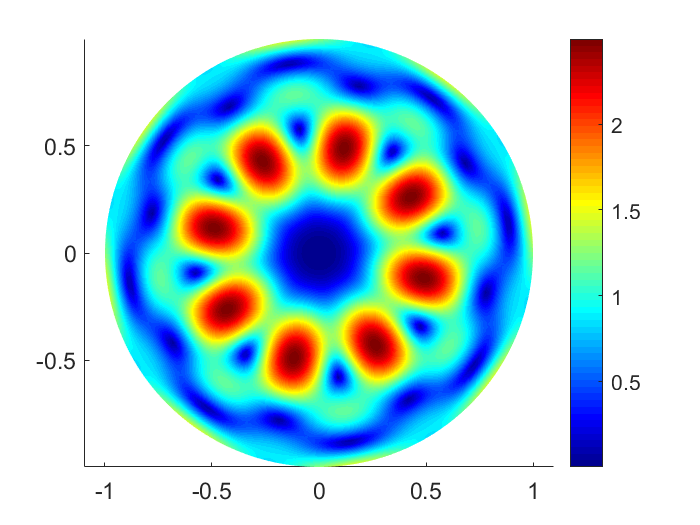}
  \includegraphics[width=0.32\textwidth]{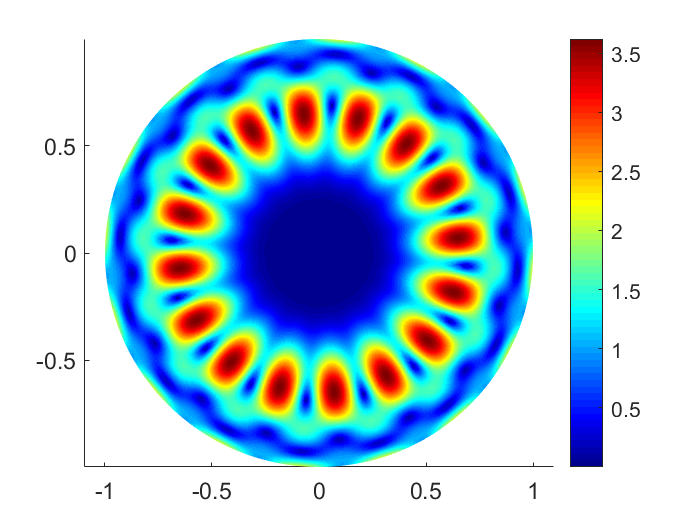}
  \includegraphics[width=0.32\textwidth]{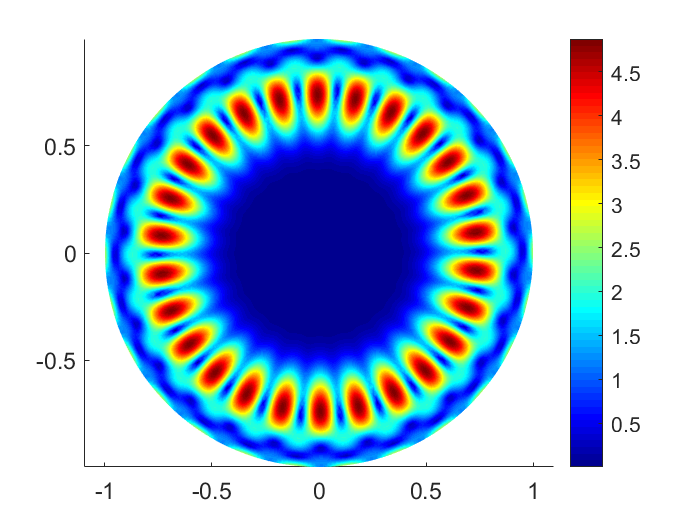}\\
  \includegraphics[width=0.32\textwidth]{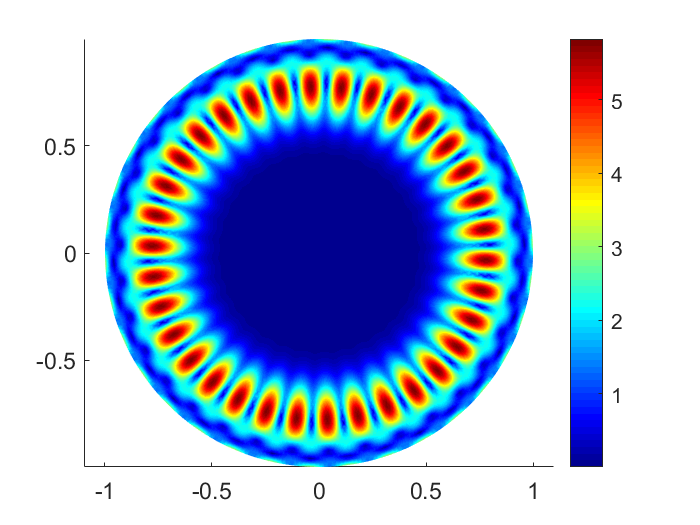}
  \includegraphics[width=0.32\textwidth]{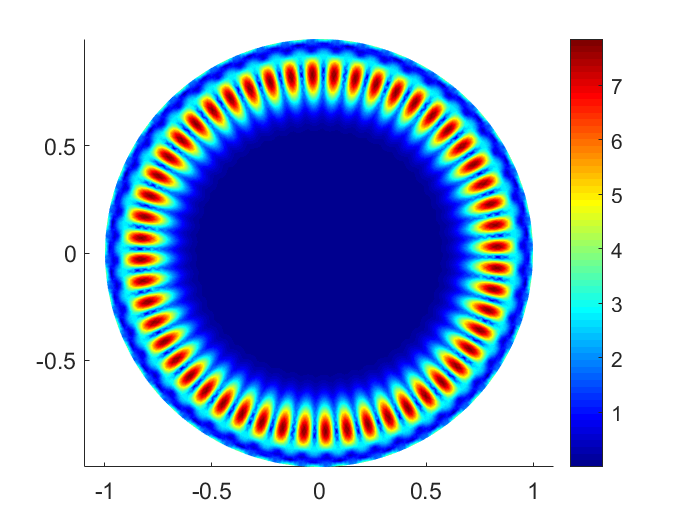}
  \includegraphics[width=0.32\textwidth]{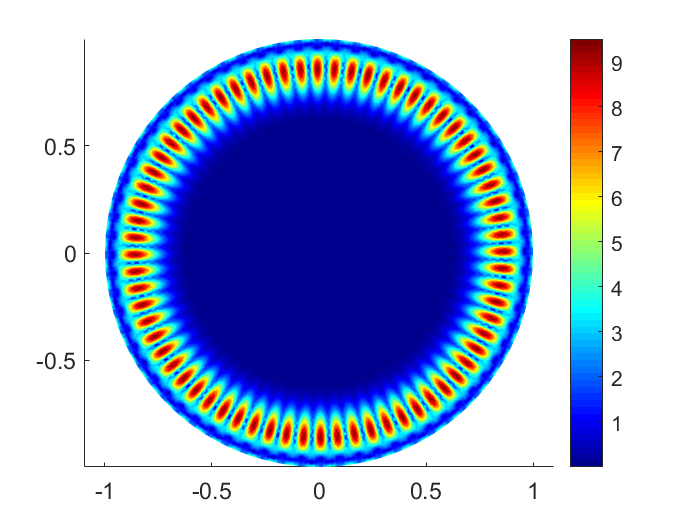}\\
\caption{The eigenfunctions $|\mathbf{v}_m|$ for bi-localized modes with $\omega_m=2.19$, $3.21$, $4.46$, $5.43$, $7.58$, and $9.48$.}
 \label{fig:radial bi v}
\end{figure}

\begin{figure}[htbp]
  \centering
  \includegraphics[width=0.32\textwidth]{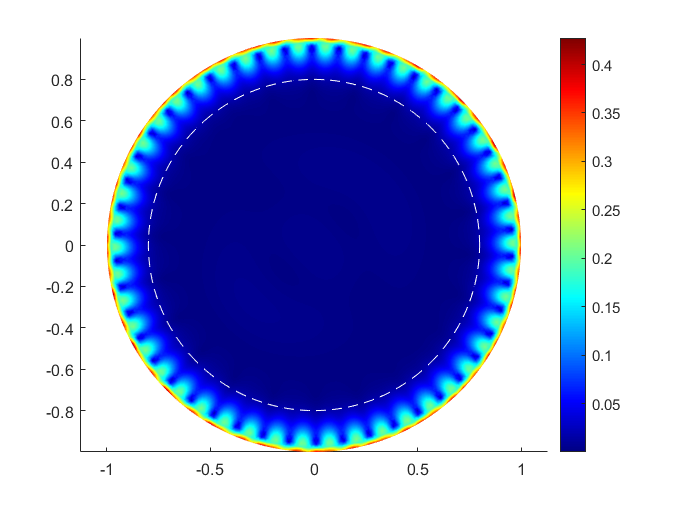}
  \includegraphics[width=0.32\textwidth]{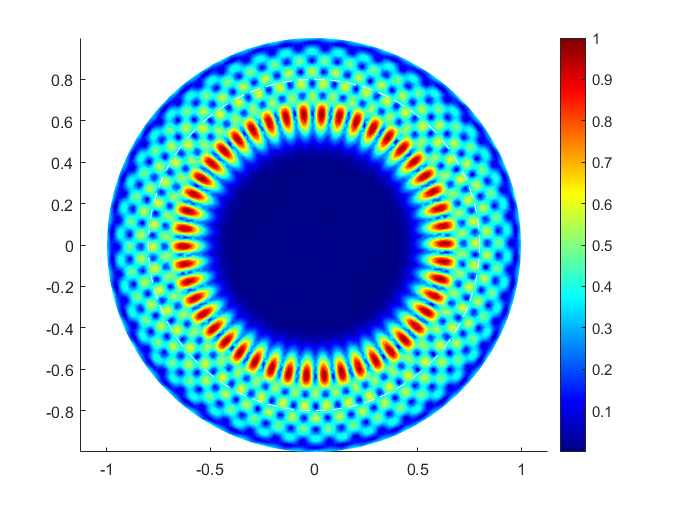}
\caption{The eigenfunctions $|\mathbf{u}_m|$ (left) and $|\mathbf{v}_m|$ (right) for piecewise constant $\widetilde{\rho}$ with $\omega_m=3.99$.}
 \label{fig:radial piecewise constant}
\end{figure}

\begin{figure}[htbp]
  \centering
  \includegraphics[width=0.32\textwidth]{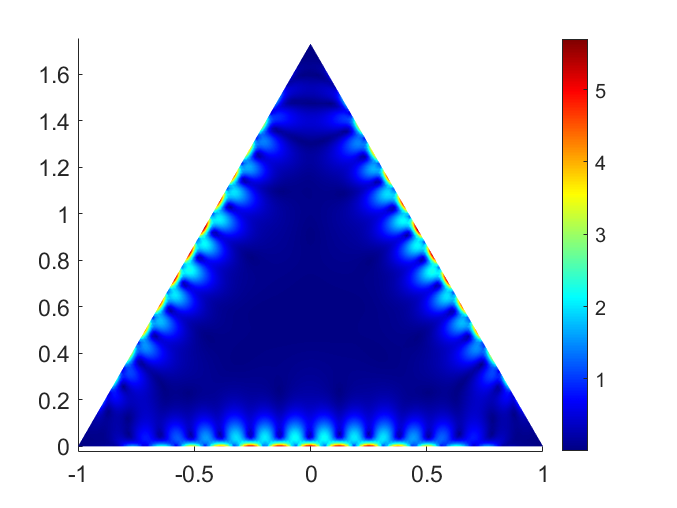}
  \includegraphics[width=0.32\textwidth]{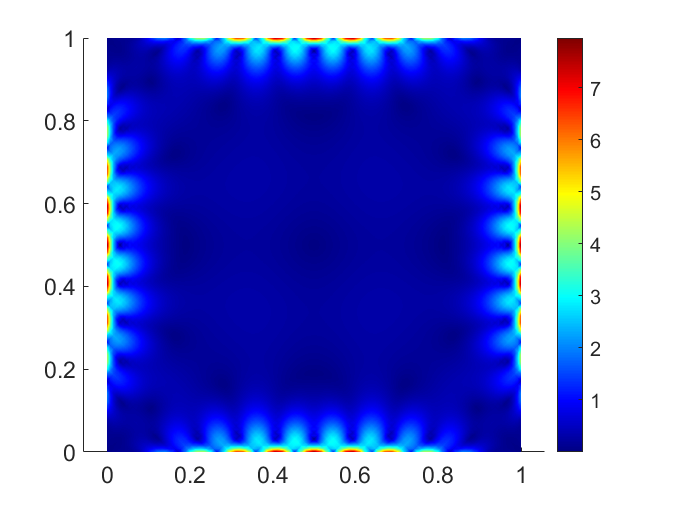}
  \includegraphics[width=0.32\textwidth]{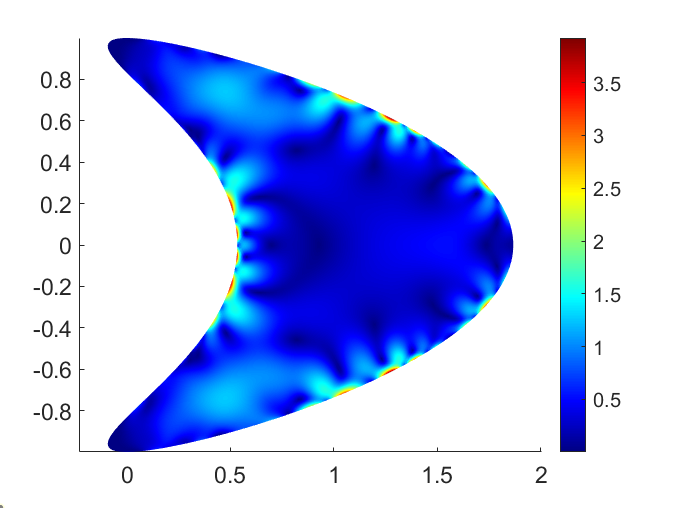}\\
  \includegraphics[width=0.32\textwidth]{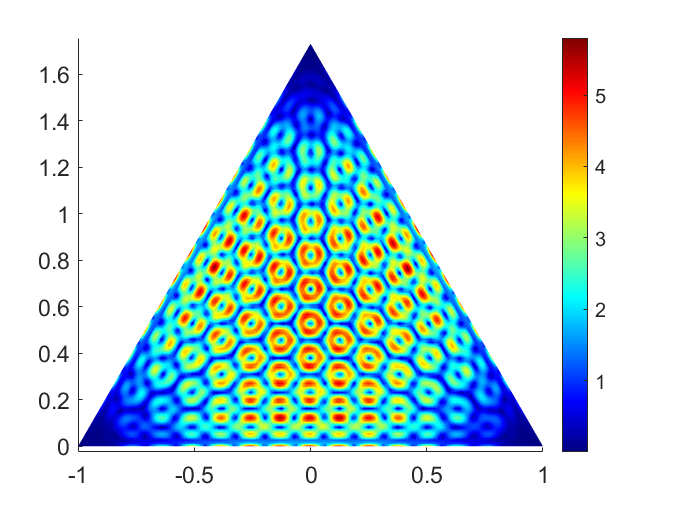}
  \includegraphics[width=0.32\textwidth]{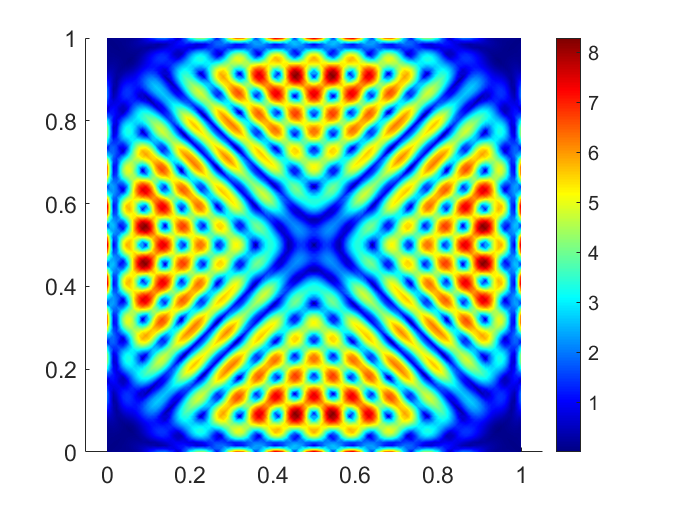}
  \includegraphics[width=0.32\textwidth]{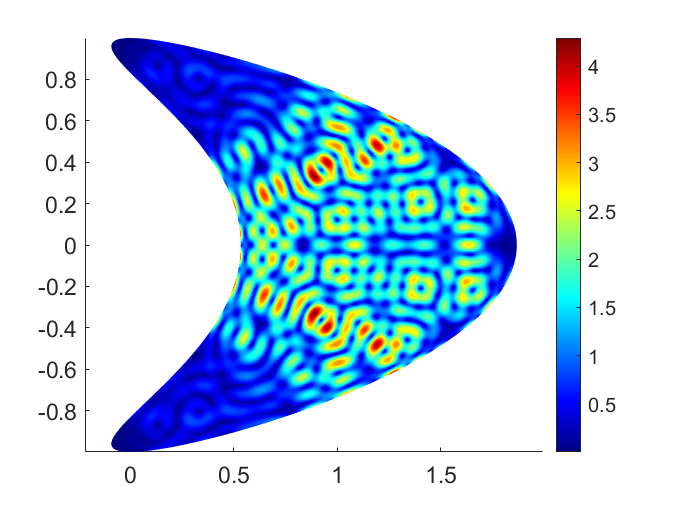}\\
\caption{The eigenfunctions $|\mathbf{u}|$ (up) and $|\mathbf{v}|$ (down) in triangle ($\omega_m=4.93$), square ($\omega_m=4.91$), and kite ($\omega_m=3.90$).}
 \label{fig:radial mono general geometry}
\end{figure}

\begin{figure}[htbp]
  \centering
  \includegraphics[width=0.48\textwidth]{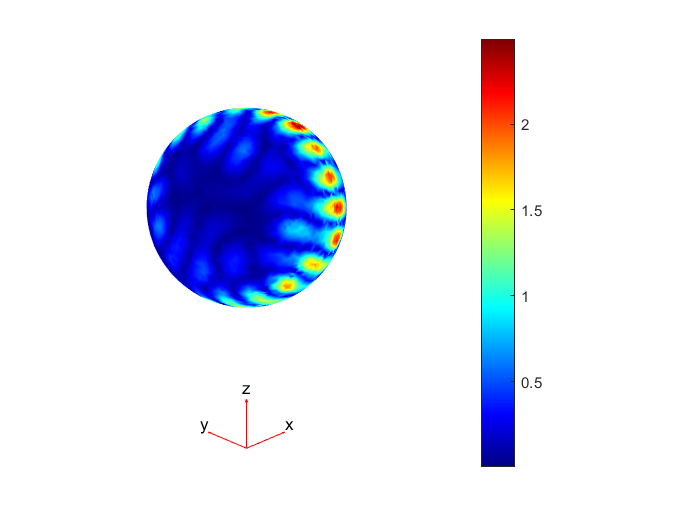}
  \includegraphics[width=0.48\textwidth]{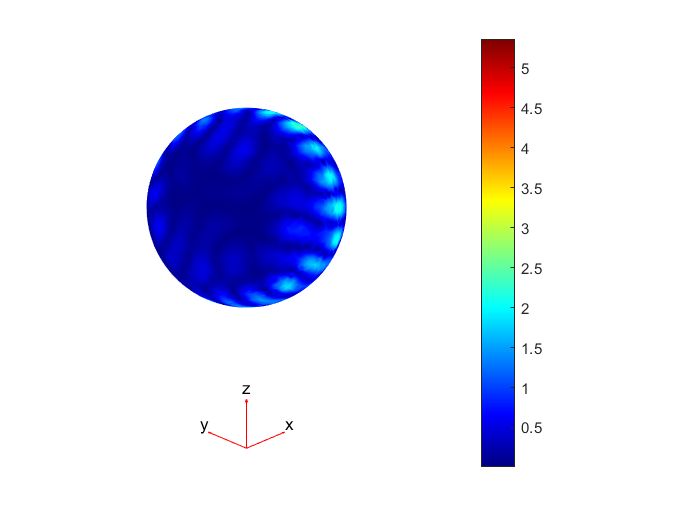}\\
  \includegraphics[width=0.48\textwidth]{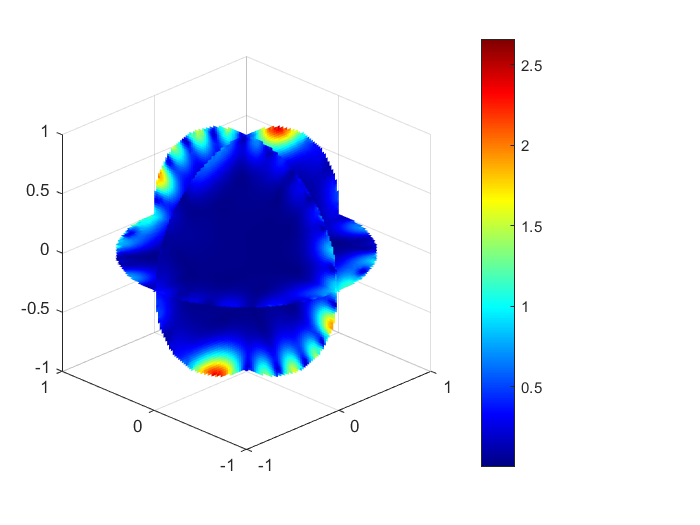}
  \includegraphics[width=0.48\textwidth]{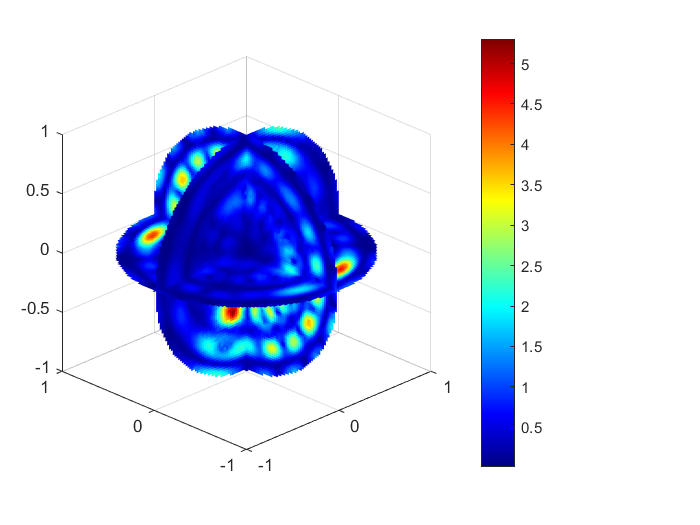}
\caption{The eigenfunctions $|\mathbf{u}|$ (left) and $|\mathbf{v}|$ (right) for mono-localized modes with $\omega_m=5.51$ in three dimensional ball.}
 \label{fig:three dimensional radial mono}
\end{figure}

\begin{figure}[htbp]
  \centering
  \includegraphics[width=0.45\textwidth]{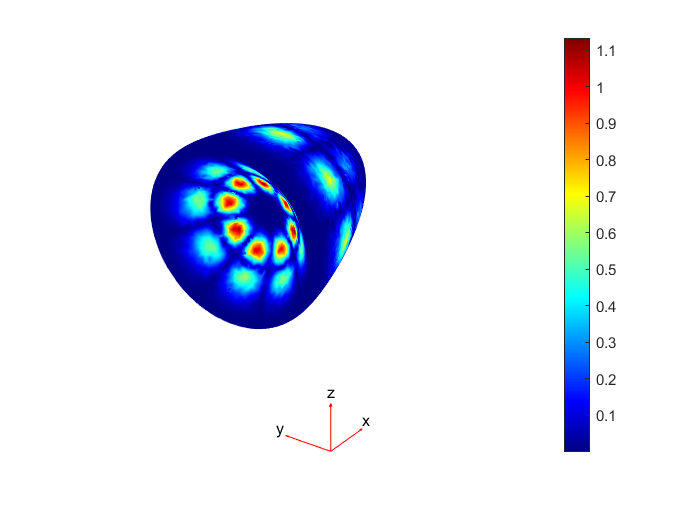}
  \includegraphics[width=0.45\textwidth]{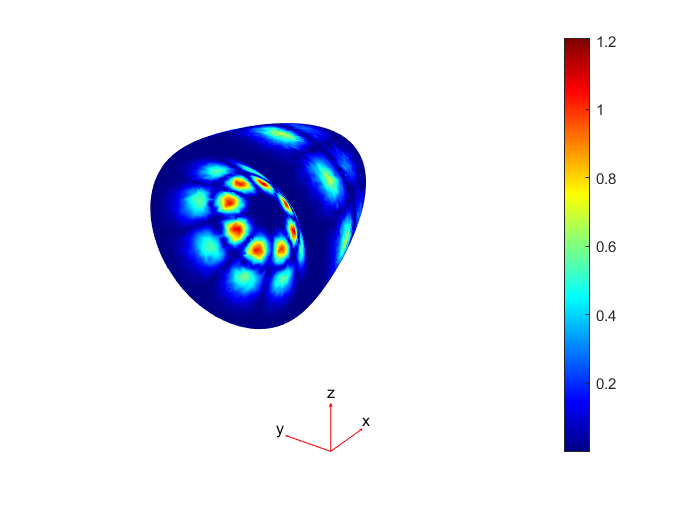}\\
  \includegraphics[width=0.45\textwidth]{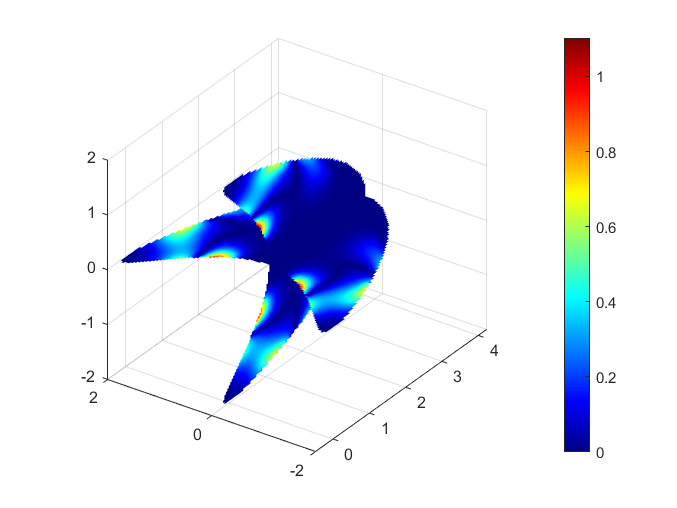}
  \includegraphics[width=0.45\textwidth]{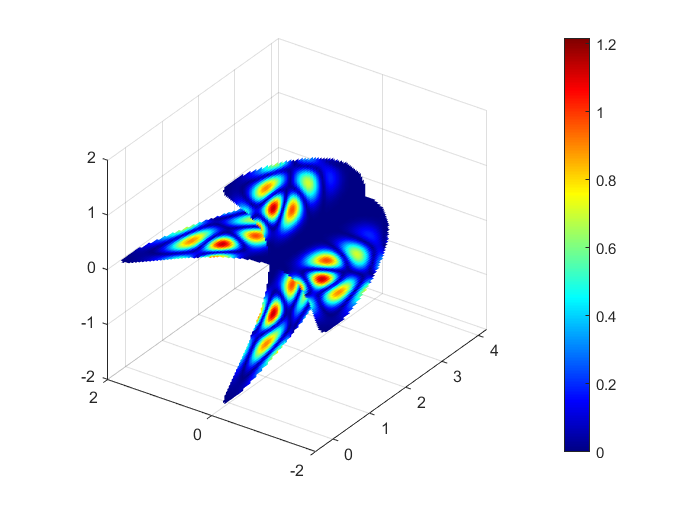}
\caption{The eigenfunctions $|\mathbf{u}|$ (left) and $|\mathbf{v}|$ (right) for mono-localized modes with $\omega_m=3.98$ in three dimensional kite.}
 \label{fig:three dimensional kite mono}
\end{figure}

\begin{figure}[htbp]
  \centering
  \includegraphics[width=0.45\textwidth]{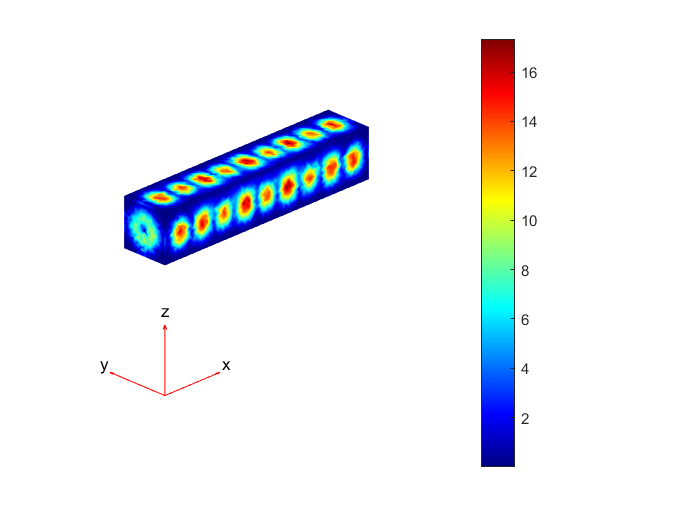}
  \includegraphics[width=0.45\textwidth]{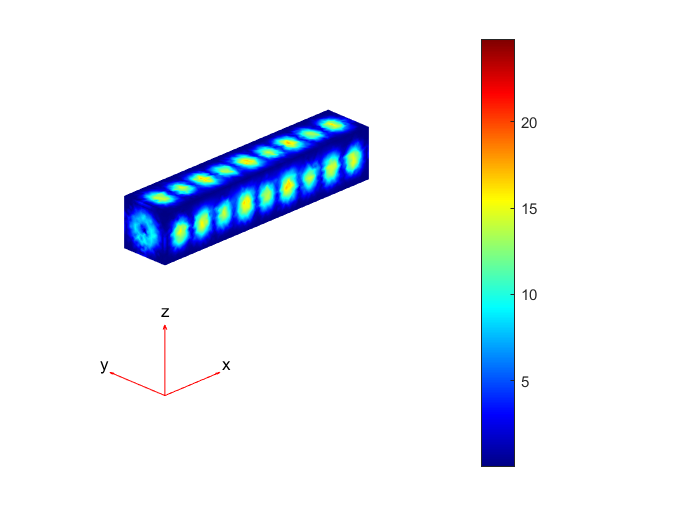}\\
  \includegraphics[width=0.45\textwidth]{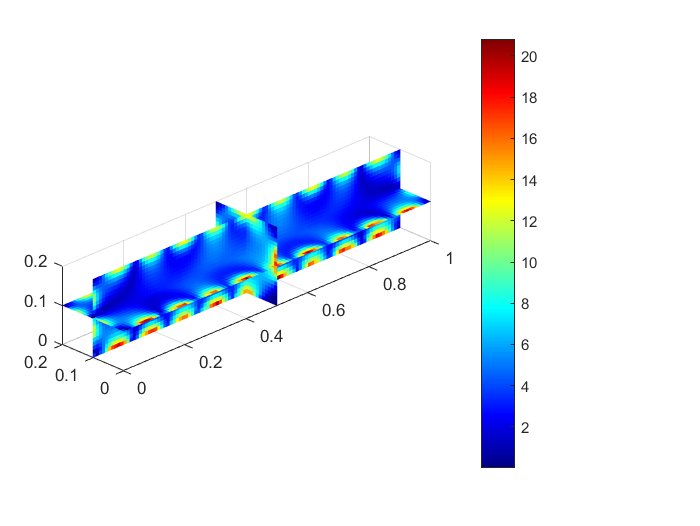}
  \includegraphics[width=0.45\textwidth]{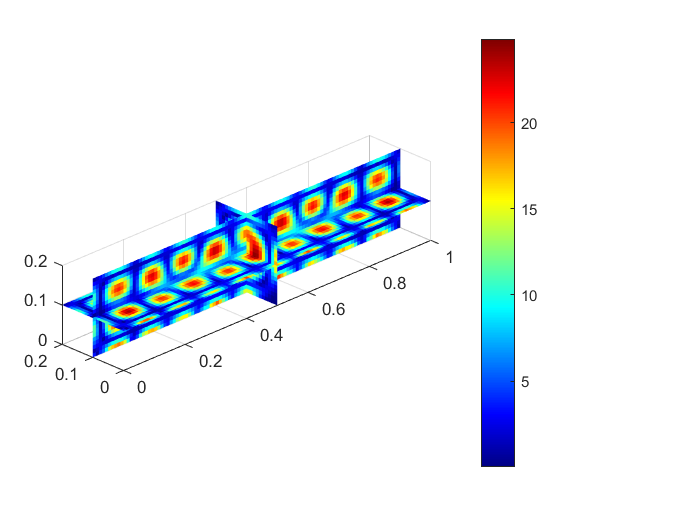}
\caption{The eigenfunctions $|\mathbf{u}|$ (left) and $|\mathbf{v}|$ (right) for mono-localized modes with $\omega_m=12.28$ in three dimensional cuboid.}
 \label{fig:three dimensional cube mono}
\end{figure}

First, we show that the boundary-localising properties hold in general geometric and parameter setups. The parameter $\rho$ and the Lam\'e parameter $\lambda$ is set to $1$ for all the cases. For the circle, Fig~\ref{fig:radial mono u} and \ref{fig:radial mono v} show the mono-localized modes and Fig~\ref{fig:radial bi u} and  \ref{fig:radial bi v} show the bi-localized modes with $\widetilde{\rho}=20$ and $\mu=1$, and Fig~\ref{fig:radial piecewise constant} shows the eigenfunctions with a piecewise constant $\widetilde{\rho}$ satisfying $\widetilde{\rho}=2*(r<0.8)+20*(r\geq 0.8)$ and $\mu=0.2$. For the other 2D cases, the eigenfunctions can be found in Fig~\ref{fig:radial mono general geometry} with $\widetilde{\rho}=20$ and $\mu=0.2$. Fig \ref{fig:three dimensional radial mono}, \ref{fig:three dimensional kite mono}, and \ref{fig:three dimensional cube mono} show the eigenfunctions for 3D cases with $\widetilde{\rho}=20$ and $\mu=1$.
On the other hand, we also note the surface-resonance phenomena of those boundary-localised modes, which are evidenced by the more severely oscillatory pattern along $\partial\Omega$ (compared to the ``natural" frequency $\omega_m$).



\begin{figure}[htbp]
  \centering
  \includegraphics[width=0.35\textwidth]{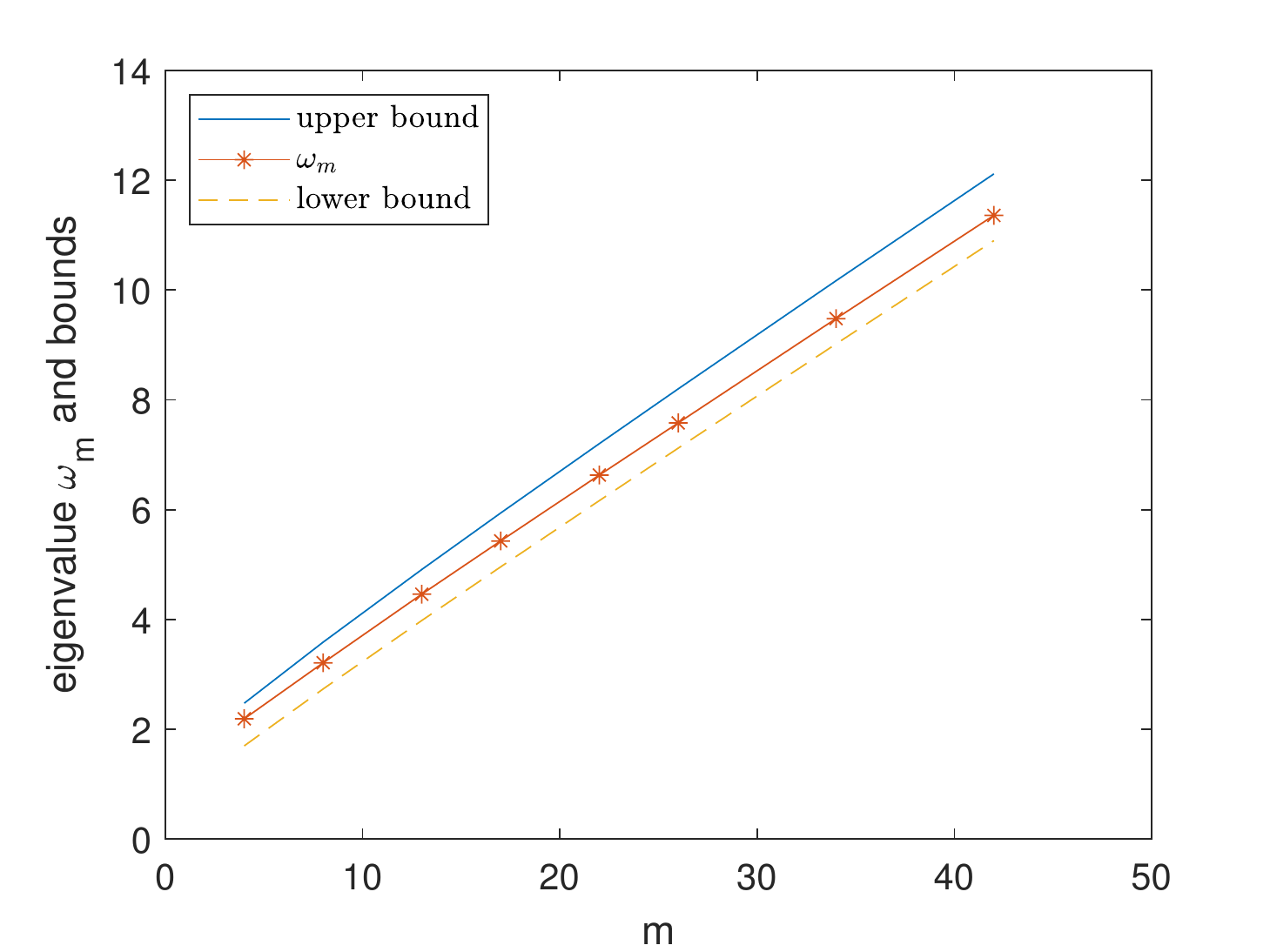}
  \includegraphics[width=0.35\textwidth]{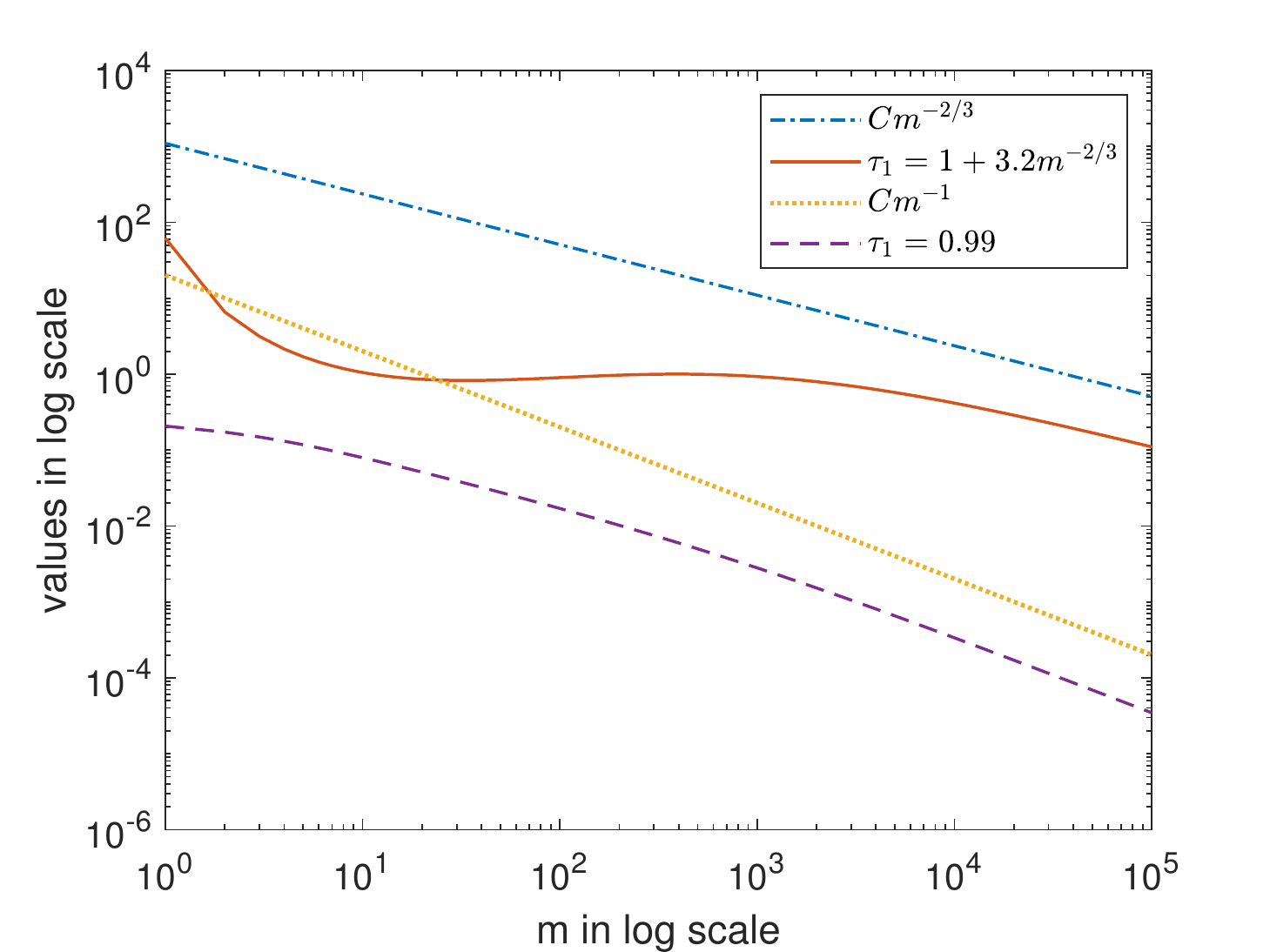}\\
  (a)\hspace{0.3\textwidth}(b)\\
  \includegraphics[width=0.35\textwidth]{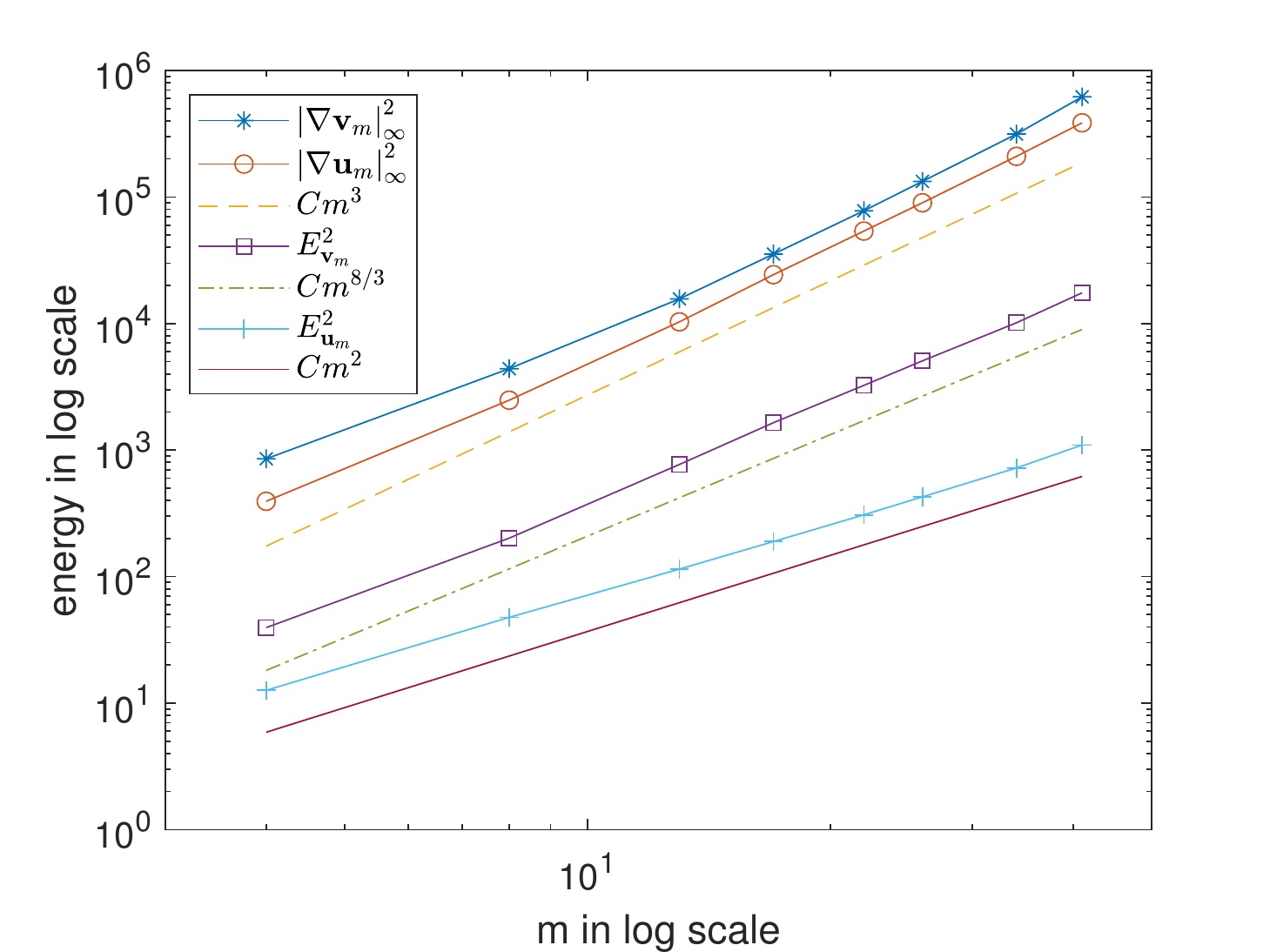}
  \includegraphics[width=0.35\textwidth]{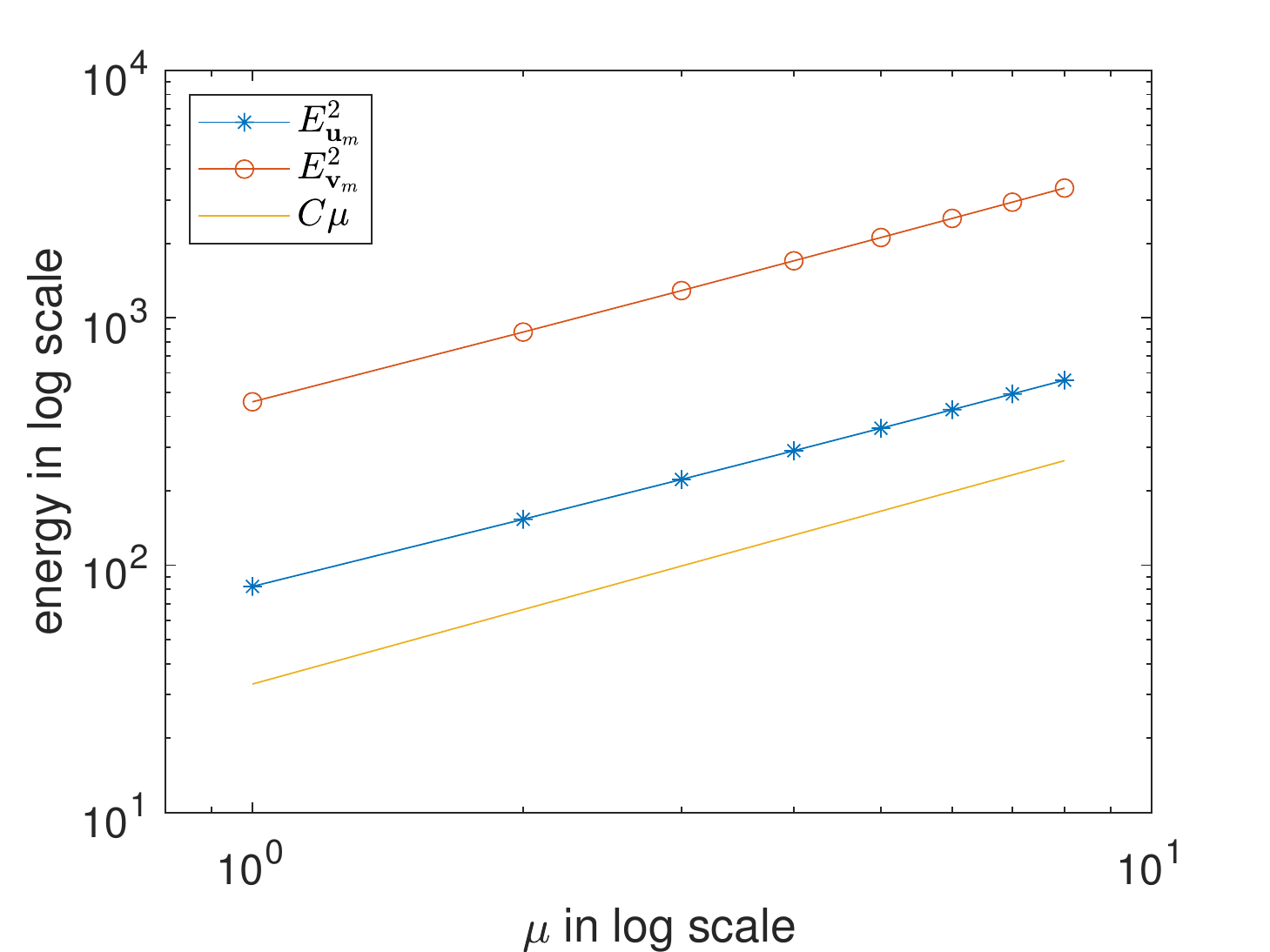}\\
  (c)\hspace{0.3\textwidth}(d)
\caption{(a) The upper and lower bounds of eigenvalue $\omega_m$ in Theorem~\ref{bilocalized};
(b) The convergence order of (\ref{eq:integral estimate}) with respect to $m$;
(c)-(d) The convergence of energy  with $m$ and $\mu$ for surface resonant waves for bi-localized modes.}
 \label{fig:bounds and convergence}
\end{figure}

Furthermore, we numerically verify {the theoretical results in Theorems~\ref{bilocalized} and} \ref{surface_resonant} of the bi-localized modes with $\rho=1$ and $\widetilde{\rho}=20$ in two dimensions. Consider the subdomain
$\Sigma(\tau, \theta_1, \theta_2)$ in (\ref{eq:theta12}) with {$\tau=2/3$}, $\theta_1=0$, and $\theta_2=\pi/3$. Talbe~\ref{tab:Surface resonant waves} shows the values of $E_{\mathbf{u}_{m}}^{2}$, $E_{\mathbf{v}_{m}}^{2}$, $\left|\nabla \mathbf{u}_m \right|^2_{\infty}$, and $\left|\nabla \mathbf{v}_m \right|^2_{\infty}$ in subdomain with different $m$ {and $\mu$, respectively}.
These quantities are increasing as $m$ {and $\mu$ increasing, respectively}. {In Fig~\ref{fig:bounds and convergence}(a), the upper and lower bounds mentioned in Theorem~\ref{bilocalized} of $\omega_m$ are shown, which is consistent with the theoretical results.} Fig~\ref{fig:bounds and convergence}(c) further shows that, {the increasing orders of $E_{\mathbf{u}_{m}}^{2}$, $E_{\mathbf{v}_{m}}^{2}, \left|\nabla \mathbf{u}_m \right|^2_{\infty}$, and $\left|\nabla \mathbf{v}_m \right|^2_{\infty}$ are $3, 10/3, 3$, and $10/3$ with respect to $m$, respectively. With respect to $m$, there is a little difference between numerical order and theoretical order. The reason is that the order of (\ref{eq:integral estimate}) is established when $m$ is large enough, which can be found in Fig~\ref{fig:bounds and convergence}(b) for $\tau_1=0.99$ and $\tau_1=1+3.2m^{-2/3}$. However, due to limiting computing resources, we only show the order of $m<42$ numerically. In addition, the orders of $E_{\mathbf{u}_{m}}^{2}$ and $E_{\mathbf{v}_{m}}^{2}$ are $1$ with respect to $\mu$ shown in Fig~\ref{fig:bounds and convergence}(d), which are consistent with Theorem~\ref{surface_resonant}}.

\begin{table}[htbp]
\caption{{Energy varies with $m$ (left) and $\mu$ (right) for surface resonant waves}.} \label{tab:Surface resonant waves}
\begin{tabular}{|r|r|r|r|r|r||r|r|r|r|}
    \hline
    \multicolumn{6}{|c||}{$\mu=\lambda=1$}&\multicolumn{4}{c|}{$\lambda=1,m=11$}\\
    \hline
     $m$ &  $\omega_m$      & $E_{\mathbf{u}_{m}}^{2}$                           &  $E_{\mathbf{v}_{m}}^{2}$
                 & $\left|\nabla \mathbf{u}_m \right|^2_{\infty}$     &  $\left|\nabla \mathbf{v}_m \right|^2_{\infty}$&
     $\mu$& $\omega_m$        & $E_{\mathbf{u}_{m}}^{2}$                           &  $E_{\mathbf{v}_{m}}^{2}$\\
    \hline   \hline
      $4$  &  $2.19$    &   $12.61$   &  $39.48$      & $393.59$    &$852.33$&
      $1$  &  $3.96$    &   $82.13$   &  $457.90 $\\
      $8$  & $3.21$    &   $47.62$   &  $201.25$   & $2477.95$   &$4398.70$&
      $2$   & $5.58$   &   $153.45$   &  $876.25$\\
      $13$  & $4.46$ &   $114.51$ &  $769.26$     & $10307.16$  & $15662.67$&
      $3$   & $6.82$   &   $222.35$   &  $1289.15$\\
      $17$  & $5.43$  &   $189.20$   &  $1644.74$  & $24317.75$  & $35378.37$ &
      $4$   & $7.87$   &   $290.41$   &  $1700.95$\\
      $22$  & $6.63$ &   $307.53$   &  $3243.68$  & $53785.40$ & $77824.26$&
      $5$   & $8.79$   &   $358.07$   &  $2112.50$\\
      $27$  & $7.58$ &   $426.62$   &  $5069.14$ & $89998.92$   & $132827.95$&
      $6$   & $9.62$   &   $425.53$   &  $2524.00$\\
      $34$  & $9.48$ &   $721.86$   &  $10135.35$ & $209526.63$   & $314377.42$&
      $7$   & $10.39$   &   $492.85$   &  $2935.51$\\
      $42$  & $11.36$ &   $1095.81$  &  $17481.54$ & $385699.50$  & $618602.54$&
      $8$   & $11.10$   &   $560.09$  &  $3347.03$\\
    \hline
\end{tabular}
\end{table}

\section*{Acknowledgment}
The work of H. Liu is supported by the Hong Kong RGC General Research Funds (projects 12302919, 12301420 and 11300821),  the NSFC/RGC Joint Research Fund (project N\_CityU101/21), the France-Hong Kong ANR/RGC Joint Research Grant, A-HKBU203/19.
The work of K. Zhang is supported in part by China Natural National Science Foundation (No.~12271207), and by the Key Laboratory of Symbolic Computation and Knowledge Engineering of Ministry of Education, Jilin University, China.
The work of J. Zhang is supported by the Natural Science Foundation of Jiangsu Province (No. BK20210540)
and The Natural Science Foundation of The Jiangsu Higher Education Institutions of China (No. 21KJB110015).


\begin{thebibliography}{99}


\bibitem{AS72}
{M.~Abramowitz and I. A.~Stegun}, {\it Handbook of Mathematical Functions: With Formulas,
Graphs, and Mathematical Tables},  US Department of Commerce,  Washington, DC, 1972.


\bibitem{B} E. Bl{\aa}sten, {\it Nonradiating sources and transmission eigenfunctions vanish at corners and edges}, SIAM J. Math. Anal., \textbf{50} (2018), no. 6, 6255--6270.

\bibitem{BL1} E. Bl{\aa}sten and H. Liu, {\it On vanishing near corners of transmission eigenfunctions}, J. Funct. Anal., \textbf{273} (2017), no. 11, 3616--3632. Addendum, arXiv: 1710.08089

\bibitem{BL} E. Bl{\aa}sten and H. Liu, {\it Scattering by curvatures, radiationless sources, transmission eigenfunctions, and inverse scattering problems}, SIAM J. Math. Anal., \textbf{53} (2021), no. 4, 3801--3837.

\bibitem{BLX} E. Bl{\aa}sten, H. Liu and J. Xiao, {\it On an electromagnetic problem in a corner and its applications}, Anal. PDE, \textbf{14} (2021), no. 7, 2207--2224.

\bibitem{CCH} F. Cakoni, D. Colton and H. Haddar, {\it Transmission eigenvalues}, Notices Amer. Math. Soc., \textbf{68} (2021), no. 9, 1499--1510.

\bibitem{CDHLW21} Y.-T. Chow, Y. Deng, Y. He, H. Liu and X. Wang, {\it Surface-localized transmission eigenstates,
super-resolution imaging and pseudo surface plasmon modes}, {SIAM J. Imaging Sci.}, \textbf{14} (2021), 946--975.

\bibitem{CDLS} Y.-T. Chow, Y. Deng, H. Liu and M. Sunkula, {\it Surface concentration of transmission eigenfunctions}, arXiv:2109.14361

\bibitem{DLWW22}
{Y. Deng, H. Liu, X. Wang and W. Wu}, {\it On Geometrical properties of electromagnetic transmission eigenfunctions and artificial mirage}, SIAM J. Appl. Math. \textbf{82} (2022), 1--24.

\bibitem{DJLZ22}
{Y. Deng, Y. Jiang, H. Liu and K. Zhang}, {\it On new surface-localized transmission eigenmodes},
Inverse Probl Imaging, \textbf{16} (2022), 595--611.


\bibitem{DCL}
H. Diao, X. Cao and H. Liu, {\it On the geometric structures of transmission eigenfunctions with a conductive boundary condition and applications}, Comm. Partial Differential Equations, \textbf{46} (2021), no. 4, 630--679.

\bibitem{DLLT}
H. Diao, H. Li, H. Liu and J. Tang, {\it Spectral properties of an acoustic-elastic transmission eigenvalue problem with applications}, arXiv:2210.16617

\bibitem{DLS}
H. Diao, H. Liu and B. Sun, {\it On a local geometric property of the generalized elastic transmission eigenfunctions and application}, Inverse Problems, \textbf{37} (2021), no. 10, Paper No. 105015, 36 pp.


\bibitem{JLZZ22}
{Y. Jiang, H. Liu, J. Zhang and K. Zhang}, {\it Boundary localization of transmission eigenfunctions in spherically stratified media},
Asymptotic Analysis, DOI: 10.3233/ASY-221794.


\bibitem{KRA06}
{I.~Krasikov}, {\it Uniform bounds for Bessel function}, Journal of Applied Analysis
\textbf{12} (2006), 83--91.

\bibitem{Liu}
H. Liu, {\it On local and global structures of transmission eigenfunctions and beyond}, J. Inverse Ill-Posed Probl., \textbf{30} (2022), no. 2, 287--305.

\bibitem{KOR02}
{B.G.~Korenev}, {\it Bessel functions and their applications}, Integral Transforms Spec. Funct., \textbf{25} (2002), 272--282.

\bibitem{MBDL} Q. Meng, Z. Bai, H. Diao and H. Liu, {\it Effective medium theory for embedded obstacles in elasticity with applications to inverse problems}, SIAM J. Appl. Math., \textbf{82} (2022), no. 2, 720--749.

\bibitem{PAR84}
{R.~Paris}, {\it An Inequality for the Bessel Function $J_{\nu} (\nu x)$}, SIAM J. Math. Anal., \textbf{16} (1984), 203--205.

\bibitem{WQ99}
{R.~Wong and C.~Qu}, {\it Best possible upper and lower bounds
for the zeros of the Bessel function $J_{\nu}$(x)}, Trans. Amer. Math. Soc., \textbf{351} (1999), 2833--2859.

\end{thebibliography}
\end{document}